\documentclass[a4paper, 12pt]{article}
\pdfoutput=1

\usepackage[margin=25mm]{geometry}

\usepackage[utf8]{inputenc}
\usepackage[T1]{fontenc}
\usepackage[british]{babel}
\usepackage[style=british]{csquotes}
\usepackage[en-GB]{datetime2}
\DTMlangsetup[en-GB]{ord=omit}
\usepackage{xcolor}

\usepackage{amsthm, thmtools}
\usepackage{mathtools}
\usepackage{titlesec} 
\usepackage{hyperref}
\hypersetup{
    bookmarksnumbered,
    colorlinks,
    linkcolor={cyan!50!blue},
    citecolor={cyan!50!blue},
    urlcolor={cyan!50!blue},
    hypertexnames=false 
}
\usepackage[capitalise,nameinlink]{cleveref}

\usepackage{tikz}
\usepackage{tikz-cd}

\usepackage[expansion=false,nopatch=footnote]{microtype}
\usepackage[lining,tabular,scale=.9]{FiraSans}
\usepackage[tt=false,semibold]{libertinus}
\usepackage[libertine,timesmathacc,smallerops]{newtxmath}
\usepackage[lining,scaled=.8]{FiraMono}
\usepackage[cal=boondox,frak=euler]{mathalpha}
\providecommand{\libertinusDisplay}{\LibertinusDisplay}

\DeclareFontEncoding{MDA}{}{}
\DeclareSymbolFont{mathdesignA}{MDA}{mdput}{m}{n}
\SetSymbolFont{mathdesignA}{bold}{MDA}{mdput}{b}{n}
\DeclareSymbolFontAlphabet{\mathbb}{mathdesignA}

\tikzcdset{
    arrow style=tikz,
    arrows={/tikz/line width=.5pt},
    diagrams={>={Straight Barb[scale=0.8]}}
}

\DeclareFontFamily{OMX}{MnSymbolE}{}
\DeclareSymbolFont{MnLargeSymbols}{OMX}{MnSymbolE}{m}{n}
\SetSymbolFont{MnLargeSymbols}{bold}{OMX}{MnSymbolE}{b}{n}
\DeclareFontShape{OMX}{MnSymbolE}{m}{n}{
    <-6>  MnSymbolE5
   <6-7>  MnSymbolE6
   <7-8>  MnSymbolE7
   <8-9>  MnSymbolE8
   <9-10> MnSymbolE9
  <10-12> MnSymbolE10
  <12->   MnSymbolE12
}{}
\DeclareFontShape{OMX}{MnSymbolE}{b}{n}{
    <-6>  MnSymbolE-Bold5
   <6-7>  MnSymbolE-Bold6
   <7-8>  MnSymbolE-Bold7
   <8-9>  MnSymbolE-Bold8
   <9-10> MnSymbolE-Bold9
  <10-12> MnSymbolE-Bold10
  <12->   MnSymbolE-Bold12
}{}
\DeclareMathDelimiter{[}{\mathopen}{MnLargeSymbols}{'000}{MnLargeSymbols}{'000}
\DeclareMathDelimiter{]}{\mathclose}{MnLargeSymbols}{'005}{MnLargeSymbols}{'005}
\DeclareMathDelimiter{\llbr}{\mathopen}{MnLargeSymbols}{'102}{MnLargeSymbols}{'102}
\DeclareMathDelimiter{\rrbr}{\mathclose}{MnLargeSymbols}{'107}{MnLargeSymbols}{'107}

\usepackage{bm}
\usepackage{cases}
\usepackage{accents}

\usepackage{setspace}
\usepackage{subdepth} 

\newcommand{\initlengths}{%
    \setlength{\abovedisplayshortskip}{3pt plus 9pt minus 3pt}%
    \setlength{\belowdisplayshortskip}{9pt plus 9pt minus 9pt}%
    \setlength{\abovedisplayskip}{9pt plus 9pt minus 9pt}%
    \setlength{\belowdisplayskip}{9pt plus 9pt minus 9pt}%
    \hfuzz 1pt%
    \tolerance 400
}

\usepackage{titling}
\pretitle{\vspace{-50pt}\setstretch{1.05}\begin{center}\LARGE\libertinusDisplay\fontdimen2\font=0.3em} 
\posttitle{\end{center}\vspace{-3pt}}

\numberwithin{paragraph}{subsection}

\newcommand{\parasep}{9pt plus 3pt minus 3pt}
\titleformat{\section}{\Large\libertinusDisplay}{\thesection}{1em}{}
\titleformat{\subsection}{\large\firamedium\boldmath}{\thesubsection}{1em}{}

\usepackage{tocloft}

\renewenvironment{abstract}{%
    \setstretch{1.15}%
    \centering\begin{minipage}{.85\textwidth}%
    \setlength{\parindent}{1.5em}%
    \centerline{\firamedium\abstractname}%
    \par\vspace{4.5pt}%
}{\end{minipage}\par\vspace{6pt}}

\makeatletter
\newcommand*{\@parabookmark}{%
  \pdfbookmark[3]{%
    \theparagraph
    \ifx\@currentlabelname\@empty
    \else
      .\space\@currentlabelname%
    \fi
  }{\theparagraph}
}
\newcommand*{\@thmbookmark}{%
  \pdfbookmark[3]{%
    \theparagraph.\space\thmt@thmname
    \ifx\@currentlabelname\@empty
    \else
      .\space\@currentlabelname%
    \fi
  }{\theparagraph}
}
\newcommand*{\parabookmark}{\@parabookmark}
\newcommand*{\thmbookmark}{\@thmbookmark}
\newcommand*{\suppressparabookmarks}{%
  \renewcommand*{\parabookmark}{}%
  \renewcommand*{\thmbookmark}{}%
}
\newcommand*{\resumeparabookmarks}{%
  \renewcommand*{\parabookmark}{\@parabookmark}%
  \renewcommand*{\thmbookmark}{\@thmbookmark}%
}

\declaretheoremstyle[
    spaceabove=\parasep, spacebelow=\parasep,
    postheadspace=.5em,
    postheadhook=\thmbookmark,
    headfont=\normalfont\bfseries,
    headpunct={},
    headformat={\NUMBER.\@\ \NAME.\@\NOTE},
    notefont=\normalfont\bfseries\boldmath,
    notebraces={}{.},
    bodyfont=\itshape,
]{theorem}
\declaretheoremstyle[
    spaceabove=\parasep, spacebelow=\parasep,
    postheadspace=.5em,
    headfont=\normalfont\bfseries,
    headpunct={},
    headformat={\NAME.\@\NOTE},
    notefont=\normalfont\bfseries\boldmath,
    notebraces={}{.},
    bodyfont=\itshape,
]{theorem*}
\declaretheoremstyle[
    spaceabove=\parasep, spacebelow=\parasep,
    postheadspace=.5em,
    postheadhook=\thmbookmark,
    headfont=\normalfont\bfseries,
    headpunct={},
    headformat={\NUMBER.\@\ \NAME.\@\NOTE},
    notefont=\normalfont\bfseries\boldmath,
    notebraces={}{.},
]{definition}
\declaretheoremstyle[
    spaceabove=\parasep, spacebelow=\parasep,
    postheadspace=.5em,
    postheadhook=\parabookmark,
    headfont=\normalfont\bfseries,
    headpunct={},
    headformat={\NUMBER.\@\NOTE},
    notefont=\normalfont\bfseries\boldmath,
    notebraces={}{.},
]{para}

\renewenvironment{proof}[1][\proofname]{\par
    \pushQED{\qed}%
    \normalfont\trivlist
    \item[\hskip\labelsep\bfseries #1\@addpunct{.}]\ignorespaces
}{%
    \popQED\endtrivlist\@endpefalse
}

\expandafter\def\csname equation*@qed\endcsname{\equation@qed}
\makeatother

\declaretheorem[sibling=paragraph, style=para, refname={\S,\S\S}]{para}
\declaretheorem[sibling=paragraph, style=theorem, name=Theorem]{theorem}
\declaretheorem[sibling=paragraph, style=theorem, name=Lemma]{lemma}

\declaretheorem[numbered=no, style=theorem*, name=Theorem]{theorem*}
\declaretheorem[numbered=no, style=theorem*, name=Lemma]{lemma*}

\declaretheorem[sibling=paragraph, style=definition, name=Example]{example}

\numberwithin{equation}{paragraph}

\crefdefaultlabelformat{#2{\upshape#1}#3}

\crefformat{equation}{#2{\upshape(#1)}#3}
\crefrangeformat{equation}{{\upshape#3(#1)#4--#5(#2)#6}}
\crefmultiformat{equation}{#2{\upshape(#1)}#3}{ and #2{\upshape(#1)}#3}{, #2{\upshape(#1)}#3}{, and~#2{\upshape(#1)}#3}

\crefformat{enumi}{#2{\upshape#1}#3}
\crefrangeformat{enumi}{{\upshape#3#1#4--#5#2#6}}
\crefmultiformat{enumi}{#2{\upshape#1}#3}{ and #2{\upshape#1}#3}{, #2{\upshape#1}#3}{, and~#2{\upshape#1}#3}

\crefformat{section}{#2{\upshape\S#1}#3}
\crefrangeformat{section}{{\upshape#3\S\S#1#4--#5#2#6}}
\crefmultiformat{section}{{\upshape#2\S\S#1#3}}{ and {\upshape#2#1#3}}{, {\upshape#2#1#3}}{, and~{\upshape#2#1#3}}
\crefformat{subsection}{#2{\upshape\S#1}#3}
\crefrangeformat{subsection}{{\upshape#3\S\S#1#4--#5#2#6}}
\crefmultiformat{subsection}{{\upshape#2\S\S#1#3}}{ and {\upshape#2#1#3}}{, {\upshape#2#1#3}}{, and~{\upshape#2#1#3}}
\crefformat{para}{#2{\upshape\S#1}#3}
\crefrangeformat{para}{{\upshape#3\S\S#1#4--#5#2#6}}
\crefmultiformat{para}{{\upshape#2\S\S#1#3}}{ and {\upshape#2#1#3}}{, {\upshape#2#1#3}}{, and~{\upshape#2#1#3}}

\crefname{figure}{Figure}{Figures}

\usepackage{enumitem}
\setlist{noitemsep}
\setlist[enumerate]{label=\textnormal{(\roman*)}}

\usepackage[labelsep=period, labelfont={bf}]{caption}

\usepackage[
    backend=biber,
    style=oxnum,
    giveninits,
    maxbibnames=5,
    maxcitenames=5,
    sorting=nyt
]{biblatex}
\DeclareFieldFormat*{title}{\textit{#1}}
\DeclareFieldFormat*{journaltitle}{#1}
\DeclareFieldFormat[article]{version}{\IfDecimal{#1}{version~}{}#1}
\DefineBibliographyStrings{english}{
  withafterword = {with an appendix by},
}

\DeclareFieldFormat{pages}{%
  \iffieldundef{bookpagination}%
    {#1}%
    {\mkpageprefix[bookpagination]{#1}}%
}

\newcommand{\calBun}{{\mathcal{B}\mkern-2mu\mathit{un}}}
\newcommand{\calHom}{{\mathcal{H}\mkern-3mu\mathit{om}}}
\newcommand{\calMap}{{\mathcal{M}\mkern-3mu\mathit{ap}\mkern1mu}}
\newcommand{\calPerf}{\mathcal{P}\mkern-2mu\mathit{er}\mkern-1mu\mathit{f}}
\newcommand{\heart}{\mathbin{\heartsuit}}
\newcommand{\llparen}{\mathopen{\mathrlap{(}\mspace{4mu}(}}
\newcommand{\longhookrightarrow}{\lhook\joinrel\longrightarrow}
\newcommand{\longsimto}{\mathrel{\overset{\smash{\raisebox{-.8ex}{$\sim$}}\mspace{3mu}}{\longrightarrow}}}
\newcommand{\lowersup}[1]{^{\raisebox{-.3ex}{$\scriptscriptstyle#1$}}}
\newcommand{\lowerSup}[1]{^{\raisebox{-.3ex}{$\scriptstyle#1$}}}
\newcommand{\rrparen}{\mathclose{)\mspace{4mu}\mathllap{)}}}
\newcommand{\simto}{\mathrel{\overset{\smash{\raisebox{-.8ex}{$\sim$}}\mspace{3mu}}{\to}}}

\newcommand{\sqrte}{\mathord{\smash{\undisplay\sqrt{e}}\mathstrut}}
\newcommand{\undisplay}{\mathchoice{\textstyle}{}{}{}}
\newcommand{\vdim}{\operatorname{vdim}}

\usepackage{transparent}
\newcommand{\bbnu}{\mathord{\mathrlap{\transparent{0}\nu}\mathchoice{\bbnutikz}{\bbnutikz}{\raisebox{-.05ex}{\scalebox{.75}{\bbnutikz}}}{}}}
\newcommand{\bbnutikz}{\tikz{%
    \draw[line width=.1ex, line cap=round, line join=round]
        (0, .9ex)
        .. controls (.2ex, .9ex) .. (.4ex, 1ex)
        -- (.64ex, .25ex)
        .. controls (.8ex, .4ex) and (1.05ex, .7ex) .. (1ex, .85ex)
        .. controls (.95ex, 1ex) and (1.06ex, 1.05ex) .. (1.1ex, .95ex)
        .. controls (1.12ex, .9ex) and (1ex, .6ex) .. (.42ex, 0ex)
        -- (.15ex, .9ex) -- cycle;
    \fill (.64ex, .25ex)
        .. controls (.8ex, .4ex) and (1.05ex, .7ex) .. (1ex, .85ex)
        .. controls (.95ex, 1ex) and (1.06ex, 1.05ex) .. (1.1ex, .95ex)
        .. controls (1.12ex, .9ex) and (1ex, .6ex) .. (.42ex, 0ex);
}}

\renewcommand{\geq}{\geqslant}
\renewcommand{\leq}{\leqslant}

\makeatletter
\def\big#1{{\hbox{$\left#1\vbox to10\p@{}\right.\n@space$}}}
\makeatother

\title{Modules and generalizations of Joyce vertex algebras}
\author{Chenjing Bu}
\date{}

\begin{document}

\initlengths

\maketitle

\begin{abstract}
    Joyce vertex algebras are vertex algebra structures
defined on the homology of certain $\mathbb{C}$-linear moduli stacks,
and are used to express wall-crossing formulae
for Joyce's homological enumerative invariants.
This paper studies the generalization of this construction
to settings that come from non-linear enumerative problems.
In the special case of orthosymplectic enumerative geometry,
we obtain twisted modules for Joyce vertex algebras.

We expect that our construction will be useful for
formulating wall-crossing formulae for
enumerative invariants for non-linear moduli stacks.
We include several variants of our construction
that apply to different flavours of enumerative invariants,
including Joyce's homological invariants,
DT4 invariants,
and a version of $K$-theoretic enumerative invariants.

\end{abstract}

{
    \setstretch{1.15}
    \hypersetup{linkcolor=black}
    \tableofcontents
}

\clearpage
\suppressparabookmarks
\section{Introduction}

\addtocounter{subsection}{1}

\begin{para}
    \emph{Joyce vertex algebras}, introduced by
    \textcite{joyce-hall,joyce-homological,joyce-surfaces},
    are a class of vertex algebras
    defined on the homology of certain moduli stacks.
    They are motivated by the study of
    wall-crossing in enumerative geometry,
    and it has been mysterious whether they are related to
    usual sources of vertex algebras,
    such as conformal field theories.

    In some cases, this construction gives familiar vertex algebras,
    such as the Heisenberg vertex algebra,
    which arises as the homology of the space~$\mathrm{BU}$,
    or certain Kac--Moody vertex algebras,
    which arise as the homology of moduli stacks
    of representations of Dynkin quivers.
    See \cref{para-eg-perf,para-quiver}
    below for these examples,
    and \textcite{latyntsev-thesis} for more details.
\end{para}

\begin{para}
    Motivated by the problem of generalizing
    Joyce's \cite{joyce-homological}
    enumerative invariants and wall-crossing formulae
    from the setting of abelian categories
    to more general algebraic stacks,
    \textcite{bu-self-dual-ii}
    considered, as a first step,
    the orthosymplectic case,
    involving moduli stacks of orthogonal or symplectic objects
    and enumerative invariants counting them.

    In this case,
    \textcite{bu-self-dual-ii}
    constructed \emph{twisted modules} for Joyce vertex algebras,
    defined on the homology
    of moduli stacks of orthogonal and symplectic objects.
    These twisted module structures were used to write down
    wall-crossing formulae for orthosymplectic enumerative invariants.
    This structure was later also studied by
    \textcite{dehority-latyntsev}.
\end{para}

\begin{para}
    \label{para-intro-gen}
    More recently,
    the author, Halpern-Leistner, Ibáñez Núñez, and Kinjo
    \cite{epsilon-i,epsilon-ii,epsilon-iii}
    developed \emph{intrinsic Donaldson--Thomas theory},
    a new framework for enumerative geometry
    designed for generalizing results from abelian categories
    to non-linear moduli stacks.
    We expect this framework to help formulating
    a generalization of Joyce's
    \cite{joyce-homological}
    formalism to a larger class of
    quasi-smooth stacks over~$\mathbb{C}$,
    including stacks that do not come from linear categories.

    Although such a generalized enumerative theory
    is still elusive,
    in this paper, we write down a natural generalization of
    Joyce's vertex algebras for general stacks,
    which we call \emph{vertex induction},
    and which also generalizes the twisted modules in \cite{bu-self-dual-ii}.
\end{para}

\begin{para}
    We expect the vertex induction to be the key operation in
    \emph{wall-crossing formulae} for these conjectural invariants,
    which relate the invariants defined for different stability conditions.
    Moreover, these wall-crossing formulae
    should have the same structure as those for
    motivic Donaldson--Thomas invariants in \cite{epsilon-iii}.
    These expectations generalize
    Joyce's work \cite{joyce-homological} in the linear case.

    Wall-crossing formulae are a powerful tool
    in enumerative geometry,
    as they provide a strong constraint on the structure of the invariants,
    and can sometimes be used for direct computations
    which can otherwise be difficult.
    See
    \textcite{gross-joyce-tanaka-2022},
    \textcite{bu-2023-curves}, and
    \textcite{bojko-lim-moreira-2024}
    for examples of this approach in the linear case,
    where Joyce vertex algebras were an essential tool.
\end{para}

\begin{para}
    We hope that the vertex induction defined in the current work,
    combined with the formalism of
    stability and wall-crossing formulae for general stacks
    in \cite{epsilon-i,epsilon-ii,epsilon-iii},
    will enable us to formulate a precise conjecture
    for the existence and behaviour of the generalized invariants
    in \cref{para-intro-gen},
    and will allow us to compute them in some cases,
    assuming the conjectural wall-crossing formulae.

    Moreover, apart from Joyce's invariants,
    several other types of enumerative invariants
    also exhibit very similar wall-crossing behaviours,
    including \emph{DT4 invariants} and certain \emph{$K$-theoretic invariants},
    which we will discuss in
    \crefrange{para-intro-dt4}{para-intro-k-inv}.
    Their wall-crossing formulae should be expressed using
    their corresponding variants of the vertex induction,
    which we also introduce in this paper.
\end{para}

\begin{para}[Joyce vertex algebras]
    Before introducing the general construction of vertex induction,
    we will first explain its special cases where we obtain
    Joyce vertex algebras,
    and modules and twisted modules for them,
    to better motivate the general construction.

    Roughly speaking,
    a \emph{vertex algebra} is a vector space~$V$,
    equipped with a unit $1 \in V$
    and multiplication maps denoted by
    \begin{equation}
        (a_1, \dotsc, a_n) \longmapsto
        a_1 (z_1) \cdots a_n (z_n)
        \in V \llbr z_1, \dotsc, z_n \rrbr \, [(z_i - z_j)^{-1}] \ ,
    \end{equation}
    which we may regard as a family of multiplications
    depending meromorphically on the formal variables~$z_i$,
    with possible poles along $z_i = z_j$ for $i \neq j$.
    This multiplication should be
    unital, associative, and commutative in suitable senses.

    In \cref{sec-va}, we recall the construction of
    Joyce vertex algebras,
    with a precise statement in \cref{thm-joyce-va}.
    We start with a moduli stack~$\mathcal{X}$
    of objects in a $\mathbb{C}$-linear abelian category~$\mathcal{A}$.
    Then the graded vector space
    $\mathrm{H}_{\bullet + 2 \vdim} (\mathcal{X}; \mathbb{Q})$
    admits a vertex algebra structure given by
    \begin{equation}
        a_1 (z_1) \cdots a_n (z_n) =
        \oplus_* \circ \tau (z)
        \Bigl(
            (a_1 \boxtimes \cdots \boxtimes a_n) \cap
            e_z^{-1} (\bbnu)
        \Bigr) \ ,
    \end{equation}
    where $\vdim$ denotes the virtual dimension of~$\mathcal{X}$,
    $a_i \in \mathrm{H}_\bullet (\mathcal{X}; \mathbb{Q})$
    are homology classes,
    $\oplus \colon \mathcal{X}^n \to \mathcal{X}$
    is the direct sum map,
    $\tau (z)$ is a translation operator
    depending on the variables~$z_i$,
    and $e_z^{-1} (\bbnu)$ is the inverse
    of the equivariant Euler class
    of the (virtual) normal bundle of the map~$\oplus$.
    This $n$-fold multiplication will be a special case
    of vertex induction applied to the moduli stack~$\mathcal{X}$.
\end{para}

\begin{para}[Modules]
    \label{para-intro-vm}
    In \cref{sec-vm},
    we construct modules and twisted modules
    for Joyce vertex algebras.

    The motivating case is from orthosymplectic enumerative geometry,
    where we start with a moduli stack~$\mathcal{X}$
    of objects in an abelian category~$\mathcal{A}$,
    equipped with a contravariant involution
    $(-)^\vee \colon \mathcal{A} \simto \mathcal{A}^\mathrm{op}$,
    which induces a $\mathbb{Z}_2$-action on~$\mathcal{X}$.
    The homotopy fixed locus~$\mathcal{X}^{\mathbb{Z}_2}$
    is then the moduli of orthogonal or symplectic objects.
    See \textcite{bu-osp-dt} for more details on this set-up.

    In this case, the involution~$(-)^\vee$
    induces a \emph{twisted involution} of the Joyce vertex algebra
    $V = \mathrm{H}_{\bullet + 2 \vdim} (\mathcal{X}; \mathbb{Q})$,
    meaning that it satisfies
    \begin{equation}
        \bigl( a_1 (z_1) \cdots a_n (z_n) \bigr)^\vee =
        a_1^\vee (-z_1) \cdots a_n^\vee (-z_n)
    \end{equation}
    for $a_1, \dotsc, a_n \in V$. The space
    $M = \mathrm{H}_{\bullet + 2 \vdim} (\mathcal{X}^{\mathbb{Z}_2}; \mathbb{Q})$
    has the structure of a \emph{twisted module} for~$V$,
    given by the action
    \begin{equation}
        \label{eq-intro-vm-twisted}
        a_1 (z_1) \cdots a_n (z_n) \, m =
        \diamond_* \circ \tau (z)
        \Bigl(
            (a_1 \boxtimes \cdots \boxtimes a_n \boxtimes m) \cap
            e_z^{-1} (\bbnu)
        \Bigr) \ ,
    \end{equation}
    where $a_i \in V$ and $m \in M$ are homology classes,
    $\diamond \colon \mathcal{X}^n \times \mathcal{X}^{\mathbb{Z}_2} \to \mathcal{X}^{\mathbb{Z}_2}$
    is the map sending
    $(x_1, \dotsc, x_n, y) \mapsto
    x_1 \oplus x_1^\vee \oplus \cdots \oplus x_n \oplus x_n^\vee \oplus y$,
    and $\bbnu$ now denotes the normal bundle of~$\diamond$.

    Here, having a twisted module
    means that we have
    \begin{equation}
        a_1 (z_1) \cdots a_n (z_n) \, m \in
        M \llbr z_1, \dotsc, z_n \rrbr
        \, [z_i^{-1}, (z_i \pm z_j)^{-1}] \ ,
    \end{equation}
    where we allow extra poles at $z_i + z_j = 0$ for $i \neq j$,
    which are not present in a usual module,
    and we require the relation
    \begin{equation}
        a (z) \, m = a^\vee (-z) \, m
    \end{equation}
    for $a \in V$ and $m \in M$.
    This forces us to allow the extra poles,
    as it implies relations such as
    $a_1 (z_1) \cdots a_n (z_n) \, m =
    a_1^\vee (-z_1) \, a_2 (z_2) \cdots a_n (z_n) \, m$.
    Note that our notion of twisted modules is
    different from the one in \textcite[\S5.6]{frenkel-ben-zvi-2004}.
\end{para}

\begin{para}[Vertex induction]
    \label{para-intro-vi}
    The vertex induction
    is the main construction of this paper,
    and generalizes the constructions of
    Joyce vertex algebras and modules above.
    We will introduce it in \cref{sec-vi},
    with precise statements in
    \cref{thm-vi,thm-vi-stack}.

    The motivating situation is
    when we are given an algebraic stack~$\mathcal{X}$ over~$\mathbb{C}$,
    and we consider its \emph{stack of graded points}
    \begin{equation}
        \mathrm{Grad} (\mathcal{X})
        = \calMap ({*} / \mathbb{G}_\mathrm{m}, \mathcal{X}) \ ,
    \end{equation}
    following \textcite{halpern-leistner-instability}.
    It is equipped with a forgetful map
    $\mathrm{Grad} (\mathcal{X}) \to \mathcal{X}$.

    For example, if~$\mathcal{X}$ is the moduli stack
    of objects in an abelian category~$\mathcal{A}$,
    then~$\mathrm{Grad} (\mathcal{X})$ is usually the moduli stack
    of $\mathbb{Z}$-graded objects in~$\mathcal{A}$,
    and the forgetful map
    $\mathrm{Grad} (\mathcal{X}) \to \mathcal{X}$
    forgets the $\mathbb{Z}$-grading.
    As another example, if $\mathcal{X} = * / G$
    for a linear algebraic group~$G$ over~$\mathbb{C}$,
    then $\mathrm{Grad} (\mathcal{X})$ is a disjoint union
    of stacks of the form~$* / L$
    for Levi subgroups $L \subset G$.
    See \cite{halpern-leistner-instability,epsilon-i}
    for details and more examples.

    In this case, the vertex induction is the operation
    \begin{align*}
        \mathrm{H}_{\bullet + 2 \vdim} (\mathcal{X}_\alpha; \mathbb{Q})
        & \longrightarrow
        \mathrm{H}_{\bullet + 2 \vdim} (\mathcal{X}; \mathbb{Q})
        \llparen z_i \rrparen \ ,
        \\
        a
        & \longmapsto
        \alpha_\star \circ \tau (z) \bigl(
            a \cap e_z^{-1} (\bbnu)
        \bigr) \ ,
    \end{align*}
    where
    $\mathcal{X}_\alpha \subset \mathrm{Grad} (\mathcal{X})$
    is a connected component;
    $\alpha_\star$ denotes pushing forward along the forgetful map
    $\mathcal{X}_\alpha \to \mathcal{X}$;
    $z_1, \dotsc, z_n$ are formal variables,
    where~$n$ depends on~$\alpha$;
    $\tau (z)$ is a certain translation operator;
    and $\bbnu$
    is the normal bundle of the map
    $\mathcal{X}_\alpha \to \mathcal{X}$.

    For example, if~$\mathcal{X}$ is the moduli stack
    of objects in an abelian category~$\mathcal{A}$,
    then each~$\mathcal{X}_\alpha$
    is isomorphic to an $n$-fold product of components of~$\mathcal{X}$
    for some~$n$, and the vertex induction gives
    the $n$-fold multiplication in the Joyce vertex algebra.
    If, moreover, $\mathcal{A}$
    is equipped with a contravariant involution~$(-)^\vee$
    as in~\cref{para-intro-vm},
    then each component $(\mathcal{X}^{\mathbb{Z}_2})_\alpha$
    is isomorphic to a component of
    $\mathcal{X}^n \times \mathcal{X}^{\mathbb{Z}_2}$
    for some~$n$, and the vertex induction for~$\mathcal{X}^{\mathbb{Z}_2}$
    gives the multiplication
    $a_1 (z_1) \cdots a_n (z_n) \, m$
    in the twisted module.
\end{para}

\begin{para}[DT4 invariants]
    \label{para-intro-dt4}
    In \cref{subsec-real-va},
    we discuss a variant of Joyce vertex algebras
    which appear in the context of \emph{DT4 invariants}.
    These are invariants counting sheaves on Calabi--Yau fourfolds,
    developed by
    \textcite{cao-leung-dt4},
    \textcite{borisov-joyce-2017},
    and \textcite{oh-thomas-2023-i,oh-thomas-ii}.
    Their existence and wall-crossing formulae
    were conjectured by
    \textcite[\S 4.4]{gross-joyce-tanaka-2022},
    and have not yet been proved at the time of this writing.

    We define a similar variant of vertex induction
    on the homology of
    \emph{oriented $n$-shifted symplectic stacks}
    for $n \in 4 \mathbb{Z} + 2$.
    For example,
    $(-2)$-shifted symplectic structures
    exist on moduli stacks of sheaves on Calabi--Yau fourfolds.
    The existence of orientations in this case is more subtle,
    but is recently shown for a class of Calabi--Yau fourfolds
    by \textcite{joyce-upmeier-bordism}.

    We expect that a generalized version of DT4 invariants should be defined
    for a general class of
    oriented $(-2)$-shifted symplectic stacks,
    especially in the non-linear cases.
    We also expect that wall-crossing formulae
    for these generalized DT4 invariants
    should be written down using the variant of vertex induction that we define.
\end{para}

\begin{para}[\textit{K}-theoretic invariants]
    \label{para-intro-k-inv}
    In \cref{subsec-k-va},
    we consider another variant of Joyce vertex algebras,
    which are multiplicative vertex algebra structures
    defined on the topological $K$-homology of moduli stacks,
    instead of the homology.
    This type of structure was originally due to \textcite{liu-k-wall-crossing},
    who also constructed $K$-theoretic enumerative invariants
    living in a version of $K$-homology,
    and proved their wall-crossing formulae
    using this multiplicative vertex algebra.

    We generalize this construction
    to define a $K$-homology version of vertex induction.
    We expect that a generalized version of Liu's $K$-theoretic invariants
    should be defined for a general class of quasi-smooth stacks,
    and that they should satisfy wall-crossing formulae
    expressed using the $K$-theoretic vertex induction.
\end{para}

\begin{para}[Acknowledgements]
    The author thanks Dominic Joyce
    for many helpful comments and suggestions,
    and Henry Liu for discussions related to the $K$-theory version.

    This work was done during the author's PhD programme
    supported by the Mathematical Institute,
    University of Oxford.
\end{para}

\begin{para}[Note]
    This work partially supersedes the author's preprint
    \cite{bu-self-dual-ii},
    where the twisted module of Joyce's vertex algebra was first constructed.
\end{para}

\resumeparabookmarks
\section{Joyce vertex algebras}

\label{sec-va}

\subsection{Vertex algebras}
\label{subsec-va}

\begin{para}
    We start by providing background on vertex algebras.
    We recall their usual definition
    in \cref{para-va-def},
    and then state a convenient alternative definition
    in \cref{para-va-def-alt}.
    We also discuss a notion of \emph{weak vertex algebras},
    which will arise in our constructions below.
\end{para}

\begin{para}[Formal power series]
    Let $K$ be a field.
    For a finite-dimensional $K$-vector space~$\Lambda$,
    let $K \llbr \Lambda \rrbr$ be the $K$-algebra of formal power series
    on~$\Lambda$ with coefficients in~$K$, defined as
    \begin{equation*}
        K \llbr \Lambda \rrbr =
        \prod_{n = 0}^\infty \mathrm{Sym}^n (\Lambda^\vee) \ .
    \end{equation*}
    Elements of $K \llbr \Lambda \rrbr$ are often written as~$f (z)$,
    with a variable $z \in \Lambda$.
    Let~$K \llparen \Lambda \rrparen$ be the fraction field of~$K \llbr \Lambda \rrbr$.
    For a $K$-vector space~$V$, define
    \begin{equation*}
        V \llbr \Lambda \rrbr =
        \prod_{n = 0}^\infty V \underset{K}{\otimes} \mathrm{Sym}^n (\Lambda^\vee) \ ,
        \qquad
        V \llparen \Lambda \rrparen =
        V \llbr \Lambda \rrbr \underset{K \llbr \Lambda \rrbr}{\otimes}
        K \llparen \Lambda \rrparen \ .
    \end{equation*}

    For finite-dimensional $K$-vector spaces~$\Lambda_1, \dotsc, \Lambda_n$,
    writing $\Lambda = \Lambda_1 \oplus \cdots \oplus \Lambda_n$,
    we have a natural isomorphism
    $K \llbr \Lambda \rrbr \simeq
    K \llbr \Lambda_1 \rrbr \cdots \llbr \Lambda_n \rrbr$,
    and an injective $K \llbr \Lambda \rrbr$-module homomorphism
    \begin{equation}
        \label{eq-def-iota}
        \iota_{\Lambda_1, \dotsc, \Lambda_n} \colon
        K \llparen \Lambda \rrparen
        \longhookrightarrow
        K \llparen \Lambda_1 \rrparen \cdots \llparen \Lambda_n \rrparen \ ,
    \end{equation}
    given by power series expansion in the given order.
    We also denote this map by
    $\iota_{z_1, \dotsc, z_n}$
    if $z_i$ is a set of coordinates of~$\Lambda_i$.
    For example, we have
    \begin{equation}
        \iota_{z, w} \Bigl(
            \frac{1}{z + w}
        \Bigr) =
        \sum_{n = 0}^\infty {}
        (-1)^n \cdot
        \frac{w^n}{z^{n + 1}} \ ,
    \end{equation}
    where the notation means that we take
    $\Lambda_1 \simeq \Lambda_2 \simeq K$,
    and $z$, $w$ are coordinates on~$\Lambda_1$, $\Lambda_2$.

    This also induces a map
    $\iota_{\Lambda_1, \dotsc, \Lambda_n} \colon
    V \llparen \Lambda \rrparen \to V \llparen \Lambda_1 \rrparen \cdots \llparen \Lambda_n \rrparen$
    for any $K$-vector space~$V$.
\end{para}

\begin{para}[Vertex algebras]
    \label{para-va-def}
    We recall the usual definition of a vertex algebra,
    taken from \textcite[Definition~1.3.1]{frenkel-ben-zvi-2004}.

    A \emph{$\mathbb{Z}$-graded vertex algebra} over~$K$
    is a quadruple $(V, 1, D, Y)$,
    consisting of a $\mathbb{Z}$-graded $K$-vector space~$V$,
    an element $1 \in V$ of degree~$0$,
    a linear map $D \colon V \to V$
    called the \emph{translation operator},
    and a linear map $Y (-, z) (-) \colon V \otimes V \to V \llparen z \rrparen$,
    satisfying the following axioms:

    \begin{enumerate}
        \item
            \label{item-va-unit}
            (\emph{Unit})
            For all $a \in V$, we have
            $Y (a, z) (1) \in a + z V \llbr z \rrbr$, and
            $Y (1, z) (a) = a$.

        \item
            (\emph{Translation})
            For any $a, b \in V$, we have
            \begin{equation}
                \label{eq-def-va-transl}
                [D, Y (a, z)] (b) =
                \frac{\partial}{\partial z} \, Y (a, z) (b),
            \end{equation}
            where $[-,-]$ denotes the commutator,
            and we have $D (1) = 0$.

        \item
            \label{item-va-locality}
            (\emph{Locality})
            For homogeneous elements $a, b, c \in V$, the elements
            \begin{align*}
                Y (a, z) \circ Y (b, w) (c)
                & \in V \llparen z \rrparen \llparen w \rrparen \ ,
                \\
                (-1)^{|a| \, |b|} \cdot
                Y (b, w) \circ Y (a, z) (c)
                & \in V \llparen w \rrparen \llparen z \rrparen
            \end{align*}
            are expansions of the same element in
            $V \llbr z, w \rrbr \, [z^{-1}, w^{-1}, (z - w)^{-1}]$,
            where $|{-}|$ denotes the $\mathbb{Z}$-grading.
    \end{enumerate}
\end{para}

\begin{para}[An equivalent definition]
    \label{para-va-def-alt}
    We introduce an alternative definition
    of a vertex algebra,
    following ideas from \textcite{kim-2011}.
    This definition is perhaps more naturally motivated than the standard one,
    and will be important for generalizing Joyce vertex algebras
    to other situations below.
    See also \cref{thm-va-func-def} below
    for a functorial description of this approach.

    A $\mathbb{Z}$-graded vertex algebra over $K$
    is equivalently the data $(V, (X_n)_{n \geq 0})$,
    consisting of a $\mathbb{Z}$-graded $K$-vector space~$V$
    and $K$-linear multiplication maps
    \begin{align*}
        X_n \colon
        V^{\otimes n}
        & \longrightarrow
        V \llbr z_1, \dotsc, z_n \rrbr \, [ (z_i - z_j)^{-1} ] \ ,
        \\
        a_1 \otimes \cdots \otimes a_n
        & \longmapsto
        X_n (a_1, \dotsc, a_n; z_1, \dotsc, z_n) \ ,
    \end{align*}
    preserving grading, where $\deg z_i = -2$,
    and we invert $z_i - z_j$ for $i \neq j$.
    In particular, this includes a map
    $X_0 \colon K \to V$, thought of as the unit.
    They should satisfy the following properties:

    \begin{enumerate}
        \item
            \label{item-va-alt-unit}
            (\emph{Unit})
            For any $a \in V$, we have
            \begin{equation}
                \label{eq-va-alt-unit}
                X_1 (a; 0) = a \ .
            \end{equation}
            More precisely, we have
            $X_1 (a; z) \in a + z V \llbr z \rrbr
            \subset V \llbr z \rrbr$.

        \item
            \label{item-va-alt-comm}
            (\emph{Commutativity})
            For any homogeneous elements $a_1, \dotsc, a_n \in V$,
            and any permutation $\sigma \in \mathfrak{S}_n$,
            we have
            \begin{equation}
                \label{eq-va-alt-comm}
                X_n (a_{\sigma (1)}, \dotsc, a_{\sigma (n)}; z_{\sigma (1)}, \dotsc, z_{\sigma (n)})
                = \pm X_n (a_1, \dotsc, a_n; z_1, \dotsc, z_n) \ ,
            \end{equation}
            where the sign is determined by the Koszul sign rule:
            it is `$-$' if and only if $\sigma$
            restricts to an odd permutation on the odd-graded elements.

        \item 
            \label{item-va-alt-assoc}
            (\emph{Associativity})
            For integers $m, n \geq 0$ and elements
            $b_1, \dotsc, b_m, a_1, \dotsc, a_n \in V$,
            we have
            \begin{multline}
                \label{eq-va-alt-assoc}
                X_{n+1} \Bigl(
                    X_m (b_1, \dotsc, b_m; w_1, \dotsc, w_m),
                    a_1, \dotsc, a_n; \
                    z_0, \dotsc, z_n
                \Bigr)
                \\[-.5ex]
                =
                \iota_{ \{ z_i \}, \{ w_j \} } \,
                X_{m+n} \bigl(
                    b_1, \dotsc, b_m, a_1, \dotsc, a_n;
                    z_0 + w_1, \dotsc, z_0 + w_m,
                    z_1, \dotsc, z_n
                \bigr) \ ,
            \end{multline}
            where $\iota_{ \{ z_i \}, \{ w_j \} }$
            is the map defined in \cref{eq-def-iota}.
    \end{enumerate}
    As is common in the literature,
    we adopt the convenient notation
    \begin{equation}
        \label{eq-va-notation}
        a_1 (z_1) \cdots a_n (z_n) =
        X_n (a_1, \dotsc, a_n; z_1, \dotsc, z_n) \ ,
    \end{equation}
    but the reader should be aware that
    this is not defined using
    the individual terms $a_i (z_i) = X_1 (a_i; z_i)$.
\end{para}

\begin{para}[Proof of the equivalence of definitions]
    \label{para-va-alt-proof}
    Given a graded vertex algebra~$V$
    as in \cref{para-va-def},
    define the maps $X_n$ in \cref{para-va-def-alt} by
    \begin{equation}
        X_n (a_1, \dotsc, a_n; z_1, \dotsc, z_n) =
        \iota _{z_1, \dotsc, z_n} ^{-1}
        [ Y (a_1, z_1) \circ \cdots \circ Y (a_n, z_n) (1) ] \ ,
    \end{equation}
    where we take the unique preimage under the embedding
    \[
        \iota_{z_1, \dotsc, z_n} \colon
        V \llbr z_1, \dotsc, z_n \rrbr \, [(z_i - z_j)^{-1}]
        \longhookrightarrow
        V \llparen z_1 \rrparen \cdots \llparen z_n \rrparen \ .
    \]
    By \textcite[Theorem~3.14]{kim-2011},
    this is well-defined,
    and satisfies the property \cref{para-va-def-alt}~\cref{item-va-alt-assoc}.
    The properties
    \cref{para-va-def-alt}~\crefrange{item-va-alt-unit}{item-va-alt-comm}
    follow from
    \cref{para-va-def}~\cref{item-va-unit},~\cref{item-va-locality}.

    Conversely, given the data $(V, (X_n)_{n \geq 0})$
    as in \cref{para-va-def-alt},
    we may define a $\mathbb{Z}$-graded vertex algebra structure on~$V$
    by setting
    \begin{equation*}
        1 = X_0 (1) \ ,
        \qquad
        D (a) =
        \frac{\partial}{\partial z} \, X_1 (a; z) \Big|_{z = 0} \ ,
        \qquad
        Y (a, z) (b) = X_2 (a, b; z, 0)
    \end{equation*}
    for all $a, b \in V$.
    Verifying the axioms is straightforward
    except for~\cref{eq-def-va-transl},
    which we prove as follows.
    Applying the associativity relation~\cref{eq-va-alt-assoc}
    multiple times, we have
    \begin{align*}
        [D, Y (a, z)] (b)
        & = \frac{\partial}{\partial w} \Bigl(
            X_1 \bigl( X_2 (a, b; z, 0), w \bigr)
            - X_2 (a, X_1 (b, w); z, 0)
        \Bigr) \Big|_{w = 0}
        \\
        & = \frac{\partial}{\partial w} \bigl(
            X_2 (a, b; z + w, w) -
            X_2 (a, b; z, w)
        \bigr) \Big|_{w = 0}
        \\
        & = \frac{\partial}{\partial z} \, X_2 (a, b; z, 0)
        = \frac{\partial}{\partial z} \, Y (a, z) (b) \ ,
    \end{align*}
    as desired.
    \qed
\end{para}

\begin{para}[Remark]
    \label{para-va-3-terms}
    In fact, in \cref{para-va-def-alt},
    it is enough to require the operators
    $X_n$ for $n = 0, 1, 2, 3$,
    and only require commutativity and associativity up to $3$ terms,
    since these are enough for converting to the usual definition,
    and we can then convert it back to obtain all the higher~$X_n$.
\end{para}

\begin{para}[Weak vertex algebras]
    \label{para-weak-va}
    We introduce a generalized version of vertex algebras,
    where for the product $a_1 (z_1) \cdots a_n (z_n)$,
    in addition to the usual poles at $z_i = z_j$ for $i \neq j$,
    we also allow arbitrary poles along any divisor.
    As we will see, this type of structure arises naturally from geometry,
    in the set-up of~\cref{subsec-va-setup}.

    A \emph{$\mathbb{Z}$-graded weak vertex algebra} over~$K$
    is the data $(V, (X_n)_{n \geq 0})$,
    consisting of a $\mathbb{Z}$-graded $K$-vector space~$V$
    and $K$-linear multiplication maps
    \begin{equation*}
        X_n \colon
        V^{\otimes n}
        \longrightarrow
        V \llparen z_1, \dotsc, z_n \rrparen \ ,
    \end{equation*}
    preserving grading, where $\deg z_i = -2$,
    satisfying the properties
    \cref{para-va-def-alt}~\cref{item-va-alt-unit,item-va-alt-comm,item-va-alt-assoc}.

    Therefore, such a weak vertex algebra is a vertex algebra
    if and only if the image of~$X_n$ lies in the subspace
    $V \llbr z_1, \dotsc, z_n \rrbr \, [(z_i - z_j)^{-1}]
    \subset V \llparen z_1, \dotsc, z_n \rrparen$.
\end{para}

\subsection{Moduli spaces}
\label{subsec-va-setup}

\begin{para}
    \label{para-va-notations}
    As we mentioned in the introduction,
    Joyce vertex algebras
    \cite{joyce-hall,joyce-homological,joyce-surfaces}
    are vertex algebra structures defined on the homology of
    certain moduli spaces.
    Here, we formulate a set of axioms for such moduli spaces.

    We use the following terminology and notations:

    \begin{itemize}
        \item
            A \emph{torus} is a Lie group~$T$
            isomorphic to $\mathrm{U} (1)^n$ for some $n \geq 0$.

        \item
            For a torus~$T$,
            let $\Lambda^T$, $\Lambda_T$
            be the weight and coweight lattices,
            $\Lambda^T = \mathrm{Hom} (T, \mathrm{U} (1))$
            and $\Lambda_T = \mathrm{Hom} (\mathrm{U} (1), T)$,
            both isomorphic to $\mathbb{Z}^{\dim T}$.

        \item
            $\mathsf{hCW}$ is the category
            whose objects are topological spaces
            that are homotopy equivalent to CW~complexes,
            and morphisms are homotopy classes of continuous maps.

        \item
            $K (X)$ is the \emph{topological $K$-theory}
            of a space $X \in \mathsf{hCW}$,
            defined as the abelian group
            $K (X) = \mathsf{hCW} (X, \mathrm{BU} \times \mathbb{Z})$.
    \end{itemize}
\end{para}

\begin{para}[The category \texorpdfstring{$\mathcal{T}$}{T}]
    \label{para-cat-t}
    For convenience of presentation,
    we introduce a category~$\mathcal{T}$
    of spaces with $\mathrm{B} T$-actions for tori~$T$,
    defined as follows:

    \begin{itemize}
        \item
            Its objects are triples~$(T, X, \odot)$,
            where~$T \simeq \mathrm{U} (1)^n$ is a torus,
            $X \in \mathsf{hCW}$ is a space, and
            \begin{equation*}
                \odot \colon \mathrm{B} T \times X \longrightarrow X
            \end{equation*}
            is a map, such that it defines a
            $\mathrm{B} T$-action on~$X$ in~$\mathsf{hCW}$.

        \item
            A morphism $f \colon (T_1, X_1, \odot_1) \to (T_2, X_2, \odot_2)$
            consists of a Lie group homomorphism
            $f^{\smash{\sharp}} \colon T_2 \to T_1$,
            together with a $\mathrm{B} T_2$-equivariant map
            $f \colon X_1 \to X_2$ in~$\mathsf{hCW}$.
    \end{itemize}
    We often abbreviate the triple~$(T, X, \odot)$ as~$X$,
    and call $\dim T$ the \emph{rank} of such an object.

    The category~$\mathcal{T}$ admits finite products,
    giving a symmetric monoidal structure on~$\mathcal{T}$.
\end{para}

\begin{para}[Weights in \textit{K}-theory]
    \label{para-k-circ}
    Let $(T, X, \odot)$ be an object of~$\mathcal{T}$.
    For each weight $\lambda \in \Lambda^T$,
    let $K (X)_\lambda \subset K (X)$
    be the subgroup of classes~$E$
    \emph{of weight~$\lambda$},
    meaning that $\odot^* (E) = L_\lambda \boxtimes E$,
    where $L_\lambda \to \mathrm{B} T$
    is the line bundle classified by the map
    $\lambda \colon \mathrm{B} T \to \mathrm{BU} (1)$.
    Let
    \begin{equation*}
        K^\circ (X) =
        \bigoplus_{\lambda \in \Lambda \lowersup{T}} K (X)_\lambda
        \subset K (X)
    \end{equation*}
    be the subgroup of classes
    which admit finite weight decompositions.
\end{para}

\begin{para}[The setting]
    \label{para-va-setting}
    We assume given the data
    $(X, \odot, \oplus, 0, \mathbb{T}_X)$, where

    \begin{itemize}
        \item
            $X \in \mathsf{hCW}$ is a space.

        \item
            $\odot \colon \mathrm{BU} (1) \times X \to X$
            is an action in~$\mathsf{hCW}$,
            giving an object~$(\mathrm{U} (1), X, \odot)$
            of~$\mathcal{T}$ of rank~$1$.

        \item
            $\oplus \colon X \times X \to X$
            and $0 \colon {*} \to X$
            are morphisms in~$\mathcal{T}$,
            with $\oplus^\sharp \colon \mathrm{U} (1) \to \mathrm{U} (1)^2$
            the diagonal map,
            defining a commutative monoid structure
            on~$X$ in~$\mathcal{T}$.

        \item
            $\mathbb{T}_X \in K (X)_0 \subset K (X)$
            is a class of weight~$0$,
            called the \emph{obstruction theory}.
    \end{itemize}
    In this case, the function
    $\operatorname{vdim} = \operatorname{rank} (\mathbb{T}_X)
    \colon \uppi_0 (X) \to \mathbb{Z}$
    is called the \emph{virtual dimension} of~$X$.

    For each $n \geq 0$, let
    $\oplus_{(n)} \colon X^n \to X$
    be the $n$-fold product using~$\oplus$,
    and let $\mathbb{T}_{X^n} =
    \sum_i \mathrm{pr}_i^* (\mathbb{T}_X) \in K (X^n)$,
    where $\mathrm{pr}_i \colon X^n \to X$
    is the $i$-th projection for $i = 1, \dotsc, n$.
    Define the (\emph{virtual}) \emph{normal bundle} of~$\oplus_{(n)}$ as the class
    \begin{equation}
        \label{eq-def-nu}
        \bbnu_{(n)} =
        \oplus_{\smash{(n)}}^* (\mathbb{T}_X) - \mathbb{T}_{X^n}
        \quad {\in} \quad
        K (X^n) \ .
    \end{equation}
    We further assume the following condition:
    \begin{itemize}
        \item
            For any integer $n \geq 0$, we have
            \begin{equation}
                \label{eq-cond-nu}
                \bbnu_{(n)} \in K^\circ (X^n)
                \quad \text{and} \quad
                (\bbnu_{(n)})_0 = 0 \ ,
            \end{equation}
            where $(-)_0$ denotes the part of $\mathrm{U} (1)^n$-weight~$0$.
    \end{itemize}
    This is a technical condition to ensure that
    the equivariant Euler class $e_z^{-1} (\bbnu_{(n)})$
    in \cref{para-equivar-euler} is well-defined,
    and is satisfied in most of our examples.
    See \cref{eg-linear-stacks,eg-dg-cat} below.
\end{para}

\begin{para}[The topological realization of a stack]
    \label{para-top-real}
    To give examples of the data
    in \cref{para-va-setting},
    we will use the \emph{topological realization} functor
    \begin{equation}
        |{-}| \colon
        \mathsf{Fun} (\mathsf{Aff}_\mathbb{C}^{\smash{\mathrm{op}}},
        \mathsf{S})
        \longrightarrow
        \mathsf{S} \ ,
    \end{equation}
    as in \textcite[\S 3.1]{blanc-2016},
    where $\mathsf{Aff}_\mathbb{C}$ is the category of affine $\mathbb{C}$-schemes,
    and~$\mathsf{S}$ is the $\infty$-category of spaces.
    This is defined as the left Kan extension of the functor
    $(-)^\mathrm{an} \colon \mathsf{Aff}_\mathbb{C} \to \mathsf{S}$
    sending an affine $\mathbb{C}$-scheme to the topological space of its analytification.

    For a derived algebraic stack~$\mathcal{X}$ over~$\mathbb{C}$,
    we denote by~$|\mathcal{X}|$
    the topological realization of its classical truncation
    defined as above.
\end{para}

\begin{example}
    \label{eg-linear-stacks}
    We give a class of examples of the data
    $(X, \odot, \oplus, 0, \mathbb{T}_X)$
    as in \cref{para-va-setting}.

    Let~$\mathcal{X}$ be a
    \emph{derived linear moduli stack} over~$\mathbb{C}$
    in the sense of \textcite[\S 7.1]{epsilon-i}
    and \cite[\S 2.4.6]{bu-davison-ibanez-nunez-kinjo-padurariu},
    so that it is equipped with a monoid structure and an action
    \begin{align*}
        \oplus \colon \mathcal{X} \times \mathcal{X}
        & \longrightarrow \mathcal{X} \ ,
        \\
        \odot \colon \mathrm{B} \mathbb{G}_\mathrm{m} \times \mathcal{X}
        & \longrightarrow \mathcal{X} \ ,
    \end{align*}
    satisfying higher coherence and compatibility conditions.
    Typical examples include the following:

    \begin{itemize}
        \item
            The derived moduli stack of finite-dimensional representations
            of a quiver, possibly with relations or a potential.

        \item
            The derived moduli stack of vector bundles
            or coherent sheaves
            on a smooth, projective $\mathbb{C}$-variety.
    \end{itemize}
    Let $X = |\mathcal{X}|$ be the topological realization
    as in \cref{para-top-real}.
    The data~$\odot, \oplus, 0$ on~$X$
    are defined by those on~$\mathcal{X}$.
    The class~$\mathbb{T}_X$ is defined as the class of
    the tangent complex $\mathbb{T}_\mathcal{X}$ of~$\mathcal{X}$.

    The condition~\cref{eq-cond-nu}
    usually holds in this case,
    since if~$\mathcal{X}$ is a moduli stack of objects
    in a $\mathbb{C}$-linear abelian category~$\mathcal{A}$,
    then we often have
    $\mathbb{T}_\mathcal{X} |_E \simeq \mathrm{Ext}_{\mathcal{A}}^{1+\bullet} (E, E)$
    at a point $E \in \mathcal{A}$, so
    $\oplus_{\smash{(n)}}^* (\mathbb{T}_\mathcal{X}) |_{E_1, \dotsc, E_n} \simeq
    \bigoplus_{i, j} \mathrm{Ext}_{\mathcal{A}}^{1+\bullet} (E_i, E_j)$
    only has $\mathbb{G}_\mathrm{m}^n$-weights
    of the form $e_i - e_j$ for $i, j \in \{ 1, \dotsc, n \}$,
    and the weight~$0$ part
    $\bigoplus_i \mathrm{Ext}_{\mathcal{A}}^{1+\bullet} (E_i, E_i)$
    agrees with $\mathbb{T}_{\mathcal{X}^n} |_{E_1, \dotsc, E_n}$.
\end{example}

\begin{example}
    \label{eg-dg-cat}
    Another class of examples of
    $(X, \odot, \oplus, 0, \mathbb{T}_X)$
    as in \cref{para-va-setting}
    are given as follows.

    Let~$\mathcal{C}$ be a $\mathbb{C}$-linear dg-category
    of finite type in the sense of
    \textcite[Definition~2.4]{toen-vaquie-2007-moduli},
    and let~$\mathcal{X}$ be the derived moduli stack
    of right proper objects in~$\mathcal{C}$,
    as in \cite[Theorem~3.6]{toen-vaquie-2007-moduli}.
    Examples include the following:

    \begin{itemize}
        \item
            The derived moduli stack of
            finite-dimensional complexes of representations
            of a quiver.
        \item
            The derived moduli stack of
            perfect complexes on a smooth, projective $\mathbb{C}$-variety.
    \end{itemize}
    Let $X = |\mathcal{X}|$ be the topological realization of~$\mathcal{X}$,
    and define the data $\odot, \oplus, 0, \mathbb{T}_X$
    using the scalar product,
    the direct sum, the zero object,
    and the tangent complex of~$\mathcal{X}$.

    In this case, the condition~\cref{eq-cond-nu} always holds,
    by the description of the tangent complex of~$\mathcal{X}$ by
    \textcite[Proposition~3.9]{brav-dyckerhoff-2021}.

    For example, if $\mathcal{X} = \calPerf$
    is the classifying stack of perfect complexes over~$\mathbb{C}$,
    as in \textcite[Proposition~3.7]{toen-vaquie-2007-moduli},
    then we have an equivalence
    $|\calPerf| \simeq \mathrm{BU} \times \mathbb{Z}$,
    by \textcite[Theorems~4.5 and~4.21]{blanc-2016}.
\end{example}

\subsection{Joyce vertex algebras}
\label{subsec-joyce-va}

\begin{para}[The equivariant Euler class]
    \label{para-equivar-euler}
    Let~$(T, X, \odot) \in \mathcal{T}$,
    and let~$E \in K^\circ (X)$
    be a class such that the weight~$0$ part~$E_0$
    is the class of a vector bundle on~$X$,
    possibly of mixed rank. Define
    \begin{equation}
        \label{eq-def-equivar-euler}
        e_z (E) =
        e (E_0) \cdot
        \prod_{\lambda \in \Lambda \lowersup{T} \setminus \{ 0 \}}
        \sum_{i = 0}^\infty
        \lambda (z)^{\operatorname{rank} E_\lambda - i} \cdot
        c_i (E_\lambda)
        \quad {\in} \quad
        \prod_{k = 0}^\infty
        \mathrm{H}^{2 k} (X; \mathbb{Q})
        ( z_1, \dotsc, z_n ) \ ,
    \end{equation}
    where $z = (z_1, \dotsc, z_n)$
    is a set of coordinates on~$\Lambda_T$,
    and $n = \dim T$.
    For each $k \geq 0$,
    the part of~$e_z (E)$ lying in~$\mathrm{H}^{2k}$
    is a rational function in~$z_1, \dotsc, z_n$,
    with possible poles at $\lambda (z) = 0$ for $\lambda \in \Lambda^T$
    such that $E_\lambda$ is not the class of a vector bundle,
    although there may not be a uniform bound
    for the orders of poles as~$k$ increases.

    We have the relation
    \begin{equation*}
        e_z (E + F) = e_z (E) \, e_z (F) \ ,
    \end{equation*}
    and in particular, if $E_0 = 0$, then
    $e_z (E) \, e_z (-E) = 1$,
    in which case we also write
    $e_z^{-1} (E) = e_z (-E)$.
\end{para}

\begin{para}[The translation operator]
    \label{para-translation-op}
    Let $(T, X, \odot) \in \mathcal{T}$.
    Define its \emph{translation operator}
    \begin{equation*}
        \tau (z) =
        \exp (z_1 D_1 + \cdots + z_n D_n)
        \colon \mathrm{H}_{\bullet} (X; \mathbb{Q})
        \longrightarrow
        \mathrm{H}_{\bullet} (X; \mathbb{Q})
        \llbr z_1, \dotsc, z_n \rrbr \ ,
    \end{equation*}
    with $\deg z_i = -2$ and $z = (z_1, \dotsc, z_n)$,
    where we choose an identification $T \simeq \mathrm{U} (1)^n$,
    and define $D_i \colon \mathrm{H}_\bullet (X; \mathbb{Q})
    \to \mathrm{H}_{\bullet + 2} (X; \mathbb{Q})$
    by
    \begin{equation*}
        D_i (a) = \odot_* (t_i \boxtimes a) \ ,
    \end{equation*}
    where $t_i \in \mathrm{H}_2 (\mathrm{B} T; \mathbb{Q})$
    is the generator of the $i$-th copy of
    $\mathrm{H}_2 (\mathrm{BU} (1); \mathbb{Q})$,
    dual to the universal first Chern class in
    $\mathrm{H}^2 (\mathrm{BU} (1); \mathbb{Q})$.

    This construction does not depend on the choice of
    identification $T \simeq \mathrm{U} (1)^n$,
    provided that the coordinates~$z_i$ on~$\Lambda_T$
    are transformed accordingly.
\end{para}

\begin{theorem}
    \label{thm-joyce-va}
    Let~$X$ be as in \cref{para-va-setting}.
    Then the assignment
    \begin{equation}
        \label{eq-def-joyce-va}
        a_1 (z_1) \cdots a_n (z_n) =
        (\oplus_{(n)})_* \circ \tau (z)
        \Bigl(
            (a_1 \boxtimes \cdots \boxtimes a_n) \cap
            e_z^{-1} (\bbnu_{(n)})
        \Bigr)
    \end{equation}
    equips the space
    $V = \mathrm{H}_{\bullet + 2 \vdim} (X; \mathbb{Q})$
    with the structure of a $\mathbb{Z}$-graded weak vertex algebra over~$\mathbb{Q}$,
    where $z = (z_1, \dotsc, z_n)$,
    and~$\tau (z)$ is the translation operator of
    the object~$X^n \in \mathcal{T}$.

    The product \cref{eq-def-joyce-va} has poles along
    $\lambda (z) = 0$
    for non-zero weights $\lambda \in \mathbb{Z}^n$
    appearing in~$\oplus_{(n)}^* (\mathbb{T}_X)$.
    In particular, if the class
    $\oplus^* (\mathbb{T}_X) \in K^\circ (X^2)$
    only has $\mathrm{U} (1)^2$-weights
    $(-1, 1),$ $(0, 0),$ and~$(1, -1),$
    then $V$ is a $\mathbb{Z}$-graded vertex algebra.
\end{theorem}

\begin{proof}
    \allowdisplaybreaks
    The property \cref{eq-va-alt-unit} holds since
    $X_1 (a; 0) = \exp (0) (a) = a$.
    The relation \cref{eq-va-alt-comm}
    holds since the definition~\cref{eq-def-joyce-va}
    is symmetric (with a sign rule) in the indices $1, \dotsc, n$.

    We now verify~\cref{eq-va-alt-assoc}.
    Write $a = a_1 \boxtimes \cdots \boxtimes a_n$
    and $b = b_1 \boxtimes \cdots \boxtimes b_m$ for short.
    Expanding the left-hand side of~\cref{eq-va-alt-assoc},
    we obtain
    \begin{align*}
        &
        {\oplus}_* \circ
        \exp \bigl( z_0 D_j + z_i D_{m + i} \bigr)
        \Bigl[
            \Bigl(
                \exp (w_j D_j)
                \bigl(
                    b \cap e_w^{-1} (\bbnu_{(m)})
                \bigr)
                \boxtimes a
            \Bigr)
            \cap e_{z_0, z}^{-1} (\bbnu_{(n+1)})
        \Bigr]
        \\[1ex]
        = {}
        &
        {\oplus}_* \circ
        \exp \bigl( z_0 D_j + z_i D_{m + i} \bigr)
        \Bigl[
            \exp (w_j D_j)
            \Bigl(
                (b \boxtimes a)
                \cap e_w^{-1} \bigl(
                    \oplus_{\smash{(m)}}^* (\mathbb{T}_X)
                    - (\mathbb{T}_1 + \cdots + \mathbb{T}_m)
                \bigr)
            \Bigr)
        \\*
        & \hspace{2em}
            {} \cap e_{z_0, z}^{-1} \Bigl(
                \oplus^* (\mathbb{T}_X) - \bigl(
                    \oplus_{\smash{(m)}}^* (\mathbb{T}_X)
                    + \mathbb{T}_{m+1} + \cdots + \mathbb{T}_{m+n}
                \bigr)
            \Bigr)
        \Bigr]
        \\[1ex]
        = {}
        &
        {\oplus}_* \circ
        \exp \bigl( z_0 D_j + z_i D_{m + i} \bigr) \circ
        \exp (w_j D_j)
        \Bigl[
            (b \boxtimes a)
            \cap e_w^{-1} \bigl(
                \oplus_{\smash{(m)}}^* (\mathbb{T}_X)
                - (\mathbb{T}_1 + \cdots + \mathbb{T}_m)
            \bigr)
        \\*
        & \hspace{2em}
            {} \cap e_{z_0 + w, z}^{-1} \Bigl(
                \oplus^* (\mathbb{T}_X) - \bigl(
                    \oplus_{\smash{(m)}}^* (\mathbb{T}_X)
                    + \mathbb{T}_{m+1} + \cdots + \mathbb{T}_{m+n}
                \bigr)
            \Bigr)
        \Bigr]
        \\[1ex]
        = {}
        &
        {\oplus}_* \circ
        \exp \bigl( z_0 D_j + z_i D_{m + i} \bigr) \circ
        \exp (w_j D_j)
        \Bigl[
            (b \boxtimes a)
            \cap e_{z_0 + w, z}^{-1} \bigl(
                \oplus^* (\mathbb{T}_X) - (
                    \mathbb{T}_1 + \cdots + \mathbb{T}_{m+n}
                )
            \bigr)
        \Bigr]
        \\[1ex]
        = {}
        &
        {\oplus}_* \circ
        \exp \bigl( (z_0 + w_j) D_j + z_i D_{m + i} \bigr)
        \Bigl(
            (b \boxtimes a)
            \cap e_{z_0 + w, z}^{-1}
            (\bbnu_{(m+n)})
        \Bigr) \ ,
    \end{align*}
    where~$\oplus$ is short for~$\oplus_{(m+n)}$,
    and $z$ is short for $z_1, \dotsc, z_n$.
    The indices $i, j$ are implicitly summed over the ranges
    $1 \leq i \leq n$ and $1 \leq j \leq m$, respectively.
    For $i = 1, \dotsc, m+n$,
    we denote by~$\mathbb{T}_i$ the pullback of~$\mathbb{T}_X$
    from the $i$-th copy of~$X$.
    The second step uses \cref{lem-swap} below,
    and the third step uses the fact that~$\bbnu_{(m)}$
    has weight zero with respect to the diagonal torus,
    associated to the variable~$z_0$.
    This computation agrees with the right-hand side of~\cref{eq-va-alt-assoc}.

    For the final statement on usual vertex algebras,
    it is enough to show that for each $n \geq 0$,
    the class $\oplus_{\smash{(n)}}^* (\mathbb{T}_X)$
    can only have $\mathrm{U} (1)^n$-weights of the form
    $e_i - e_j$ for $i, j \in \{ 1, \dotsc, n \}$,
    where $e_i = (0, \dotsc, 0, 1, 0, \dotsc, 0)$
    with~$1$ at the $i$-th position.
    This is because the assumption implies that
    for any weight $(k_1, \dotsc, k_n)$ appearing
    in~$\oplus_{\smash{(n)}}^* (\mathbb{T}_X)$
    and any partition $\{ 1, \dotsc, n \} = I \sqcup J$,
    writing $k_I = \sum_{i \in I} k_i$ and $k_J = \sum_{j \in J} k_j$,
    we have $(k_I, k_J) \in \{ (-1, 1), (0, 0), (1, -1) \}$,
    which implies the desired property.
\end{proof}

\begin{lemma}
    \label{lem-swap}
    Let $(T, X, \odot) \in \mathcal{T}$,
    and let~$E \in K^\circ (X)$,
    such that its weight~$0$ part~$E_0$
    is the class of a vector bundle on~$X$.
    Then for any $a \in \mathrm{H}_\bullet (X; \mathbb{Q})$, we have
    \begin{equation}
        \tau (w) (a)
        \cap e_z (E)
        = \iota_{z, w}
        \bigl(
            \tau (w) (a \cap e_{z + w} (E))
        \bigr)
        \quad \in \quad
        \mathrm{H}_\bullet (X; \mathbb{Q})
        \llparen z \rrparen \llbr w \rrbr \ ,
    \end{equation}
    where $w, z$ are sets of~$n$ variables, with $n = \dim T$,
    and $\iota_{z, w}$ denotes the expansion map from
    $\mathrm{H}_\bullet (X; \mathbb{Q}) \llbr w \rrbr (z + w)$
    to $\mathrm{H}_\bullet (X; \mathbb{Q}) \llparen z \rrparen \llbr w \rrbr$.
\end{lemma}

\begin{proof}
    Choose an identification $T \simeq \mathrm{U} (1)^n$,
    and let $\delta_i \in \mathrm{H}^2 (\mathrm{B} T; \mathbb{Z})$
    be the universal first Chern class in
    the $i$-th copy of~$\mathrm{BU} (1)$,
    where $i = 1, \dotsc, n$.
    We have
    \begin{align}
        {\odot}^* \, e_z (E)
        & =
        \prod_{\lambda \in \Lambda \lowersup{T}}
        \sum_{i = 0}^\infty \lambda (z)^{\operatorname{rank} E_\lambda - i} \cdot
        c_i (L_\lambda \boxtimes E_\lambda)
        \notag \\
        & =
        \prod_{\lambda \in \Lambda \lowersup{T}}
        \sum_{i = 0}^\infty \lambda (z)^{\operatorname{rank} E_\lambda - i} \cdot
        \sum_{j = 0}^i {}
        \binom{\operatorname{rank} E_\lambda - j}{i - j} \,
        \lambda (\delta)^{i - j} \boxtimes c_j (E_\lambda)
        \notag \\
        \label{eq-pf-lem-swap}
        & =
        \prod_{\lambda \in \Lambda \lowersup{T}}
        \sum_{j = 0}^\infty {}
        \lambda (z + \delta)^{\operatorname{rank} E_\lambda - j} \cdot
        c_j (E_\lambda)
        = e_{z + \delta} (E) \ ,
    \end{align}
    where we write $\delta = (\delta_1, \dotsc, \delta_n)$,
    and we expand $e_{z+\delta} (E)$ in non-negative powers of~$\delta_i$.
    In the final line, we pull back~$E$ along the projection
    to the second factor, $\mathrm{B} T \times X \to X$.

    Let $t_i \in \mathrm{H}_2 (\mathrm{B} T; \mathbb{Z})$
    be the dual of~$\delta_i$ as in~\cref{para-translation-op},
    and equip $\mathrm{H}_\bullet (\mathrm{B} T; \mathbb{Q})$
    with a ring structure given by
    the group structure of~$\mathrm{B} T$ in $\mathsf{hCW}$.
    Then
    \begin{align*}
        \tau (w) (a) \cap e_z (E)
        & =
        \odot_* \Bigl[
            \bigl( \exp (w_i t_i) \boxtimes a \bigr) \cap
            \odot^* \, e_z (E)
        \Bigr]
        \\
        & =
        \odot_* \Bigl[
            (a \cap e_{z + \partial / \partial t} (E))
            \exp (w_i t_i)
        \Bigr]
        \\
        & =
        \odot_* \Bigl[
            \iota_{z, w} (a \cap e_{z + w} (E))
            \exp (w_i t_i)
        \Bigr]
        = \iota_{z, w}
        \bigl( \tau (w) ( a \cap e_{z + w} (E) ) \bigr) \ ,
    \end{align*}
    where $z + \partial / \partial t$ denotes the set of variables
    $(z_1 + \partial / \partial t_1, \dotsc, z_n + \partial / \partial t_n)$,
    and we expand $e_{z + \partial / \partial t} (E)$
    in non-negative powers of $\partial / \partial t_i$.
    The second step uses \cref{eq-pf-lem-swap}
    and the fact that $\delta_i$ acts as $\partial / \partial t_i$
    on $\mathrm{H}_\bullet (\mathrm{B} T; \mathbb{Q})$.
    This proves the lemma.
\end{proof}

\begin{para}[Morphisms of Joyce vertex algebras]
    \label[theorem]{thm-joyce-va-mor}
    Finally, we mention a construction of morphisms of Joyce vertex algebras,
    as in \textcite[\S 2.5]{gross-joyce-tanaka-2022}.
\end{para}

\begin{theorem*}
    Let $X, X'$ be two spaces as in \cref{para-va-setting},
    and let $f \colon X \to X'$ be a map in~$\mathsf{hCW}$,
    respecting the monoid structures and $\mathrm{BU} (1)$-actions.
    Assume that the class
    \begin{equation*}
        \mathbb{T}_f = \mathbb{T}_X - f^* (\mathbb{T}_{X'})
    \end{equation*}
    is the class of a vector bundle on~$X$,
    possibly of mixed rank.
    Then the map
    \begin{align*}
        Y_f \colon
        \mathrm{H}_{\bullet + 2 \operatorname{vdim}} (X; \mathbb{Q})
        & \longrightarrow
        \mathrm{H}_{\bullet + 2 \operatorname{vdim}} (X'; \mathbb{Q}) \ ,
        \\
        a
        & \longmapsto
        f_* (a \cap e (\mathbb{T}_f))
    \end{align*}
    is a morphism of\/ $\mathbb{Z}$-graded weak vertex algebras,
    where $e$ denotes the Euler class.

    Moreover, if\/ $f \colon X \to X'$ and $g \colon X' \to X''$
    are maps as above,
    then $Y_{g \circ f} = Y_g \circ Y_f$.
\end{theorem*}

\begin{proof}
    Let $a_1, \dotsc, a_n \in \mathrm{H}_\bullet (X; \mathbb{Q})$.
    Then
    \begin{align*}
        & \phantom{{} = {}}
        Y_f (a_1) (z_1) \cdots Y_f (a_n) (z_n)
        \\[1ex]
        & =
        (\oplus_{(n)})_* \circ \tau (z)
        \smash{
            \Bigl(
                (f^{\times n})_* \Bigl(
                    (a_1 \boxtimes \cdots \boxtimes a_n)
                    \cap e (\mathbb{T}_{\smash{f}}^{\boxplus n})
                \Bigr)
                \cap e_z (\mathbb{T}_{\smash{X'}}^{\boxplus n} - \oplus_{\smash{(n)}}^* (\mathbb{T}_{X'}))
            \Bigr)
        }
        \\[1ex]
        & =
        (\oplus_{(n)})_* \circ \tau (z) \circ (f^{\times n})_*
        \smash{
            \Bigl(
                (a_1 \boxtimes \cdots \boxtimes a_n)
                \cap e_z \Bigl(
                    \mathbb{T}_{\smash{X}}^{\boxplus n}
                    - (f^{\times n})^* \circ \oplus_{\smash{(n)}}^* (\mathbb{T}_{X'})
                \Bigr)
            \Bigr)
        }
        \\[1ex]
        & =
        f_* \circ (\oplus_{(n)})_* \circ \tau (z)
        \smash{
            \Bigl(
                (a_1 \boxtimes \cdots \boxtimes a_n)
                \cap e_z (-\bbnu_{(n)} + \oplus_{\smash{(n)}}^* (\mathbb{T}_f))
            \Bigr)
        }
        \\[1ex]
        & =
        f_* \circ (\oplus_{(n)})_*
        \smash{
            \Bigl(
                \tau (z) \Bigl(
                    (a_1 \boxtimes \cdots \boxtimes a_n)
                    \cap e_z (-\bbnu_{(n)})
                \Bigr)
                \cap e (\oplus_{\smash{(n)}}^* (\mathbb{T}_f))
            \Bigr)
        }
        \\[1ex]
        & =
        Y_f \bigl( a_1 (z_1) \cdots a_n (z_n) \bigr) \ ,
    \end{align*}
    where the second last step follows from \cref{lem-swap}
    by setting $z = 0$ in the lemma,
    which is possible since $\mathbb{T}_f$ is a vector bundle,
    so $e_z (\oplus_{\smash{(n)}}^* (\mathbb{T}_f))$
    does not have negative powers of~$z$.

    The final statement is elementary.
\end{proof}

\section{Modules}

\label{sec-vm}

We present a construction of modules for Joyce vertex algebras from moduli spaces,
which we state in \cref{thm-joyce-vm}.
We also consider an important variant of this construction
in \cref{thm-joyce-vm-tw},
arising from orthosymplectic enumerative geometry,
where the Joyce vertex algebra admits an involution,
and we obtain a twisted module for this vertex algebra.

\subsection{Vertex algebra modules}

\begin{para}[Modules]
    \label{para-vm-def}
    We recall the definition of vertex algebra modules from
    \textcite[\S 5.1]{frenkel-ben-zvi-2004}.

    Let~$(V, 1, D, Y)$ be a $\mathbb{Z}$-graded vertex algebra over~$K$,
    as in \cref{para-va-def}.
    Then a \emph{$V$-module}
    is a $\mathbb{Z}$-graded $K$-vector space~$M$,
    equipped with an operation
    $Y^M (-, z) (-) \colon V \otimes M \to M \llparen z \rrparen$,
    preserving grading with $\deg z = -2$,
    satisfying the following conditions:

    \begin{enumerate}
        \item
            \label{item-vm-unit}
            (\emph{Unit})
            For any $m \in M$, we have
            $Y^M (1, z) \, m = m$.

        \item
            \label{item-vm-assoc}
            (\emph{Associativity})
            For any $a, b \in V$ and $m \in M$, the elements
            \begin{align*}
                Y^M (a, z) \circ Y^M (b, w) (m)
                & \in M \llparen z \rrparen \llparen w \rrparen \ ,
                \\
                Y^M (Y (a, z - w) (b), w) (m)
                & \in M \llparen w \rrparen \llparen z - w \rrparen
            \end{align*}
            are expansions of the same element in
            $M \llbr z, w \rrbr \, [z^{-1}, w^{-1}, (z - w)^{-1}]$.
    \end{enumerate}
\end{para}

\begin{para}[An equivalent definition]
    \label{para-vm-def-alt}
    We can also give an alternative definition
    of a vertex algebra module,
    in the style of \cref{para-va-def-alt}.

    Let~$(V, (X_n)_{n \geq 0})$ be a $\mathbb{Z}$-graded vertex algebra over~$K$,
    in the sense of \cref{para-va-def-alt}.
    Then a \emph{$V$-module} is a $\mathbb{Z}$-graded $K$-vector space~$M$,
    equipped with operations
    \begin{align*}
        X_n^M \colon
        V^{\otimes n} \otimes M
        & \longrightarrow
        M \llbr z_1, \dotsc, z_n \rrbr \, [ z_i^{-1}, (z_i - z_j)^{-1} ] \ ,
        \\
        a_1 \otimes \cdots \otimes a_n \otimes m
        & \longmapsto
        X_n^M (a_1, \dotsc, a_n, m; z_1, \dotsc, z_n) \ ,
    \end{align*}
    preserving grading, where $\deg z_i = -2$,
    and we invert $z_i - z_j$ for $i \neq j$.
    They should satisfy the following properties:

    \begin{enumerate}
        \item
            \label{item-vm-alt-unit}
            (\emph{Unit})
            We have $X_0^M = \mathrm{id}_M$.

        \item
            \label{item-vm-alt-assoc}
            (\emph{Associativity})
            For integers $k, n \geq 0$ and elements
            $a_1, \dotsc, a_n, b_1, \dotsc, b_k \in V$, we have
            \begin{align}
                & X_{n+1}^M \Bigl(
                    X_k (b_1, \dotsc, b_k; w_1, \dotsc, w_k);
                    a_1, \dotsc, a_n, m;
                    z_0, \dotsc, z_n
                \Bigr)
                \notag
                \\[-.5ex]
                \label{eq-vm-alt-assoc-1}
                & \hspace{2em} =
                \iota_{ \{ z_i \}, \{ w_j \} } \,
                X_{n+k}^M \bigl(
                    b_1, \dotsc, b_k,
                    a_1, \dotsc, a_n, m;
                    z_0 + w_1, \dotsc, z_0 + w_k,
                    z_1, \dotsc, z_n
                \bigr) \ ,
                \\[1ex]
                & X_{n}^M \Bigl(
                    a_1, \dotsc, a_n, X_k^M (b_1, \dotsc, b_k, m; w_1, \dotsc, w_k);
                    z_1, \dotsc, z_n
                \Bigr)
                \notag
                \\[-.5ex]
                \label{eq-vm-alt-assoc-2}
                & \hspace{2em} =
                \iota_{ \{ z_i \}, \{ w_j \} } \,
                X_{n+k}^M \bigl(
                    a_1, \dotsc, a_n, b_1, \dotsc, b_k, m;
                    z_1, \dotsc, z_n,
                    w_1, \dotsc, w_k
                \bigr) \ ,
            \end{align}
            where the maps $\iota_{ \{ z_i \}, \{ w_j \} }$
            are defined in \cref{eq-def-iota}.
    \end{enumerate}
    Again, this definition may look more complicated than
    \cref{para-vm-def},
    but it will be more convenient for our purposes,
    such as for defining weaker notions of modules below.

    We adopt the convenient notation
    \begin{equation}
        a_1 (z_1) \cdots a_n (z_n) \, m =
        X^M_n (a_1, \dotsc, a_n, m; z_1, \dotsc, z_n) \ ,
    \end{equation}
    with the same caveats as those below \cref{eq-va-notation}
    for vertex algebras.
\end{para}

\begin{para}[Proof of the equivalence of definitions]
    We essentially follow the proof of
    \textcite[Theorem~3.14]{kim-2011},
    adapting it to the case of modules.

    Given a $V$-module~$M$ in the sense of \cref{para-vm-def},
    define the maps $X_n^M$ by
    \begin{equation}
        \label{eq-vm-proof-1}
        X_n^M (a_1, \dotsc, a_n, m; z_1, \dotsc, z_n) =
        \iota_{z_1, \dotsc, z_n} ^{-1} \bigl(
            Y^M (a_1, z_1) \circ \cdots \circ Y^M (a_n, z_n) (m)
        \bigr) \ ,
    \end{equation}
    where we take the unique preimage under the embedding
    \begin{equation}
        \label{eq-vm-proof-2}
        \iota_{z_1, \dotsc, z_n} \colon
        M \llbr z_1, \dotsc, z_n \rrbr \, [z_i^{-1}, (z_i - z_j)^{-1}]
        \longhookrightarrow
        M \llparen z_1 \rrparen \cdots \llparen z_n \rrparen \ .
    \end{equation}
    To see that such a preimage exists,
    note that for fixed elements $a_1, \dotsc, a_n, m$ as above,
    there exists $N > 0$ such that in the expression
    \begin{equation}
        \label{eq-vm-proof-3}
        \biggl(
            \prod_{1 \leq i < j \leq n} {} (z_i - z_j)
        \biggr)^N \cdot
        Y^M (a_1, z_1) \circ \cdots \circ Y^M (a_n, z_n) (m) \ ,
    \end{equation}
    the order of the operators $Y^M (a_i, z_i)$
    can be permuted freely without changing the result,
    up to an appropriate sign.
    Indeed, this is true when $n = 2$,
    a proof of which can be found in
    \textcite[Remark~5.1.4]{frenkel-ben-zvi-2004},
    and the general case follows from
    swapping two adjacent operators at a time.

    This property implies that the expression
    \cref{eq-vm-proof-3}
    must lie in the intersection of the spaces
    $M \llparen z_{\sigma (1)} \rrparen \cdots \llparen z_{\sigma (n)} \rrparen
    \subset M \llbr z_i^{\pm 1} \rrbr$
    for all $\sigma \in \mathfrak{S}_n$,
    which is $M \llbr z_i \rrbr \, [z_i^{-1}]$.
    Consequently, the element
    $Y^M (a_1, z_1) \circ \cdots \circ Y^M (a_n, z_n) (m)$
    lies in the image of \cref{eq-vm-proof-2}.

    From here, verifying the axioms in \cref{para-vm-def-alt} is straightforward.

    Conversely, given the data $(M, (X_n^M)_{n \geq 0})$, define
    \begin{equation}
        Y^M (a, z) (m) =
        X^M_1 (a, m; z)
    \end{equation}
    for all $a \in V$ and $m \in M$.
    Then this defines a vertex algebra module structure on~$M$,
    where the property \cref{para-vm-def}~\cref{item-vm-unit}
    holds since
    $Y^M (1, z) (m) = X^M_1 (X_0 (1), m; z) = X^M_0 (m) = m$.
    \qed
\end{para}

\begin{para}[Remark]
    In \cref{para-vm-def-alt},
    it is enough to require the operators
    $X_n^M$ for $n = 0, 1, 2$,
    and the higher ones can be recovered from them,
    in a similar way to \cref{para-va-3-terms}.
\end{para}

\begin{para}[Weak modules]
    \label{para-weak-vm}
    Similarly to weak vertex algebras
    introduced in \cref{para-weak-va},
    we define a \emph{weak module}
    over a $\mathbb{Z}$-graded weak vertex algebra $(V, (X_n)_{n \geq 0})$
    to be a $\mathbb{Z}$-graded vector space~$M$,
    equipped with multiplication maps
    \begin{equation*}
        X^M_n \colon
        V^{\otimes n} \otimes M
        \longrightarrow
        M \llparen z_1, \dotsc, z_n \rrparen \ ,
    \end{equation*}
    preserving grading, where $\deg z_i = -2$,
    satisfying the properties
    \cref{para-vm-def-alt}~\crefrange{item-vm-alt-unit}{item-vm-alt-assoc}.

    For example, every weak vertex algebra is a weak module over itself.

    Note that if~$V$ is a vertex algebra,
    then a weak module is a usual module
    if and only if the image of~$X^M_n$ lies in the subspace
    $M \llbr z_1, \dotsc, z_n \rrbr \, [z_i^{-1}, (z_i - z_j)^{-1}]
    \subset M \llparen z_1, \dotsc, z_n \rrparen$.
\end{para}

\begin{para}[Twisted modules]
    \label{para-twisted-modules}
    We introduce \emph{twisted modules}
    for vertex algebras equipped with involutions,
    which are a special class of weak modules.
    They will play an important role
    in the theory of wall-crossing
    for orthosymplectic enumerative invariants.
    Note that these are \emph{different} from
    the notion of twisted modules in, for example,
    \textcite[\S5.6]{frenkel-ben-zvi-2004}.

    Let~$V$ be a $\mathbb{Z}$-graded vertex algebra over~$K$.
    Suppose that~$V$ is equipped with
    a \emph{twisted involution}
    $(-)^\vee \colon V \simto V^\mathrm{op}$,
    where $V^\mathrm{op}$ is the vertex algebra
    with the same underlying graded vector space as~$V$,
    with multiplication
    \begin{equation}
        \label{eq-va-invol}
        a_1 (z_1) \cdots a_n (z_n)
        \text{ in } V^\mathrm{op}
        \quad {=} \quad
        a_1 (-z_1) \cdots a_n (-z_n)
        \text{ in } V \ .
    \end{equation}
    We require that $(-)^\vee$ is an isomorphism of
    $\mathbb{Z}$-graded vertex algebras,
    and $(-)^{\vee \vee} = \mathrm{id}_V$.

    Define a \emph{twisted module} for~$V$
    as a weak module $(M, (X_n^M)_{n \geq 0})$,
    such that the image of each~$X^M_n$ lies in the subspace
    \begin{equation*}
        M \llbr z_1, \dotsc, z_n \rrbr \, [z_i^{-1}, (z_i \pm z_j)^{-1}]
        \subset M \llparen z_1, \dotsc, z_n \rrparen \ ,
    \end{equation*}
    where we invert $z_i \pm z_j$ for $i \neq j$, and such that
    \begin{equation}
        \label{eq-vm-invol}
        a^\vee (z) \, m =
        a (-z) \, m
    \end{equation}
    for $a \in V$ and $m \in M$.

    We will construct examples of twisted modules
    in \cref{subsec-joyce-vm} below.
\end{para}

\begin{para}[Residues]
    Let $V$ be a $\mathbb{Z}$-graded vertex algebra.
    Recall that the quotient $L = V / D (V)$,
    where~$D$ is the translation operator in \cref{para-va-def},
    admits a $\mathbb{Z}$-graded Lie algebra structure given by
    \begin{equation}
        [a, b] = \mathrm{res}_{z = w} \, a (z) \, b (w) \ ,
    \end{equation}
    with sign rules implemented for odd elements
    in the axioms of a Lie algebra.
    A usual module for~$V$ gives rise to a module for~$L$,
    defined by $a \cdot m = \mathrm{res}_{z = 0} \, a (z) \, m$.

    Now, suppose that~$V$ is equipped with an involution $(-)^\vee$
    as in \cref{para-twisted-modules},
    and let~$M$ be a twisted module for~$V$.
    Then the Lie algebra $L = V / D (V)$
    admits an induced involution
    $(-)^\vee \colon L \simto L^\mathrm{op}$,
    where~$L^\mathrm{op}$ is the Lie algebra with the opposite Lie bracket
    $[a, b]^\mathrm{op} = -[a, b]$.
    Moreover, the assignment
    \begin{equation}
        a \heart m = \mathrm{res}_{z = 0} \, a (z) \, m
    \end{equation}
    for $a \in V$ and $m \in M$
    establishes~$M$ as a \emph{twisted module} for~$L$,
    in that we have
    \begin{align}
        a^\vee \heart m
        & =
        -a \heart m \ ,
        \\
        \label{eq-lie-tw-mod}
        a \heart (b \heart m)
        - (-1)^{|a| \cdot |b|} \cdot b \heart (a \heart m)
        & =
        [a, b] \heart m
        - [a^\vee, b] \heart m
    \end{align}
    for homogeneous elements $a, b \in L$ and $m \in M$,
    where $|{-}|$ denotes the degree.
    The identity \cref{eq-lie-tw-mod} can be seen by
    applying the residue theorem to
    $a (z) \, b (w) \, m$, giving
    \begin{equation}
        \mathrm{res}_{z = 0} \, \mathrm{res}_{w = 0} \, a (z) \, b (w) \, m
        = \mathrm{res}_{w = 0} \, \bigl(
            \mathrm{res}_{z = 0} +
            \mathrm{res}_{z = w} +
            \mathrm{res}_{z = -w}
        \bigr) \,
        a (z) \, b (w) \, m \ ,
    \end{equation}
    where each of the four terms corresponds to a term
    in \cref{eq-lie-tw-mod},
    and we use the fact that
    $a (z) \, b (w) \, m = a^\vee (-z) \, b (w) \, m$.

    Note that this structure of a twisted module
    for a Lie algebra also appears in
    motivic Donaldson--Thomas theory
    for orthosymplectic structure groups,
    described in \textcite[\S 5.2.2]{bu-osp-dt}.
    Also, as noted there, when the coefficient ring contains~$1/2$,
    a twisted module for~$L$ is equivalent to a usual module
    for the subalgebra of~$L$ consisting of elements $a \in L$ with
    $a = -a^\vee$, with the action
    $a \cdot m = (1/2) (a \heart m)$.
\end{para}

\subsection{Modules for Joyce vertex algebras}
\label{subsec-joyce-vm}

\begin{para}[The setting]
    \label{para-vm-setting}
    We assume given the data
    $(X, \odot, \oplus, 0, \mathbb{T}_X)$
    as in \cref{para-va-setting},
    We also assume given the extra data
    $(Y, \diamond, \mathbb{T}_Y)$, where

    \begin{itemize}
        \item
            $Y \in \mathsf{hCW}$ is a space,
            regarded as an object of~$\mathcal{T}$ of rank~$0$.

        \item
            $\diamond \colon X \times Y \to Y$
            is a map,
            exhibiting~$Y$ as an $X$-module in~$\mathsf{hCW}$,
            or equivalently, in~$\mathcal{T}$.

        \item
            $\mathbb{T}_Y \in K (Y)$ is a class,
            called the \emph{obstruction theory}.
    \end{itemize}
    For each integer $n \geq 0$, let
    \begin{equation}
        \diamond_{(n)} \colon X^n \times Y \longrightarrow Y
    \end{equation}
    denote the $n$-fold module multiplication map.
    Define the \emph{normal bundle} of~$\diamond_{(n)}$ by
    \begin{equation}
        \bbnu_{\diamond, \smash{(n)}}
        = \diamond_{\smash{(n)}}^* (\mathbb{T}_{Y})
        - \mathbb{T}_{X^n \times Y}
        \in K (X^n \times Y) \ ,
    \end{equation}
    where
    $\mathbb{T}_{X^n \times Y}
    = \mathrm{pr}_{\smash{X^n}}^* (\mathbb{T}_{X^n}) + \mathrm{pr}_Y^* (\mathbb{T}_{Y})$,
    where $\mathrm{pr}_{X^n}$, $\mathrm{pr}_Y$
    are the projections to $X^n$ and~$Y$.
    We further assume the following condition:
    \begin{itemize}
        \item
            For any integer $n \geq 0$, we have
            \begin{equation}
                \label{eq-nu-diamond-weight-cond}
                \bbnu_{\diamond, \smash{(n)}} \in K^\circ (X^n \times Y)
                \quad \text{and} \quad
                (\bbnu_{\diamond, \smash{(n)}})_0 = 0 \ ,
            \end{equation}
            where $(-)_0$ denotes the part of $\mathrm{U} (1)^n$-weight~$0$.
    \end{itemize}
\end{para}

\begin{theorem}
    \label{thm-joyce-vm}
    Let~$X, Y$ be as in \cref{para-vm-setting}. Then the assignment
    \begin{equation}
        \label{eq-def-vm-joyce}
        a_1 (z_1) \cdots a_n (z_n) \, m =
        (\diamond_{(n)})_* \circ
        \exp (z_1 D_1 + \cdots + z_n D_n)
        \Bigl(
            (a_1 \boxtimes \cdots \boxtimes a_n \boxtimes m)
            \cap e_z^{-1} (\bbnu_{\diamond, \smash{(n)}})
        \Bigr)
    \end{equation}
    defines a weak module structure on
    $M = \mathrm{H}_{\bullet + 2 \operatorname{vdim}} (Y; \mathbb{Q})$
    for the weak vertex algebra
    $V = \mathrm{H}_{\bullet + 2 \operatorname{vdim}} (X; \mathbb{Q})$
    in \cref{thm-joyce-va}.

    The product
    \cref{eq-def-vm-joyce}
    has poles along $\lambda (z) = 0$
    for non-zero weights $\lambda \in \mathbb{Z}^n$
    appearing in~$\diamond_{(n)}^* (\mathbb{T}_Y)$.
    In particular, if the class
    $\oplus^* (\mathbb{T}_X) \in K^\circ (X^2)$
    only has $\mathrm{U} (1)^2$-weights
    $(-k, k)$ with $|k| \leq 1$,
    so that~$V$ is a vertex algebra,
    and if the class
    $\diamond_{\smash{(2)}}^* (\mathbb{T}_Y) \in K^\circ (X^2 \times Y)$
    only has $\mathrm{U} (1)^2$-weights
    $(k, \ell)$ with $|k| + |\ell| \leq 1$,
    then $M$~is a module for~$V$.
\end{theorem}

\begin{proof}
    The proof is very similar to that of
    \cref{thm-joyce-va},
    and an analogous computation
    verifies the axioms for a weak module.
    See also \cref{thm-vi} below for a
    more general proof.

    For the final statement on usual modules,
    let $n \geq 0$, and consider the class
    $\diamond_{\smash{(n)}}^* (\mathbb{T}_Y) \in K^\circ (X^n \times Y)$.
    It is enough to show that it only has
    $\mathrm{U} (1)^n$-weights of the form
    $e_i - e_j$
    for $0 \leq i \leq j \leq n$,
    where we set $e_0 = 0$, and for $i = 1, \dotsc, n$,
    we set $e_i = (0, \dotsc, 0, 1, 0, \dotsc, 0)$
    with~$1$ at the $i$-th position.
    But the assumption implies that
    for any $I, J \subset \{ 1, \dotsc, n \}$
    with $I \cap J = \varnothing$,
    writing $k_I = \sum_{i \in I} k_i$ and $k_J = \sum_{j \in J} k_j$,
    we have $|k_I| + |k_J| \leq 1$,
    which implies the desired property.
\end{proof}

\begin{para}[An involution-twisted version]
    \label{para-vm-tw-setting}
    We now discuss a version that often arises
    from moduli spaces involving the orthogonal and symplectic groups,
    which will be important in studying enumerative invariants
    for such moduli spaces.

    Let $(X, Y)$ be as in \cref{para-vm-setting}.
    We further assume given the following data:
    \begin{itemize}
        \item
            An involution $(-)^\vee \colon X \simto X$,
            such that $(-)^{\vee \vee} = \mathrm{id}$ in $\mathsf{hCW}$.
    \end{itemize}
    It should satisfy the following conditions:
    \begin{itemize}
        \item
            $(-)^\vee$ preserves the monoid structure~$\oplus$.
        \item
            $(-)^\vee$ reverses the $\mathrm{BU} (1)$-action~$\odot$,
            in that it is equivariant with respect to the map
            $(-)^{-1} \colon \mathrm{BU} (1) \to \mathrm{BU} (1)$.
        \item
            $(-)^\vee$ preserves the module action~$\diamond$, in that
            ${\diamond} \circ ( (-)^\vee \times \mathrm{id}_Y ) = {\diamond}$.
    \end{itemize}
    Again, these are conditions in $\mathsf{hCW}$,
    and do not require higher coherence.
\end{para}

\begin{example}
    \label{eg-vm-tw}
    Our main source of examples
    of the situation of \cref{para-vm-tw-setting}
    is orthosymplectic enumerative geometry,
    as in \textcite{bu-osp-dt}.

    Let~$\mathcal{X}$ be a \emph{self-dual derived linear moduli stack} over~$\mathbb{C}$,
    that is a derived linear moduli stack~$\mathcal{X}$
    as in \cref{eg-linear-stacks},
    equipped with an involution
    $(-)^\vee \colon \mathcal{X} \simto \mathcal{X}$,
    satisfying conditions similar to those in~\cref{para-vm-tw-setting},
    as in \cite[\S 2.2.5]{bu-osp-dt}.
    Typical examples include the following:

    \begin{itemize}
        \item
            The derived moduli stack of representations
            of a \emph{self-dual quiver}~$Q$,
            possibly with potential,
            where being self-dual means that we are given
            a contravariant involution
            $(-)^\vee \colon Q \simto Q^\mathrm{op}$
            of~$Q$,
            which induces an involution of the stack.
            See~\cref{subsec-quiver} below for details,
            and see also \cite[\S 6.1]{bu-osp-dt}.

        \item
            The derived moduli stack of objects
            in a quasi-abelian subcategory
            $\mathcal{A} \subset \mathsf{Perf} (Z)$,
            where~$Z$ is a smooth projective $\mathbb{C}$-variety,
            such that~$\mathcal{A}$ is stable under the dual functor
            $(-)^\vee = \mathbb{R} \calHom (-, L) [s] \colon
            \mathsf{Perf} (Z) \simto \mathsf{Perf} (Z)^\mathrm{op}$,
            where we choose a line bundle $L \to Z$
            and an integer~$s$.
            See also \cite[\S 6.2]{bu-osp-dt} for this set-up.
    \end{itemize}
    In this case, the derived fixed locus of the involution
    $\mathcal{X}^\mathrm{sd} = \mathcal{X}^{\mathbb{Z}_2}$
    is the moduli stack of \emph{self-dual objects},
    whose points correspond to pairs $(x, \phi)$ with
    $x \in \mathcal{X}$ and $\phi \colon x \simto x^\vee$,
    such that $\phi = \phi^\vee$.
    Such objects can often be interpreted as
    orthogonal or symplectic objects in a linear category,
    where the difference between being orthogonal and symplectic
    is encoded in the higher coherence data
    of the involution on~$\mathcal{X}$ as a $\mathbb{Z}_2$-action.
    Again, see \cite{bu-osp-dt} for details.

    We have an action
    $\diamond \colon \mathcal{X} \times \mathcal{X}^\mathrm{sd}
    \to \mathcal{X}^\mathrm{sd}$
    given by $(x, y) \mapsto x \oplus y \oplus x^\vee$,
    where $x \oplus x^\vee$ is equipped with the obvious
    self-dual structure.

    We now set $X = |\mathcal{X}|$
    and $Y = |\mathcal{X}^\mathrm{sd}|$
    to be the topological realizations,
    and with the involution on~$X$
    and the $X$-action on~$Y$
    given by the involution on~$\mathcal{X}$
    and its action on~$\mathcal{X}^\mathrm{sd}$.
    Let~$\mathbb{T}_X$, $\mathbb{T}_Y$
    be the classes of the tangent complexes
    $\mathbb{T}_\mathcal{X}$, $\mathbb{T}_{\smash{\mathcal{X}^\mathrm{sd}}}$,
    where we note that
    $\mathbb{T}_{\smash{\mathcal{X}^\mathrm{sd}}}$
    is given by the $\mathbb{Z}_2$-invariant part of
    the pullback of~$\mathbb{T}_\mathcal{X}$
    to~$\mathcal{X}^\mathrm{sd}$.

    Then, assuming that~$X$ satisfies
    the conditions in \cref{para-va-setting},
    all the extra conditions in
    \cref{para-vm-setting,para-vm-tw-setting}
    will be automatically satisfied,
    where the weight condition
    \cref{eq-nu-diamond-weight-cond}
    follows from the similar condition~\cref{eq-cond-nu}
    for~$\mathbb{T}_\mathcal{X}$
    and the above description of~$\mathbb{T}_{\smash{\mathcal{X}^\mathrm{sd}}}$.
\end{example}

\begin{example}
    \label{eg-vm-dg-cat}
    We also consider a derived version of
    \cref{eg-vm-tw},
    similar to \cref{eg-dg-cat}.

    Namely, let~$\mathcal{C}$ be a $\mathbb{C}$-linear dg-category
    as in \cref{eg-dg-cat},
    and suppose we are given a contravariant involution
    $(-)^\vee \colon \mathcal{C} \simto \mathcal{C}^\mathrm{op}$,
    inducing a $\mathbb{Z}_2$-action on the derived moduli stack,
    $(-)^\vee \colon \mathcal{X} \simto \mathcal{X}$.
    Examples of~$\mathcal{C}$ include the following:
    \begin{itemize}
        \item
            The derived category of representations
            of a self-dual quiver,
            with the involution induced by
            the self-dual structure of the quiver.

        \item
            The category of perfect complexes
            on a smooth projective $\mathbb{C}$-variety~$Z$,
            with the involution
            $(-)^\vee = \mathbb{R} \calHom (-, L) [s]$
            as in \cref{eg-vm-tw}.
    \end{itemize}
    Again, the derived fixed locus
    $\mathcal{X}^\mathrm{sd} = \mathcal{X}^{\mathbb{Z}_2}$
    is the moduli stack of self-dual objects,
    and we have an action
    $\diamond \colon \mathcal{X} \times \mathcal{X}^\mathrm{sd}
    \to \mathcal{X}^\mathrm{sd}$
    given by $(x, y) \mapsto x \oplus y \oplus x^\vee$.

    Set $X = |\mathcal{X}|$
    and $Y = |\mathcal{X}^\mathrm{sd}|$,
    and let $\mathbb{T}_X$, $\mathbb{T}_Y$
    be the classes of the tangent complexes
    $\mathbb{T}_\mathcal{X}$, $\mathbb{T}_{\smash{\mathcal{X}^\mathrm{sd}}}$.
    Then all the conditions in
    \cref{para-va-setting,para-vm-setting,para-vm-tw-setting}
    are automatically satisfied.
\end{example}

\begin{theorem}
    \label{thm-joyce-vm-tw}
    Let~$X, Y$ be as in \cref{para-vm-tw-setting}.
    Assume the class
    $\oplus^* (\mathbb{T}_X) \in K^\circ (X^2)$
    only has $\mathrm{U} (1)^2$-weights
    $(-k, k)$ with $|k| \leq 1$,
    so that~$V = \mathrm{H}_{\bullet + 2 \operatorname{vdim}} (X; \mathbb{Q})$ is a vertex algebra
    by \cref{thm-joyce-va},
    and the class
    $\diamond_{\smash{(2)}}^* (\mathbb{T}_Y) \in K^\circ (X^2 \times Y)$
    only has $\mathrm{U} (1)^2$-weights
    $(k, \ell)$ with $|k| + |\ell| \leq 2$.

    Then the product~\cref{eq-def-vm-joyce} establishes
    $M = \mathrm{H}_{\bullet + 2 \operatorname{vdim}} (Y; \mathbb{Q})$
    as a twisted module for the vertex algebra~$V$.
\end{theorem}

\begin{proof}
    It follows from \cref{thm-joyce-vm}
    that~$M$ is a weak module for~$V$,
    and the relations
    \crefrange{eq-va-invol}{eq-vm-invol}
    follow from the constructions.
    It remains to verify the restriction on poles,
    which is equivalent to that the class
    $\diamond_{\smash{(n)}}^* (\mathbb{T}_Y) \in K^\circ (X^n \times Y)$
    only has $\mathrm{U} (1)^n$-weights of the form
    $e_i \pm e_j$
    for $0 \leq i \leq j \leq n$,
    with notations as in the proof of
    \cref{thm-joyce-vm}.
    But the assumption implies that
    for any $I, J \subset \{ 1, \dotsc, n \}$
    with $I \cap J = \varnothing$,
    we have $|k_I| + |k_J| \leq 2$,
    again with notations as in the proof of
    \cref{thm-joyce-vm},
    and an elementary argument shows that
    this implies the desired property on $\mathrm{U} (1)^n$-weights.
\end{proof}

\begin{para}[Morphisms of modules]
    \label[theorem]{thm-joyce-vm-mor}
    Finally, we state a result
    analogous to \cref{thm-joyce-va-mor},
    which gives a construction of morphisms
    between modules for Joyce vertex algebras.
\end{para}

\begin{theorem*}
    Let $(X, Y)$ and $(X', Y')$ be as in \cref{para-vm-setting},
    $f \colon X \to X'$ be a map as in \cref{thm-joyce-va-mor},
    and $g \colon Y \to Y'$ a map compatible with
    the actions of\/~$X$ and~$X'$,
    such that the class
    \begin{equation*}
        \mathbb{T}_g = \mathbb{T}_Y - g^* (\mathbb{T}_{Y'})
    \end{equation*}
    is the class of a vector bundle on~$Y$,
    possibly of mixed rank.
    Then the map
    \begin{align*}
        Y_g \colon
        \mathrm{H}_{\bullet + 2 \operatorname{vdim}} (Y; \mathbb{Q})
        & \longrightarrow
        \mathrm{H}_{\bullet + 2 \operatorname{vdim}} (Y'; \mathbb{Q}) \ ,
        \\
        a
        & \longmapsto
        g_* (a \cap e (\mathbb{T}_g))
    \end{align*}
    is compatible with the weak vertex algebra module structures,
    in the sense that
    \begin{equation*}
        Y_g \bigl( a_1 (z_1) \cdots a_n (z_n) \, m \bigr) =
        Y_f (a_1) (z_1) \cdots Y_f (a_n) (z_n) \, Y_g (m)
    \end{equation*}
    for any $a_1, \dotsc, a_n \in \mathrm{H}_\bullet (X; \mathbb{Q})$
    and $m \in \mathrm{H}_\bullet (Y; \mathbb{Q})$.
    In particular, if\/ $X = X'$ and $f = \mathrm{id}_X$,
    then $g$ is a homomorphism of weak modules.
\end{theorem*}

The proof is very similar to that of
\cref{thm-joyce-va-mor},
and we omit it here.
See \cref{thm-vi} below for a proof
of a more general result.

\section{Vertex induction}

\label{sec-vi}

This section presents the main construction of this paper,
the vertex induction,
which unifies and simultaneously generalizes
\cref{thm-joyce-va,thm-joyce-va-mor,thm-joyce-vm,thm-joyce-vm-tw,thm-joyce-vm-mor}.

We start by introducing a \emph{paracategory of vertex spaces},
in which vertex algebras are algebra objects,
which we prove in \cref{thm-va-func-def}.
This construction is somewhat similar to that in
\textcite{borcherds-1998},
although not precisely the same.
This will then provide a functorial approach to
reformulating the constructions in \crefrange{sec-va}{sec-vm},
and we will use this to describe the vertex induction
in \cref{thm-vi}.
We specialize to the case of algebraic stacks
in \cref{thm-vi-stack},
which is our motivating case.

\subsection{Vertex algebras as algebra objects}

\begin{para}[Paracategories]
    A \emph{paracategory} is, roughly speaking,
    a `category with partially defined composition'.
    The following definition is adapted from
    \textcite[\S 2]{hermida-mateus-2003-paracategories}.
    Compare the notion of a \emph{relaxed multilinear category}
    in~\textcite[\S 4]{borcherds-1998}.

    A \emph{paracategory}~$\mathcal{C}$ is a directed graph
    $s, t \colon \mathcal{C}_1 \rightrightarrows \mathcal{C}_0$,
    together with a partially defined function
    $\circ_n \colon \mathcal{C}_n \dashrightarrow \mathcal{C}_1$ for all $n \geq 0$,
    meaning a function from a subset
    of the $n$-fold fibre product set
    $\mathcal{C}_n = \mathcal{C}_1 \times_{s, \mathcal{C}_0, t} \cdots \times_{s, \mathcal{C}_0, t} \mathcal{C}_1$
    to~$\mathcal{C}_1$,
    satisfying the following axioms:

    \begin{enumerate}
        \item
            The map $\circ_0 \colon \mathcal{C}_0 \to \mathcal{C}_1$ is fully defined.
            That is, \emph{identity morphisms} always exist.
        \item
            The map $\circ_1 \colon \mathcal{C}_1 \to \mathcal{C}_1$
            is fully defined, and is the identity map.
        \item
            For $(x_1, \dotsc, x_n) \in \mathcal{C}_{k_1} \times_{\mathcal{C}_0} \cdots \times_{\mathcal{C}_0} \mathcal{C}_{k_n}$,
            we have
            \begin{equation}
                \circ_{k_1 + \cdots + k_n} (x_1, \dotsc, x_n) =
                \circ_n (\circ_{k_1} x_1, \dotsc, \circ_{k_n} x_n) \ ,
            \end{equation}
            in the sense that whenever one side of the equality is defined,
            so is the other, and both sides are equal.
    \end{enumerate}
    Any (small) category can be seen as a paracategory in the obvious way.

    Note that the axioms imply that composition with an isomorphism is always defined,
    since for example, if $x \in \mathcal{C}_1$ is an isomorphism,
    and $y \in \mathcal{C}_1$ is such that $s (y) = t (x)$, then
    $y = y \circ (x \circ x^{-1}) = (y \circ x) \circ x^{-1}$,
    showing that $y \circ x$ is defined.

    A \emph{functor of paracategories} is a map of directed graphs $f \colon \mathcal{C} \to \mathcal{D}$,
    such that for any $x \in \mathcal{C}_k$, if $\circ_k x$ is defined,
    then $\circ_k f (x)$ is defined and is equal to $f (\circ_k x)$.
    Note that such functors can always be composed.

    A \emph{natural isomorphism} of functors $f, g \colon \mathcal{C} \to \mathcal{D}$
    consists of an isomorphism $\eta_x \colon f (x) \simto g (x)$ for any $x \in \mathcal{C}_0$,
    such that the naturality squares commute.
    Natural isomorphisms can always be composed.
\end{para}

\begin{para}[Monoidal paracategories]
    A (\emph{symmetric}) \emph{monoidal structure} on a paracategory~$\mathcal{C}$
    consists of an object $1 \in \mathcal{C}_0$ and a functor
    $\otimes \colon \mathcal{C} \times \mathcal{C} \to \mathcal{C}$,
    equipped with natural isomorphisms satisfying the usual axioms.

    One can also define (\emph{symmetric}) \emph{monoidal functors}
    between (symmetric) monoidal paracategories in the usual way.

    For a (coloured) operad~$\mathcal{O}$
    and a symmetric monoidal paracategory~$\mathcal{C}$,
    an \emph{$\mathcal{O}$-algebra} in~$\mathcal{C}$ is a symmetric monoidal functor $\mathcal{O}^\otimes \to \mathcal{C}$,
    where $\mathcal{O}^\otimes$ is the underlying
    symmetric monoidal category of~$\mathcal{O}$,
    as in \textcite[Remark~2.2.4.6]{lurie-ha}.

    For example, an \emph{associative algebra} in~$\mathcal{C}$
    can be described as an object $x \in \mathcal{C}$,
    with $n$-ary multiplication maps $m_n \colon x^{\otimes n} \to x$
    for all $n \geq 0$, satisfying the usual axioms.
    A \emph{commutative algebra} is an associative algebra
    such that the maps $m_n$ are equivariant under permutations.
\end{para}

\begin{para}[Vertex spaces]
    \label{para-vertex-spaces}
    We now define a symmetric monoidal paracategory
    of \emph{vertex spaces},
    where vertex algebras are algebra objects.

    Define a \emph{vertex space} over a field~$K$
    as a triple $(\Lambda, V, \tau)$,
    consisting of the following:
    \begin{itemize}
        \item
            A finite-dimensional $K$-vector space~$\Lambda$,
            whose dimension is called the \emph{rank}
            of the vertex space.
        \item
            A $\mathbb{Z}$-graded $K$-vector space~$V$,
            thought of as a space of \emph{fields} on~$\Lambda$.
        \item
            A $K$-linear map
            $\tau (z) \colon V \to V \llbr \Lambda \rrbr$,
            called the \emph{translation operator},
            where~$z$ is a set of coordinates for~$\Lambda$,
            and $\deg \Lambda^\vee = -2$,
            thought of as translating fields along~$z$.
            It should satisfy
            \begin{align}
                \label{eq-tau-unit}
                \tau (0) & = \mathrm{id}_V \ ,
                \\
                \label{eq-tau-add}
                \tau (z) \circ \tau (w) & = \tau (z + w) \ .
            \end{align}
            In other words, $\tau$ is a representation
            of the formal group $\hat{\mathbb{G}}_\mathrm{a}^r$
            of the additive group of~$\Lambda$.
    \end{itemize}
    We often write $V$ for the vertex space,
    and write $\Lambda_V$ and~$\tau_V$ for $\Lambda$ and~$\tau$.
\end{para}

\begin{para}[Maps of vertex spaces]
    \label{para-map-vertex-spaces}
    A \emph{map of vertex spaces} $f \colon V \to V'$
    consists of a $K$-linear map $f^{\smash{\sharp}} \colon \Lambda_{V'} \to \Lambda_V$,
    and a $K$-linear map
    \begin{equation*}
        f (z) \colon V \longrightarrow V' \llparen \Lambda_V \rrparen \ ,
    \end{equation*}
    where $\deg \Lambda_V^\vee = -2$, such that
    \begin{align}
        \label{eq-fs-map-1}
        f (z) \circ \tau_V (w) & =
        \iota_{z, w} \, f (z + w) \ ,
        \\
        \label{eq-fs-map-2}
        \tau_{V'} (z') \circ f (z) & =
        \iota_{z', z} \, f (f^{\smash{\sharp}} (z') + z) \ ,
    \end{align}
    for $z, w \in \Lambda_V$ and $z' \in \Lambda_{V'}$.
    Here, $f (z)$ can be thought of as translating fields
    along the vector~$z$, then pulling them back to~$\Lambda_{V'}$ along~$f^\sharp$.

    The \emph{identity map} $\mathrm{id} \colon V \to V$
    is given by $\mathrm{id}^\sharp = \mathrm{id}_{\Lambda_V}$ and $\mathrm{id} (z) = \tau_V (z)$.

    Note that the relation~\cref{eq-fs-map-1} might tempt one
    think that $f (z)$ is determined by $f (0)$ via $f (z) = f (0) \circ \tau_V (z)$.
    However, this is not the case, as $f (z)$ is often singular at $z = 0$,
    in which case $f (0)$ is undefined.
\end{para}

\begin{para}[Composition]
    \label{para-composition-vertex-spaces}
    Maps of vertex spaces \emph{cannot} always be composed.
    Given a diagram
    \begin{equation}
        \begin{tikzcd}
            V \ar[r, "f"']
            \ar[rr, bend left=20, shift left=1, "h"]
            & V' \ar[r, "g"']
            & V''
        \end{tikzcd}
    \end{equation}
    of maps of vertex spaces,
    we say that $h$ is the \emph{composition} of~$f$ and~$g$,
    written $h = g \circ f$,
    if $h^\sharp = f^\sharp \circ g^\sharp$, and we have
    \begin{equation}
        g (z') \circ f (z) =
        \iota_{z', z} \, h (f^\sharp (z') + z)
    \end{equation}
    for $z \in \Lambda_V$ and $z' \in \Lambda_{V'}$.
    Note that~$h$ is unique if it exists,
    and always exists if~$f$ or~$g$ is the identity.
    Composition is associative and unital
    whenever all involved compositions are defined.

    Moreover, this generalizes to define $n$-fold compositions of maps of vertex spaces
    for any $n \geq 0$. This defines the structure of a paracategory.

    We denote by $\mathsf{VS}_K$ the paracategory of vertex spaces over~$K$.
\end{para}

\begin{para}[The monoidal structure]
    Given vertex spaces $V$ and~$V'$,
    we define a vertex space structure on $V \otimes V'$ by
    \begin{align}
        \Lambda_{V \otimes V'} & = \Lambda_V \oplus \Lambda_{V'} \ ,
        \\
        \tau_{V \otimes V'} (z, z') & = \tau_V (z) \otimes \tau_{V'} (z') \ .
    \end{align}

    The tensor product defines a symmetric monoidal structure on~$\mathsf{VS}_K$,
    and has a unit given by the vertex space $K$, with $\Lambda_K = 0$ and $\tau_K = \mathrm{id}$.
    We include a sign rule when identifying $V \otimes V'$ with $V' \otimes V$
    on the odd graded pieces.
\end{para}

\begin{theorem}
    \label{thm-va-func-def}
    A $\mathbb{Z}$-graded weak vertex algebra over a field\/~$K$,
    in the sense of\/ \cref{para-weak-va},
    is equivalently a commutative algebra object in~$\mathsf{VS}_K$ of rank~$1$.
    Moreover, a $\mathbb{Z}$-graded weak module for such an algebra,
    in the sense of\/ \cref{para-weak-vm},
    is equivalently a module object in~$\mathsf{VS}_K$ of rank~$0$.
\end{theorem}

\begin{proof}
    These follow from the definitions.
\end{proof}

\subsection{Vertex induction}
\label{subsec-vi}

\begin{para}
    We now introduce the vertex induction,
    which is a simultaneous generalization of
    the constructions of
    Joyce vertex algebras and modules
    in \crefrange{sec-va}{sec-vm}.

    As mentioned in the introduction,
    we expect this construction to be useful for enumerative geometry,
    especially for formulating wall-crossing formulae
    for enumerative invariants of non-linear moduli stacks.
\end{para}

\begin{para}[The paracategory~\texorpdfstring{$\mathcal{T}'$}{T'}]
    \label{para-cat-t-prime}
    To help with stating our main construction in
    \cref{thm-vi},
    we define an auxiliary paracategory~$\mathcal{T}'$ as follows.
    Recall the category~$\mathcal{T}$ from
    \cref{para-cat-t}.
    \begin{itemize}
        \item
            An object of~$\mathcal{T}'$
            is a quadruple $(T_X, X, \odot_X, \mathbb{T}_X)$,
            where $(T_X, X, \odot_X) \in \mathcal{T}$,
            and~$\mathbb{T}_X \in K^\circ (X)$ is a class of weight~$0$.

        \item
            A morphism
            $f \colon (T_X, X, \odot_X, \mathbb{T}_X) \to (T_Y, Y, \odot_Y, \mathbb{T}_Y)$
            in~$\mathcal{T}'$ is a morphism
            $f \colon (T_X, X, \odot_X) \to (T_Y, Y, \odot_Y)$ in~$\mathcal{T}$,
            such that the normal bundle
            \begin{equation*}
                \bbnu_f = f^* (\mathbb{T}_Y) - \mathbb{T}_X
            \end{equation*}
            satisfies the following conditions:
            \begin{enumerate}[label=(\alph*)]
                \item
                    \label{item-nu-k-circ}
                    We have $\bbnu_f \in K^\circ (X)$.
                    Equivalently, we have
                    $f^* (\mathbb{T}_Y) \in K^\circ (X)$.

                \item
                    The weight~$0$ part $(-\bbnu_f)_0$
                    is the class of a vector bundle on~$X$.
            \end{enumerate}
    \end{itemize}
    Note that composition in~$\mathcal{T}'$ is not fully defined,
    since these conditions are not preserved by composition in general:
    For example, consider the composition
    $X \to Y \to Z$, where $X = Y = Z = \mathrm{BU} (1)$,
    the maps are identity maps,
    $T_X = \mathrm{U} (1)$ with the tautological action
    of $\mathrm{B} T_X$ on~$X$,
    $T_Y = T_Z = \{ 1 \}$,
    $\mathbb{T}_X = \mathbb{T}_Y = a$,
    and $\mathbb{T}_Z = 0$
    for some $a \in K (X) \setminus K^\circ (X)$.
\end{para}

\begin{theorem}
    \label{thm-vi}
    There is a functor of paracategories
    \begin{align*}
        V \colon \mathcal{T}'
        & \longrightarrow \mathsf{VS}_\mathbb{Q} \ ,
        \\
        X
        & \longmapsto
        (\Lambda_X = \Lambda_{T_X} \ , \
        V_X = \mathrm{H}_{\bullet + 2 \operatorname{vdim}} (X; \mathbb{Q}) \ , \
        \tau_X ) \ ,
        \\
        \smash{(X \overset{f}{\to} Y)}
        & \longmapsto
        f_* \circ
        \tau_X (z) \circ
        \bigl( (-) \cap e_z^{-1} (\bbnu_f) \bigr) \ ,
    \end{align*}
    where the translation operator~$\tau_X$
    is defined in \cref{para-translation-op}.
\end{theorem}

\begin{proof}
    The proof is very similar to that of
    \cref{thm-joyce-va},
    and the key is in proving that~$V$ preserves composition.
    Let $X \overset{\smash[t]{f}}{\to} Y \overset{\smash[t]{g}}{\to} Z$
    be composable morphisms in~$\mathcal{T}'$.
    Then for any $a \in V_X$, we have
    \begin{align*}
        &
        g_* \circ \tau_Y (z) \Bigl[
            f_* \circ \tau_X (w)
            \Bigl(
                a \cap e_w^{-1} (\bbnu_f)
            \Bigr)
            \cap e_z^{-1} (\bbnu_g)
        \Bigr]
        \\
        = {}
        &
        g_* \circ f_* \circ \tau_X (f^\sharp (z))
        \Bigl[
            \tau_X (w)
            \Bigl(
                a \cap e_w^{-1} (\bbnu_f)
            \Bigr)
            \cap e_z^{-1} (f^* (\bbnu_g))
        \Bigr]
        \\
        = {}
        &
        g_* \circ f_* \circ \tau_X (f^\sharp (z)) \circ \tau_X (w)
        \Bigl(
            a \cap e_w^{-1} (\bbnu_f)
            \cap e_{f^{\smash{\sharp}} (z) + w}^{-1} (f^* (\bbnu_g))
        \Bigr)
        \\
        = {}
        &
        (g \circ f)_* \circ \tau_X (f^\sharp (z) + w)
        \Bigl(
            a \cap e_{f^{\smash{\sharp}} (z) + w}^{-1} (\bbnu_{g \circ f})
        \Bigr) \ ,
    \end{align*}
    where the third line uses \cref{lem-swap},
    and the last step uses the fact that
    $\bbnu_f$ has weight~$0$ with respect to~$T_Y$,
    and the fact that
    $\bbnu_{g \circ f} = \bbnu_f + f^* (\bbnu_g)$.
\end{proof}

This generalizes the main constructions in \crefrange{sec-va}{sec-vm},
in the sense that the associativity and commutativity axioms
of vertex algebras and modules
are now encoded as the functoriality of the functor~$V$,
in view of \cref{thm-va-func-def}.
See also \cref{para-vi-stack-appl} below
for more motivations of this construction from enumerative geometry.

\subsection{Vertex induction for stacks}
\label{subsec-vi-stacks}

\begin{para}
    \label{para-cl-setup}
    We now discuss the vertex induction
    applied to algebraic stacks over~$\mathbb{C}$,
    which is the case that we expect to be useful for enumerative geometry.

    Throughout, let~$\mathcal{X}$ be a
    derived algebraic stack over~$\mathbb{C}$,
    locally finitely presented,
    such that its classical truncation~$\mathcal{X}_\mathrm{cl}$
    is a classical algebraic $1$-stack over~$\mathbb{C}$,
    quasi-separated, with affine stabilizers and separated inertia.
\end{para}

\begin{para}
    We use the formalism of the \emph{component lattice}
    of an algebraic stack,
    following \textcite{epsilon-i}.
    The component lattice is a combinatorial structure
    attached to any algebraic stack,
    and is key to generalizing existing results
    in enumerative geometry for linear moduli stacks
    to non-linear moduli stacks.
    It provides the indexing set
    for decomposition-type theorems,
    or for terms in wall-crossing formulae, etc.
    See \cite{epsilon-ii,epsilon-iii,bu-davison-ibanez-nunez-kinjo-padurariu}
    for these applications.
\end{para}

\begin{para}[The component lattice]
    \label{para-component-lattice}
    For a finite rank free $\mathbb{Z}$-module~$\Lambda$,
    let $T_\Lambda = \operatorname{Spec} (\mathbb{C} \Lambda^\vee)
    \simeq \mathbb{G}_\mathrm{m}^{\operatorname{rank} \Lambda}$
    be the torus with coweight lattice~$\Lambda$,
    where $\mathbb{C} \Lambda^\vee$ is the group algebra of~$\Lambda^\vee$.
    Define the \emph{stack of\/ $\Lambda^\vee$-graded points} of~$\mathcal{X}$
    as the derived mapping stack
    \begin{equation*}
        \mathrm{Grad}^\Lambda (\mathcal{X}) =
        \calMap (* / T_\Lambda, \mathcal{X}) \ .
    \end{equation*}
    It is also a derived algebraic stack
    satisfying the conditions in \cref{para-cl-setup},
    and this construction is contravariant in~$\Lambda$.
    We also write
    $\mathrm{Grad} (\mathcal{X}) = \mathrm{Grad}^{\mathbb{Z}} (\mathcal{X})$.

    The \emph{component lattice} of~$\mathcal{X}$
    is the functor
    \begin{align*}
        \mathrm{CL} (\mathcal{X}) \colon
        \mathsf{Lat} (\mathbb{Z})^\mathrm{op}
        & \longrightarrow
        \mathsf{Set} \ ,
        \\
        \Lambda
        & \longmapsto
        \uppi_0
        ( \mathrm{Grad}^\Lambda (\mathcal{X}) ) \ ,
    \end{align*}
    where~$\mathsf{Lat} (\mathbb{Z})$ is the category of
    finite rank free $\mathbb{Z}$-modules,
    and $\uppi_0 (-)$ denotes the set of connected components.
    This does not depend on the derived structure,
    in that we have $\mathrm{CL} (\mathcal{X}) \simeq \mathrm{CL} (\mathcal{X}_\mathrm{cl})$.

    A presheaf on~$\mathsf{Lat} (\mathbb{Z})$
    is also called a \emph{formal lattice},
    so that $\mathrm{CL} (\mathcal{X})$ is a formal lattice.
    Its \emph{underlying set} is the set
    $| \mathrm{CL} (\mathcal{X}) | = \mathrm{CL} (\mathcal{X}) (\mathbb{Z})
    \simeq \uppi_0 ( \mathrm{Grad} (\mathcal{X}) )$.

    For example, if $\mathcal{X} = V / G$
    is a quotient stack,
    where~$G$ is a linear algebraic group over~$\mathbb{C}$
    and~$V$ is a $G$-representation,
    then the component lattice $\mathrm{CL} (\mathcal{X})$
    is the quotient formal lattice $\Lambda_T / W$, where
    $\Lambda_T = \mathrm{Hom} (\mathbb{G}_\mathrm{m}, T)$
    is the coweight lattice of the maximal torus $T \subset G$,
    and $W = \mathrm{N}_G (T) / \mathrm{Z}_G (T)$
    is the Weyl group.
    In general, $\mathrm{CL} (\mathcal{X})$
    is usually glued from copies of~$\mathbb{Z}^n$ for various~$n$,
    along their automorphisms and maps between them.
    See \cite{epsilon-i} for details and more examples.
\end{para}

\begin{para}[Special faces]
    \label{para-special-faces}
    As in \cite{epsilon-i},
    define a \emph{face} of~$\mathcal{X}$
    as a morphism of formal lattices
    $\alpha \colon \Lambda \to \mathrm{CL} (\mathcal{X})$
    for some $\Lambda \in \mathsf{Lat} (\mathbb{Z})$,
    regarded as a formal lattice via the Yoneda embedding.
    Such faces naturally form a category
    $\mathsf{Face} (\mathcal{X})$.
    For a face $\alpha \in \mathsf{Face} (\mathcal{X})$, write
    \begin{equation*}
        \mathcal{X}_\alpha \subset \mathrm{Grad}^\Lambda (\mathcal{X})
    \end{equation*}
    for the connected component corresponding to~$\alpha$.
    For a morphism of faces $f \colon \alpha \to \beta$,
    there is an induced morphism
    $f^\circ \colon \mathcal{X}_\beta \to \mathcal{X}_\alpha$,
    giving a functor
    \begin{equation}
        \label{eq-face-functor}
        \mathcal{X}_{(-)} \colon
        \mathsf{Face} (\mathcal{X})^\mathrm{op}
        \longrightarrow \mathsf{dSt}_\mathbb{C} \ ,
    \end{equation}
    where $\mathsf{dSt}_\mathbb{C}$ is the $\infty$-category of
    derived stacks over~$\mathbb{C}$.

    A \emph{special face} of~$\mathcal{X}$
    is a face~$\alpha$ such that one cannot enlarge~$\alpha$
    without changing~$\mathcal{X}_\alpha$.
    More precisely, $\alpha$ is special
    if for any morphism of faces $f \colon \alpha \to \beta$
    such that~$f^\circ$ is an isomorphism,
    it admits a retraction $g \colon \beta \to \alpha$,
    so that $g \circ f = \mathrm{id}_\alpha$.
    (Note that \cite{epsilon-i} uses rationalized faces,
    but the integral version defined here also works,
    which we use for convenience of presentation,
    and taking the rationalization gives an equivalence
    of the two versions.)

    For example, if $\mathcal{X} = V / G$
    as in \cref{para-component-lattice},
    then the special faces of~$\mathcal{X}$
    are precisely the intersections of hyperplanes in~$\Lambda_T$
    dual to weights of~$V$ and roots of~$G$ in
    $\Lambda^T = (\Lambda_T)^\vee$.
    For such a face
    $\alpha \colon \Lambda \hookrightarrow \Lambda_T
    \twoheadrightarrow \Lambda_T / W$,
    the stack~$\mathcal{X}_\alpha \simeq V^\alpha / L_\alpha$
    is the quotient of the fixed locus $V^\alpha \subset V$
    of the corresponding subtorus $T_\alpha \subset T$
    by the Levi subgroup
    $L_\alpha = \mathrm{Z}_G (T_\alpha) \subset G$.

    Let $\mathsf{Face}^\mathrm{sp} (\mathcal{X}) \subset \mathsf{Face} (\mathcal{X})$
    be the full subcategory of special faces.
    There is a \emph{special closure} functor
    \begin{equation*}
        (-)^\mathrm{sp} \colon \mathsf{Face} (\mathcal{X}) \longrightarrow \mathsf{Face}^\mathrm{sp} (\mathcal{X}) \ ,
    \end{equation*}
    which is left adjoint to the inclusion.
    The functor \cref{eq-face-functor}
    factors through this functor,
    and in particular, we have
    $\mathcal{X}_{\alpha^\mathrm{sp}} \simeq \mathcal{X}_\alpha$
    for any $\alpha \in \mathsf{Face} (\mathcal{X})$,
    so~$\alpha^\mathrm{sp}$ is a canonical replacement of~$\alpha$
    without changing the stack~$\mathcal{X}_\alpha$.

    By the \emph{finiteness theorem}
    \cite[Theorem~6.2.3]{epsilon-i},
    if~$\mathcal{X}$ is quasi-compact and satisfies a very mild condition
    called \emph{quasi-compact graded points},
    then~$\mathcal{X}$ has only finitely many special faces.
\end{para}

\begin{para}[Tangent weights]
    \label{para-tangent-weights}
    For any face
    $\alpha \colon \Lambda \to \mathrm{CL} (\mathcal{X})$,
    write
    $\alpha^\star (\mathbb{T}_\mathcal{X})
    = \mathrm{tot}_\alpha^* (\mathbb{T}_\mathcal{X})
    \in \mathsf{Perf} (\mathcal{X}_\alpha)$,
    where $\mathbb{T}_\mathcal{X}$ is the tangent complex of~$\mathcal{X}$,
    and $\mathrm{tot}_\alpha \colon \mathcal{X}_\alpha \to \mathcal{X}$
    is the forgetful morphism.
    The canonical $* / T_\Lambda$-action on~$\mathcal{X}_\alpha$
    induces a $\Lambda^\vee$-grading of
    $\alpha^\star (\mathbb{T}_\mathcal{X})$,
    and we have
    \begin{equation}
        \label{eq-tangent-x-alpha}
        \mathbb{T}_{\mathcal{X}_\alpha} \simeq
        \alpha^\star (\mathbb{T}_\mathcal{X})_0 \ ,
    \end{equation}
    where $(-)_0$ denotes the degree~$0$ part
    with respect to the $\Lambda^\vee$-grading.

    Define the set of \emph{tangent weights} of~$\mathcal{X}$
    on~$\alpha$ by
    \begin{equation}
        \mathrm{wt} (\mathcal{X}, \alpha) =
        \{ \lambda \in \Lambda^\vee \mid \alpha^\star (\mathbb{T}_\mathcal{X})_\lambda \neq 0 \} \ .
    \end{equation}
    We say that~$\mathcal{X}$ has \emph{finite tangent weights}
    if this set is finite for all faces~$\alpha$,
    or equivalently, for all special faces~$\alpha$.
\end{para}

\begin{para}[A functor to \texorpdfstring{$\mathcal{T}'$}{T'}]
    \label{para-face-to-t-circ}
    For any stack~$\mathcal{X}$
    as in \cref{para-cl-setup},
    with finite tangent weights
    as in \cref{para-tangent-weights},
    and any face $\alpha \colon \Lambda \to \mathrm{CL} (\mathcal{X})$,
    there is an associated object
    \begin{equation*}
        (|T_\Lambda|, |\mathcal{X}_\alpha|, \odot, \mathbb{T}_{\mathcal{X}_\alpha}) \in \mathcal{T}' \ ,
    \end{equation*}
    where
    $|T_\Lambda| \simeq (\mathbb{C}^\times)^n \simeq \mathrm{U} (1)^n$
    in $\mathsf{hCW}$ for $n = \operatorname{rank} \Lambda$,
    and $\odot$ is induced by the canonical
    $* / T_\Lambda$-action on~$\mathcal{X}_\alpha$.
    This defines a functor
    \begin{equation*}
        |\mathcal{X}_{(-)}| \colon
        \mathsf{Face} (\mathcal{X})^\mathrm{op}
        \longrightarrow \mathcal{T}' \ ,
    \end{equation*}
    which factors through the special closure functor
    $(-)^\mathrm{sp} \colon \mathsf{Face} (\mathcal{X}) \to \mathsf{Face}^\mathrm{sp} (\mathcal{X})$.
\end{para}

\begin{theorem}
    \label{thm-vi-stack}
    Let $\mathcal{X}$ be a derived algebraic stack over~$\mathbb{C}$
    as in \cref{para-cl-setup},
    and assume that it has finite tangent weights
    as in \cref{para-tangent-weights}.
    Then there is a functor
    \begin{align*}
        \mathsf{Face}^\mathrm{sp} (\mathcal{X})^\mathrm{op}
        & \longrightarrow \mathsf{VS}_\mathbb{Q} \ ,
        \\
        \alpha
        & \longmapsto
        (\Lambda \ , \
        V_\alpha =
        \mathrm{H}_{\bullet + 2 \operatorname{vdim}} (\mathcal{X}_\alpha; \mathbb{Q}) \ , \
        \tau_\alpha) \ ,
        \\
        \smash{(\alpha \overset{f}{\to} \beta)}
        & \longmapsto
        (f^\circ)_* \circ
        \tau_{\beta} (z) \circ
        \bigl( (-) \cap e_z^{-1} (\bbnu_{f^\circ}) \bigr) \ ,
    \end{align*}
    where~$\Lambda$ is the underlying $\mathbb{Z}$-lattice of~$\alpha$,
    $\tau_\alpha = \tau_{|\mathcal{X}_\alpha|}$
    is the translation operator in \cref{para-translation-op},
    $f^\circ$ is defined in \cref{para-special-faces},
    and
    $\bbnu_{f^\circ} = (f^\circ)^* (\mathbb{T}_{\mathcal{X}_\alpha}) - \mathbb{T}_{\mathcal{X}_\beta}
    \in K^\circ (|\mathcal{X}_\beta|)$.
\end{theorem}

\begin{proof}
    This follows from \cref{thm-vi}
    and \cref{para-face-to-t-circ}.
\end{proof}

\begin{para}
    \label{para-vi-stack-appl}
    We expect that the maps of vertex spaces associated to~$f$
    in \cref{thm-vi-stack},
    or the \emph{vertex induction maps},
    together with their residues along their poles in the~$z$ variables,
    will be useful in studying the structure of enumerative invariants.
    In particular, we expect that the wall-crossing formulae
    for a conjectural version of Joyce's homological invariants,
    discussed in \cref{para-intro-gen},
    should be expressed using these maps,
    and should have the same form as those for
    intrinsic Donaldson--Thomas invariants
    to appear in \cite{epsilon-iii}.
\end{para}

\section{Variants}

\subsection{A real version}
\label{subsec-real-va}

\begin{para}
    We introduce a generalization of the vertex induction,
    where we allow the obstruction theory
    to be an oriented real $K$-theory class
    instead of a complex $K$-theory class.
    This will recover the constructions in
    \crefrange{sec-va}{sec-vi}
    if the oriented real $K$-theory class comes from a complex $K$-theory class.

    We expect that this construction will be important in
    studying \emph{DT\/$4$ invariants}
    for general $(-2)$-shifted symplectic stacks,
    especially in the non-linear cases,
    where the stack is not the moduli stack of objects in a linear category.
    The linear case is discussed in
    \textcite[\S 4.4]{gross-joyce-tanaka-2022}.
\end{para}

\begin{para}[Real \texorpdfstring{$K$}{K}-theory]
    \label{para-ko}
    For a space $X \in \mathsf{hCW}$,
    define the \emph{real topological $K$-theory} of~$X$
    as the commutative ring
    \begin{equation}
        \mathit{KO} (X) =
        \mathsf{hCW} (X, \mathrm{BO} \times \mathbb{Z}) \ .
    \end{equation}
    There is a forgetful map
    $\mathit{KO} (X) \to K (X)$
    given by complexification.

    Given a class $E \in \mathit{KO} (X)$,
    an \emph{orientation} of~$E$ is a trivialization
    of the $\mathbb{Z}_2$-bundle on~$X$
    classified by the composition
    $X \to \mathrm{BO} \times \mathbb{Z} \to \mathrm{B} \mathbb{Z}_2$,
    with the second map induced by the projection
    $\det \colon \mathrm{O} \to \mathbb{Z}_2$
    on each component.
    Such an orientation exists if and only if the class
    $w_1 (E) \in \mathrm{H}^1 (X; \mathbb{Z}_2)$
    is zero, in which case the orientations
    form an $\mathrm{H}^0 (X; \mathbb{Z}_2)$-torsor.

    Recall the category~$\mathcal{T}$ from \cref{para-cat-t}.
    For an object $(T, X, \odot) \in \mathcal{T}$,
    define a subset
    \begin{equation*}
        \mathit{KO}^\circ (X) \subset \mathit{KO} (X)
    \end{equation*}
    of elements that map to the subgroup
    $K^\circ (X) \subset K (X)$
    defined in \cref{para-k-circ}.
\end{para}

\begin{para}[The square root Euler class]
    For a space~$X$ and an $\mathrm{SO} (n)$-bundle $E$ on~$X$,
    recall the \emph{square root Euler class}
    \begin{equation*}
        \sqrte (E) \in \mathrm{H}^{n} (X; \mathbb{Z}) \ ,
    \end{equation*}
    which is the Euler class of the underlying
    oriented real vector bundle,
    and is only non-zero if~$n$ is even.
    It satisfies
    $\sqrte (E)^2 = (-1)^n \cdot e (E)$,
    where $e (E) = c_n (E)$ is the top Chern class
    of the associated complex vector bundle.
\end{para}

\begin{para}[The equivariant square root Euler class]
    \label{para-sqrt-equivar-euler}
    Let~$(T, X, \odot) \in \mathcal{T}$,
    with~$X$ connected.
    Let $E \in \mathit{KO}^\circ (X)$ be a class
    such that the weight~$0$ part~$E_0$
    is the class of an $\mathrm{O} (n)$-bundle on~$X$.
    Suppose that~$E$ is equipped with an orientation,
    so in particular, $E_0$ is also an $\mathrm{SO} (n)$-bundle.

    We define an \emph{equivariant square root Euler class}
    $\sqrte_z (E)$,
    which is a real analogue of the class $e_z (E)$ in \cref{para-equivar-euler}.

    Consider the hyperplane arrangement on~$\Lambda_T \otimes \mathbb{R}$
    given by hyperplanes dual to non-zero $T$-weights appearing in~$E$.
    For an open chamber
    $\sigma \subset \Lambda_T \otimes \mathbb{R}$
    of this hyperplane arrangement, define
    \begin{equation}
        \sqrte_z (E) =
        \sqrte (E_0) \cdot
        \prod_{\langle \sigma, \lambda \rangle > 0}
        e_z (E_\lambda)
        \quad \in \quad
        \prod_{k = 0}^\infty
        \mathrm{H}^{2k} (X; \mathbb{Q})
        ( z_1, \dotsc, z_n ) \ ,
    \end{equation}
    where $E_0$ is equipped with the orientation
    induced by writing $E_0 = E - \sum_{\lambda \neq 0} E_\lambda$,
    where the last term has the orientation given by
    $\sum_{\langle \sigma, \lambda \rangle > 0} E_\lambda$,
    that is,
    factoring $\sum_{\lambda \neq 0} E_\lambda \colon X \to \mathrm{BO} \times \mathbb{Z}$
    as $X \to \mathrm{BU} \times \mathbb{Z} \to \mathrm{BSO} \times \mathbb{Z} \to \mathrm{BO} \times \mathbb{Z}$,
    with the first map given by
    $\sum_{\langle \sigma, \lambda \rangle > 0} E_\lambda$,
    and the second map given by $E \mapsto E + E^\vee$.

    One can verify that
    $\sqrte_z (E)$ does not depend on the choice of~$\sigma$,
    while the orientation of~$E_0$ depends on~$\sigma$.
    The class $\sqrte_z (E)$ is only non-zero if~$E$ has even rank,
    and we have the relations
    \begin{align}
        \sqrte_z (E^\mathrm{op})
        & =
        -\sqrte_z (E) \ ,
        \\
        \sqrte_z (E)^2
        & =
        (-1)^{\operatorname{rank} E / 2} \cdot
        e_z (E) \ ,
    \end{align}
    where~$E^\mathrm{op}$ denotes~$E$ with the opposite orientation.

    When~$X$ is not connected,
    we define $\sqrte_z (E)$ on each connected component of~$X$ as above.

    We have the relation
    \begin{equation}
        \sqrte_z (E + F) = \sqrte_z (E) \, \sqrte_z (F) \ ,
    \end{equation}
    where $E + F$ is equipped with the induced orientation,
    which does not depend on the order of~$E$ and~$F$
    if they have even rank.
    In particular, if $E_0 = 0$, then
    $\sqrte_z (E) \, \sqrte_z (-E) = 1$,
    in which case we also write
    $\sqrte_z^{\,-1} (E) = \sqrte_z (-E)$.
\end{para}

\begin{para}[Orientation of the normal bundle]
    \label{para-nu-orientation}
    For a map of spaces $f \colon X \to Y$
    and classes $\mathbb{T}_X \in \mathit{KO} (X)$ and
    $\mathbb{T}_Y \in \mathit{KO} (Y)$,
    equipped with orientations,
    the \emph{normal bundle}
    $\bbnu_f = f^* (\mathbb{T}_Y) - \mathbb{T}_X \in \mathit{KO} (X)$
    is equipped with the induced orientation,
    defined by the property that the sum
    $f^* (\mathbb{T}_Y) = \bbnu_f + \mathbb{T}_X$
    preserves orientation.

    This depends on the order of the sum,
    in that if we used $\mathbb{T}_X + \bbnu_f$ instead,
    the orientation would be changed by a factor of
    $(-1)^{\operatorname{rank} \mathbb{T}_X \cdot \operatorname{rank} \bbnu_f}$.
    The orientation is reversed whenever
    the orientation of either~$\mathbb{T}_X$ or~$\mathbb{T}_Y$ is reversed.
\end{para}

\begin{para}[The paracategory~\texorpdfstring{$\mathcal{T}''$}{T''}]
    Recall the category~$\mathcal{T}'$
    from \cref{para-cat-t-prime}.
    We define a similar paracategory~$\mathcal{T}''$,
    using real $K$-theory, as follows:
    \begin{itemize}
        \item
            An object is a quintuple $(T_X, X, \odot_X, \mathbb{T}_X, o_X)$,
            where $(T_X, X, \odot_X) \in \mathcal{T}$,
            and~$\mathbb{T}_X \in \mathit{KO}^\circ (X)$ is a class of weight~$0$,
            and~$o_X$ is an orientation of~$\mathbb{T}_X$.

        \item
            A morphism
            $f \colon (T_X, X, \odot_X, \mathbb{T}_X) \to (T_Y, Y, \odot_Y, \mathbb{T}_Y)$
            is a morphism
            $f \colon (T_X, X, \odot_X) \to (T_Y, Y, \odot_Y)$ in~$\mathcal{T}$,
            such that the normal bundle
            $\bbnu_f = f^* (\mathbb{T}_Y) - \mathbb{T}_X$
            satisfies $\bbnu_f \in \mathit{KO}^\circ (X)$,
            and the part $(-\bbnu_f)_0$ of weight~$0$
            is the class of a vector bundle on~$X$.
    \end{itemize}
    There is a forgetful functor $\mathcal{T}'' \to \mathcal{T}'$,
    although this will not be used below.
\end{para}

\begin{example}[Shifted symplectic stacks]
    \label{eg-symp-stacks}
    Let~$\mathcal{X}$ be an $n$-shifted symplectic stack over~$\mathbb{C}$,
    with $n \in 4 \mathbb{Z} + 2$.
    An important special case is when $n = -2$,
    such as when~$\mathcal{X}$ is a moduli stack of coherent sheaves on
    a Calabi--Yau fourfold.

    In this case, the symplectic structure
    $\mathbb{T}_\mathcal{X} \simto \mathbb{L}_\mathcal{X} [n]$
    gives rise to an orthogonal structure
    on the complex $\mathbb{T}_\mathcal{X} [n/2]$, and hence a class
    $\mathbb{T}_{X} = -[\mathbb{T}_\mathcal{X} [n/2]]
    \in \mathit{KO} (X)$,
    where $X = | \mathcal{X} |$ is the topological realization.
    An orientation of~$\mathbb{T}_X$,
    as defined in \cref{para-ko},
    is equivalent to an isomorphism
    $\det (\mathbb{T}_{\mathcal{X}}) \simto \mathcal{O}_{\mathcal{X}}$
    that squares to the isomorphism
    $\det (\mathbb{T}_{\mathcal{X}})^{\otimes 2} \simeq \mathcal{O}_{\mathcal{X}}$
    given by the shifted symplectic structure,
    which agrees with
    \textcite[Definition~2.12]{borisov-joyce-2017}.
\end{example}

\begin{theorem}
    \label{thm-real-joyce-va}
    There is a functor of paracategories
    \begin{align*}
        V \colon \mathcal{T}''
        & \longrightarrow \mathsf{VS}_\mathbb{Q} \ ,
        \\
        X
        & \longmapsto
        V_X = \mathrm{H}_{\bullet + \operatorname{vdim}} (X; \mathbb{Q}) \ ,
        \\
        \smash{(X \overset{f}{\to} Y)}
        & \longmapsto
        f_* \circ
        \tau_X (z) \circ
        \bigl( (-) \cap \sqrte_z^{\,-1} (\bbnu_f) \bigr) \ ,
    \end{align*}
    where~$V_X$ is equipped with the translation operator~$\tau_X$
    defined in \cref{para-translation-op},
    and~$\bbnu_f$ is equipped with the induced orientation
    as in \cref{para-nu-orientation}.

    In particular, we have the following:
    \begin{enumerate}
        \item
            \label{item-real-va}
            Let~$X$ be as in \cref{para-va-setting},
            but with $K (-)$, $K^\circ (-)$ replaced by
            $\mathit{KO} (-)$, $\mathit{KO}^\circ (-)$ instead,
            and assume that $\mathbb{T}_X \in \mathit{KO} (X)$
            is equipped with an orientation.
            Then there is a $\mathbb{Z}$-graded weak vertex algebra structure
            on $\mathrm{H}_{\bullet + \vdim} (X; \mathbb{Q})$, defined by
            \begin{multline}
                \label{eq-def-real-va}
                a_1 (z_1) \cdots a_n (z_n) =
                (-1)^{\sum_{1 \leq i < j \leq n} |a_i| \vdim_j} \cdot {}
                \\
                (\oplus_{(n)})_* \circ \tau (z)
                \Bigl(
                    (a_1 \boxtimes \cdots \boxtimes a_n) \cap
                    \sqrte_z^{\,-1} (\bbnu_{(n)})
                \Bigr) \ ,
            \end{multline}
            where~$a_i \in \mathrm{H}_\bullet (X; \mathbb{Q})$
            are homogeneous elements,
            each supported on a single connected component of\/~$X$,
            and $\vdim_i$ denotes the rank of\/~$\mathbb{T}_X$
            on the support of~$a_i$.

        \item
            \label{item-real-vm}
            Let~$X, Y$ be as in \cref{para-vm-setting},
            but with $K (-)$, $K^\circ (-)$ replaced by
            $\mathit{KO} (-)$, $\mathit{KO}^\circ (-)$ instead,
            and assume that $\mathbb{T}_X$, $\mathbb{T}_Y$
            are equipped with orientations.
            Then there is a weak module structure
            on $\mathrm{H}_{\bullet + \vdim} (Y; \mathbb{Q})$
            for the weak vertex algebra
            $\mathrm{H}_{\bullet + \vdim} (X; \mathbb{Q})$
            in \cref{item-real-va}, defined by
            \begin{multline}
                \label{eq-def-real-vm}
                a_1 (z_1) \cdots a_n (z_n) \, m =
                (-1)^{\sum_{1 \leq i < j \leq n + 1} |a_i| \vdim_j} \cdot {}
                \\
                (\diamond_{(n)})_* \circ \tau (z)
                \Bigl(
                    (a_1 \boxtimes \cdots \boxtimes a_n \boxtimes m) \cap
                    \sqrte_z^{\,-1} (\bbnu_{\diamond, \smash{(n)}})
                \Bigr) \ ,
            \end{multline}
            where $a_i \in \mathrm{H}_\bullet (X; \mathbb{Q})$
            and $m \in \mathrm{H}_\bullet (Y; \mathbb{Q})$
            are homogeneous elements
            supported on single components,
            $\vdim_i$ is as in \cref{item-real-va},
            and $\vdim_{n + 1}$ is the rank of\/~$\mathbb{T}_Y$
            on the support of~$m$.

        \item
            \label{item-real-vi-stack}
            Let~$\mathcal{X}$ be as in \cref{thm-vi-stack},
            equipped with an $n$-shifted symplectic structure
            for $n \in 4 \mathbb{Z} + 2$, and an orientation.
            Then there is a functor
            \begin{align}
                \mathsf{Face}^\mathrm{sp} (\mathcal{X})^\mathrm{op}
                & \longrightarrow \mathsf{VS}_\mathbb{Q} \ ,
                \notag
                \\
                \alpha
                & \longmapsto
                V_\alpha =
                \mathrm{H}_{\bullet + \operatorname{vdim}} (\mathcal{X}_\alpha; \mathbb{Q}) \ ,
                \notag
                \\
                \smash{(\alpha \overset{f}{\to} \beta)}
                & \longmapsto
                (f^\circ)_* \circ
                \tau_{\beta} (z) \circ
                \bigl( (-) \cap \sqrte_z^{\,-1} (\bbnu_{f^\circ}) \bigr) \ .
            \end{align}
    \end{enumerate}
\end{theorem}

\begin{proof}
    For the main statement,
    the proof of \cref{thm-vi}
    can be adapted to this situation without much change,
    and we only need to check the real version of \cref{lem-swap},
    that is, for $(T, X, \odot) \in \mathcal{T}$
    and $E \in \mathit{KO}^\circ (X)$,
    with an orientation,
    such that~$E_0$ is the class of a vector bundle, we have
    \begin{equation}
        \tau (w) (a)
        \cap \sqrte_z (E)
        = \iota_{z, w}
        \bigl(
            \tau (w) (a \cap \sqrte_{z + w} (E))
        \bigr)
    \end{equation}
    for any $a \in \mathrm{H}_\bullet (X; \mathbb{Q})$.
    But it is straightforward to adapt the proof of \cref{lem-swap}
    to this situation.

    The statement \cref{item-real-va}
    follows from the main statement,
    where we precompose the vertex induction map
    $V_{X^n} \to V_X$ with the map
    $V_{\smash{X}}^{\otimes n} \to V_{X^n}$ given by
    \begin{equation*}
        a_1 \otimes \cdots \otimes a_n
        \longmapsto
        (-1)^{\sum_{1 \leq i < j \leq n} |a_i| \vdim_j} \cdot
        a_1 \boxtimes \cdots \boxtimes a_n \ ,
    \end{equation*}
    with notations as in the statement of~\cref{item-real-va},
    with the sign inserted so that it is equivariant under
    (signed) permutations of the factors.
    The sign is needed because in~$V_X$, compared to homology,
    the parity of elements is reversed on components with odd virtual dimension.

    The statement \cref{item-real-vm} holds analogously.

    For \cref{item-real-vi-stack},
    as in \cref{para-face-to-t-circ},
    there is a functor
    \begin{equation*}
        |\mathcal{X}_{(-)}| \colon
        \mathsf{Face} (\mathcal{X})^\mathrm{op}
        \longrightarrow \mathcal{T}'' \ ,
    \end{equation*}
    sending a face
    $\alpha \colon \Lambda \to \mathrm{CL} (\mathcal{X})$
    to the object
    $(|T_\Lambda|, |\mathcal{X}_\alpha|, \odot, \mathbb{T}_{\mathcal{X}_\alpha}, o_{\mathcal{X}_\alpha}) \in \mathcal{T}''$,
    where we choose an arbitrary orientation of~$\mathbb{T}_{\mathcal{X}_\alpha}$
    for each~$\alpha$, which exists by \cref{eq-tangent-x-alpha},
    as the class of~$\mathbb{T}_\mathcal{X}$ is orientable and
    $\alpha^\star (\mathbb{T}_\mathcal{X})_\lambda =
    \alpha^\star (\mathbb{T}_\mathcal{X})_{-\lambda}^\vee$
    in $K (|\mathcal{X}_\alpha|)$.
\end{proof}

\begin{para}
    \label{eg-real-va-symp}
    We can also apply
    \cref{thm-real-joyce-va}~\crefrange{item-real-va}{item-real-vm}
    to shifted symplectic stacks as follows.

    Suppose that~$\mathcal{X}$ is an oriented $n$-shifted symplectic stack over~$\mathbb{C}$,
    with $n \in 4 \mathbb{Z} + 2$,
    and that it is also equipped with the structure of
    a derived linear moduli stack,
    as in \cref{eg-linear-stacks}.
    Suppose that the weight condition~\cref{eq-cond-nu} is satisfied,
    which is usually the case, as discussed in \cref{eg-linear-stacks}.
    Then \cref{eq-def-real-va} defines a vertex algebra structure on
    $\mathrm{H}_{\bullet + \vdim} (\mathcal{X}; \mathbb{Q})$.
    If, moreover, $\mathcal{X}$ is equipped with a self-dual structure
    as in \cref{eg-vm-tw},
    such that its $\mathbb{Z}_2$-action preserves the
    shifted symplectic structure,
    then the fixed locus~$\mathcal{X}^\mathrm{sd}$
    is also $n$-shifted symplectic,
    and if it is also orientable,
    then \cref{eq-def-real-vm} defines a twisted module structure on
    $\mathrm{H}_{\bullet + \vdim} (\mathcal{X}^\mathrm{sd}; \mathbb{Q})$.

    For example, when~$\mathcal{X}$ is an open substack of
    the moduli stack of perfect complexes
    on a Calabi--Yau fourfold~$Z$,
    an orientation is given by
    \textcite[Theorem~13.7]{joyce-upmeier-bordism},
    under a topological assumption on~$Z$.
    In this case, self-dual structures may be obtained as in
    \cref{eg-vm-tw},
    but we do not know if orientations exist on~$\mathcal{X}^\mathrm{sd}$.

    Alternatively, we could also consider derived versions
    as in \cref{eg-dg-cat,eg-vm-dg-cat}.
    Namely, if~$\mathcal{C}$ and~$\mathcal{X}$
    are as in \cref{eg-dg-cat},
    and if~$\mathcal{C}$ is equipped with a $(2-n)$-Calabi--Yau structure,
    with $n \in 4 \mathbb{Z} + 2$,
    then~$\mathcal{X}$ has an $n$-shifted symplectic structure
    by \textcite[Theorem~5.6]{brav-dyckerhoff-2021}.
    If it is orientable,
    then \cref{eq-def-real-va} defines a vertex algebra structure on
    $\mathrm{H}_{\bullet + \vdim} (\mathcal{X}; \mathbb{Q})$,
    and if~$\mathcal{C}$ is self-dual as in
    \cref{eg-vm-dg-cat},
    where the involution is compatible with the Calabi--Yau structure,
    then~$\mathcal{X}^\mathrm{sd}$ is also $n$-shifted symplectic,
    and if it is also orientable,
    then \cref{eq-def-real-vm} defines a twisted module structure on
    $\mathrm{H}_{\bullet + \vdim} (\mathcal{X}^\mathrm{sd}; \mathbb{Q})$.
\end{para}

\subsection{A \textit{K}-theory version}
\label{subsec-k-va}

\begin{para}
    We introduce a $K$-theory version of vertex induction,
    generalizing Liu's~\cite{liu-k-wall-crossing} construction
    of a $K$-theory version of Joyce vertex algebras.
    We expect this to be useful for
    generalizing the $K$-theoretic invariants of \textcite{liu-k-wall-crossing}
    to more general quasi-smooth stacks.

    An important difference between the $K$-theory version
    and the homology version is that the former is
    \emph{multiplicative} in its nature,
    rather than additive.
    This is reflected in the fact that we obtain
    multiplicative vertex algebras and vertex spaces,
    rather than the usual versions.
\end{para}

\begin{para}[Multiplicative vertex algebras]
    \label{para-mult-va}
    Define a \emph{multiplicative vertex algebra} over~$K$
    to be the data $(V, (X_n)_{n \geq 0})$,
    consisting of a $K$-vector space~$V$
    and $K$-linear multiplication maps
    \begin{align*}
        X_n \colon
        V^{\otimes n}
        & \longrightarrow
        V \llbr x_1, \dotsc, x_n \rrbr \, [(x_i - x_j)^{-1}] \ ,
        \\
        a_1 \otimes \cdots \otimes a_n
        & \longmapsto
        X_n (a_1, \dotsc, a_n; z_1, \dotsc, z_n) \ ,
    \end{align*}
    where we write $z_i = 1 + x_i$,
    satisfying the following properties:

    \begin{enumerate}
        \item
            (\emph{Unit})
            For any $a \in V$, we have
            \begin{equation}
                \label{eq-mult-va-unit}
                X_1 (a; 1) = a \ .
            \end{equation}

        \item
            (\emph{Commutativity})
            For any homogeneous elements $a_1, \dotsc, a_n \in V$,
            and any permutation $\sigma \in \mathfrak{S}_n$,
            we have
            \begin{equation}
                \label{eq-mult-va-comm}
                X_n (a_{\sigma (1)}, \dotsc, a_{\sigma (n)}; z_{\sigma (1)}, \dotsc, z_{\sigma (n)})
                = X_n (a_1, \dotsc, a_n; z_1, \dotsc, z_n) \ .
            \end{equation}

        \item
            (\emph{Associativity})
            For integers $m, n \geq 0$ and elements
            $b_1, \dotsc, b_m, a_1, \dotsc, a_n \in V$,
            we have
            \begin{multline}
                \label{eq-mult-va-assoc}
                X_{n+1} \Bigl(
                    X_m (b_1, \dotsc, b_m; w_1, \dotsc, w_m),
                    a_1, \dotsc, a_n; \
                    z_0, \dotsc, z_n
                \Bigr)
                \\[-.5ex]
                =
                \iota_{ \{ z_i - 1 \}, \{ w_j - 1 \} } \,
                X_{m+n} \bigl(
                    b_1, \dotsc, b_m, a_1, \dotsc, a_n;
                    z_0 w_1, \dotsc, z_0 w_m,
                    z_1, \dotsc, z_n
                \bigr) \ .
            \end{multline}
    \end{enumerate}
    Similarly, we define \emph{weak multiplicative vertex algebras}
    as above, but with the codomain of~$X_n$ enlarged to
    $V \llparen x_1, \dotsc, x_n \rrparen$.

    As usual, we abbreviate the product
    $X_n (a_1, \dotsc, a_n; z_1, \dotsc, z_n)$
    as $a_1 (z_1) \cdots a_n (z_n)$.
\end{para}

\begin{para}[Multiplicative modules]
    We have a similar notion of a \emph{module}
    over a multiplicative vertex algebra~$(V, (X_n)_{n \geq 0})$,
    defined as a $K$-vector space~$M$
    equipped with $K$-linear multiplication maps
    \begin{align*}
        X^M_n \colon
        V^{\otimes n} \otimes M
        & \longrightarrow
        M \llbr x_1, \dotsc, x_n \rrbr \,
        [x_i^{-1}, (x_i - x_j)^{-1}] \ ,
        \\
        a_1 \otimes \cdots \otimes a_n \otimes m
        & \longmapsto
        X^M_n (a_1, \dotsc, a_n, m; z_1, \dotsc, z_n) \ ,
    \end{align*}
    where we write $z_i = 1 + x_i$, with the following properties:
    \begin{enumerate}
        \item
            (\emph{Unit})
            We have $X_0^M = \mathrm{id}_M$.

        \item
            (\emph{Associativity})
            For integers $k, n \geq 0$ and elements
            $a_1, \dotsc, a_n, b_1, \dotsc, b_k \in V$, we have
            \begin{align}
                & X_{n+1}^M \Bigl(
                    X_k (b_1, \dotsc, b_k; w_1, \dotsc, w_k);
                    a_1, \dotsc, a_n, m;
                    z_0, \dotsc, z_n
                \Bigr)
                \notag
                \\[-.5ex]
                & \hspace{2em} =
                \iota_{ \{ z_i - 1 \}, \{ w_j - 1 \} } \,
                X_{n+k}^M \bigl(
                    b_1, \dotsc, b_k,
                    a_1, \dotsc, a_n, m;
                    z_0 w_1, \dotsc, z_0 w_k,
                    z_1, \dotsc, z_n
                \bigr) \ ,
                \\[1ex]
                & X_{n}^M \Bigl(
                    a_1, \dotsc, a_n, X_k^M (b_1, \dotsc, b_k, m; w_1, \dotsc, w_k);
                    z_1, \dotsc, z_n
                \Bigr)
                \notag
                \\[-.5ex]
                & \hspace{2em} =
                \iota_{ \{ z_i - 1 \}, \{ w_j - 1 \} } \,
                X_{n+k}^M \bigl(
                    a_1, \dotsc, a_n, b_1, \dotsc, b_k, m;
                    z_1, \dotsc, z_n,
                    w_1, \dotsc, w_k
                \bigr) \ .
            \end{align}
    \end{enumerate}

    Define a \emph{weak module} over a weak multiplicative vertex algebra
    in the same way as above,
    but with the codomain of~$X_n^M$ enlarged to
    $M \llparen x_1, \dotsc, x_n \rrparen$.

    Similarly to \cref{para-twisted-modules},
    we also have a notion of \emph{twisted modules}.
    Namely, if a multiplicative vertex algebra~$(V, (X_n)_{n \geq 0})$
    is equipped with a \emph{twisted involution},
    meaning an involution $(-)^\vee \colon V \simto V^\mathrm{op}$
    such that
    \begin{equation}
        \bigl( a_1 (z_1) \cdots a_n (z_n) \bigr)^\vee =
        a_1^\vee (z_1^{-1}) \cdots a_n^\vee (z_n^{-1})
    \end{equation}
    for all $a_1, \dotsc, a_n \in V$,
    then we define a \emph{twisted module} for~$V$
    as a weak module such that the image of~$X^M_n$ lies in the subspace
    $M \llbr x_1, \dotsc, x_n \rrbr \, [x_i^{-1}, (x_i - x_j)^{-1}, (x_i + x_j + x_i x_j)^{-1} : i \neq j]$,
    and such that
    $a^\vee (z) \, m = a (z^{-1}) \, m$
    for all $a \in V$ and $m \in M$.
    Here, inverting $x_i + x_j + x_i x_j$
    means we allow poles also at $z_i z_j = 1$ for $i \neq j$.

    As usual, we abbreviate
    $X^M_n (a_1, \dotsc, a_n, m; z_1, \dotsc, z_n)$
    as $a_1 (z_1) \cdots a_n (z_n) \, m$.
\end{para}

\begin{para}[Topological \texorpdfstring{$K$}{K}-homology]
    \label{para-k-homology}
    We now define the underlying vector space
    where our multiplicative vertex algebras and vertex spaces will live.

    Let $X \in \mathsf{hCW}$ be a space.
    Define the \emph{topological $K$-homology} of~$X$ by
    \begin{equation}
        K_\circ (X) =
        \uppi_0 (\mathit{KU} \wedge \Sigma^\infty X_+) \ ,
        \qquad
        K_\circ (X; \mathbb{Q}) =
        K_\circ (X) \otimes \mathbb{Q} \ ,
    \end{equation}
    where $\mathit{KU}$ is the $K$-theory spectrum,
    $X_+ = X \sqcup \{ * \}$ is~$X$ with an extra base point,
    and $\Sigma^\infty$ denotes taking the suspension spectrum.
    This is the zeroth homology group
    in the generalized homology theory represented by $\mathit{KU}$;
    see \textcite[Part~III, \S 6]{adams-1974} for background on this.

    For example,
    we have
    \begin{equation}
        \label{eq-k-homology-bu}
        K_\circ (\mathrm{BU} (1))
        \simeq \mathbb{Z} [\ell] \ ,
        \qquad
        \ell^k \cap (-) = \frac{1}{k!} \cdot
        \Bigl( \frac{\partial}{\partial L} \Bigr)^k \ ,
    \end{equation}
    where $\ell^k \cap (-)$ acts on
    $K (\mathrm{BU} (1); \mathbb{Z})
    \simeq \mathbb{Z} \llbr L - 1 \rrbr$,
    and $L \to \mathrm{BU} (1)$ is the universal line bundle.
    See \textcite[Part~II, Lemma 2.14]{adams-1974} for this fact;
    compare \textcite[Proposition~2.3.5]{liu-k-wall-crossing}.

    The following lemma characterizes a finiteness property
    of topological $K$-homology,
    which will be useful in dealing with convergence issues
    in \cref{para-wedge-z} later.
\end{para}

\begin{lemma}
    \label{lem-k-hom-fin}
    Let $X \in \mathsf{hCW}$ be a space.
    Then for any class $a \in K_\circ (X; \mathbb{Q})$,
    there exists $N \in \mathbb{N}$,
    such that $a \cap E = 0$
    for all $E \in K (X; \mathbb{Q})$
    with $\mathrm{ch}_i (E) = 0$ for all $i < N$.
\end{lemma}

\begin{proof}
    Choose a realization of~$X$ as a CW complex.
    Then~$a$ is supported on a finite subcomplex $U \subset X$,
    since by its definition,
    $K_\circ (X)$ is a colimit of
    $K_\circ (U)$ over finite subcomplexes $U \subset X$.
    We may then choose $N = \dim U + 1$,
    so that $a \cap E =
    i_* (a \cap E|_U) = 0$,
    because $E|_U$ is in the kernel of the isomorphism
    $\mathrm{ch} \colon K (U; \mathbb{Q}) \simto \mathrm{H}^{2 \bullet} (U; \mathbb{Q})$,
    and is hence zero.
\end{proof}

\begin{para}[The multiplicative translation operator]
    \label{para-translation-op-mul}
    Let $(T, X, \odot) \in \mathcal{T}$.
    Let $z = (z_1, \dotsc, z_n)$
    be a set of coordinates on~$\Lambda_T$,
    and write $x_i = z_i - 1$.

    Following ideas of \textcite[\S 3.3.6]{liu-k-wall-crossing},
    define the \emph{multiplicative translation operator}
    \begin{align*}
        D (z) \colon
        K_\circ (X; \mathbb{Q})
        & \longrightarrow
        K_\circ (X; \mathbb{Q}) \, \llbr x_i \rrbr \ ,
        \\
        a
        & \longmapsto
        \odot_* \biggl( {}
            \sum_{k_1, \dotsc, k_n \geq 0} {}
            x^k \cdot (\ell^k \boxtimes a)
        \biggr) \ ,
    \end{align*}
    where $x^k = x_1^{\smash{k_1}} \cdots x_n^{\smash{k_n}}$,
    $\ell^k = \ell_1^{\smash{k_1}} \boxtimes \cdots \boxtimes \ell_n^{\smash{k_n}}$,
    and $\ell_i^{\smash{k_i}} \in K_\circ (\mathrm{BU} (1))$
    is the element in~\cref{eq-k-homology-bu}.

    By the computation in
    \textcite[Lemma~3.3.7]{liu-k-wall-crossing}, we have
    \begin{equation}
        \label{eq-translation-op-mul}
        D (z) (a) \cap E_\lambda =
        z^{\lambda} \cdot D (z) (a \cap E_\lambda)
    \end{equation}
    for $a \in K_\circ (X; \mathbb{Q})$
    and a class $E_\lambda \in K^\circ (X)$
    of weight $\lambda \in \Lambda^T$,
\end{para}

\begin{para}[Multiplicative vertex spaces]
    Similarly to \cref{para-vertex-spaces},
    we also define \emph{multiplicative vertex spaces},
    which we use to describe multiplicative versions
    of the constructions in~\cref{sec-vi}.

    A \emph{multiplicative vertex space} over a field~$K$
    is a triple $(\Lambda, V, D)$,
    consisting of a finite-dimensional $K$-vector space~$\Lambda$,
    a $K$-vector space~$V$, and a $K$-linear map
    $D (z) \colon V \to V \llbr \Lambda \rrbr$,
    called the \emph{translation operator},
    where we set $z = (z_1, \dotsc, z_n)$
    for $z_i = 1 + x_i$,
    for $(x_1, \dotsc, x_n)$ a set of coordinates for~$\Lambda$,
    satisfying
    \begin{align}
        D (1) & = \mathrm{id}_V \ ,
        \\
        D (z) \circ D (w) & = D (z w) \ .
    \end{align}
    In other words, $D$ is a representation
    of the formal group $\hat{\mathbb{G}}_\mathrm{m}^n$.

    We define maps of multiplicative vertex spaces
    similarly to \cref{para-map-vertex-spaces},
    and define compositions similarly to \cref{para-composition-vertex-spaces}.
    We obtain a paracategory
    $\mathsf{VS}^\mathrm{mul}_K$
    of multiplicative vertex spaces over~$K$.
\end{para}

\begin{para}[The equivariant exterior power]
    \label{para-wedge-z}
    Let $(T, X, \odot) \in \mathcal{T}$,
    and $E \in K^\circ (X)$ be a class
    such that $E_0$ is the class of a vector bundle.

    Heuristically, we wish to define the equivariant exterior power
    ${\wedge}_{-z} (E)$ by
    \begin{equation}
        \label{eq-def-wedge-z}
        {\wedge}_{-z} (E) =
        \prod_{\lambda \in \Lambda \lowersup{T}} {}
        \biggl( {}
            \sum_{k \geq 0} {}
            (-z^\lambda)^{k} \cdot {\wedge}^k (E_\lambda)
        \biggr) \ ,
    \end{equation}
    where $z = (z_1, \dotsc, z_n)$
    is a set of coordinates on~$\Lambda_T$,
    and we write
    $z^\lambda = z_1^{\smash{\lambda_1}} \cdots z_n^{\smash{\lambda_n}}$,
    etc.
    However, the expression~\cref{eq-def-wedge-z}
    does not make sense yet,
    since the sum can be infinite
    when~$E_\lambda$ is not the class of a vector bundle.

    To fix this, we expand the sum
    in the variables $x_i = z_i - 1$ as
    \begin{multline}
        \label{eq-def-wedge-x}
        {\wedge}_{-z} (E) =
        \biggl( {}
            \sum_{k \geq 0} {}
            (-1)^k \cdot {\wedge}^k (E_0)
        \biggr) \cdot {}
        \\[-1ex]
        \prod_{\lambda \in \Lambda \lowersup{T} \setminus \{ 0 \}} {}
        \biggl( {}
            \bigl( 1 - (1 + x)^\lambda \bigr)^{\operatorname{rank} E_\lambda} \cdot
            \sum_{k \geq 0} {} \biggl(
                \frac{-(1 + x)^\lambda}{1 - (1 + x)^\lambda}
            \biggr)^k \cdot
            {\vee}^k (E_\lambda)
        \biggr)
        \quad {\in} \quad
        K^\circ (X) \llparen x_i \rrparen^\wedge
        \ ,
    \end{multline}
    where $E_0$ is a vector bundle by assumption,
    $(1 + x)^\lambda = \prod_i {} (1 + x_i)^{\smash{\lambda_i}}$,
    $K^\circ (X) \llparen x_i \rrparen^\wedge$
    denotes the completion of $K^\circ (X) \llparen x_i \rrparen$
    with respect to $I^N (X) \llparen x_i \rrparen$,
    where $I^N (X) \subset K^\circ (X)$
    is the ideal of classes~$E$
    with $\mathrm{ch}_i (E) = 0$ for all $i < N$,
    and
    \begin{equation}
        \label{eq-def-vee}
        \vee^k (E_\lambda) =
        \sum_{i = 0}^k {} (-1)^{k-i} \cdot
        \binom{\operatorname{rank} E_\lambda - i}{k - i} \cdot
        {\wedge}^i (E_\lambda) \ ,
    \end{equation}
    which satisfies $\vee^k (E_\lambda) \in I^k (X)$.
    For example, if $E_\lambda = L_1 + \cdots + L_r$
    is a sum of line bundles, then
    $\vee^k (E_\lambda) =
    \sum_{1 \leq i_1 < \cdots < i_k \leq r} {}
    (L_{i_1} - 1) \cdots (L_{i_k} - 1)$.

    We take~\cref{eq-def-wedge-x}
    as the definition of ${\wedge}_{-z} (E)$ from now on.
    By \cref{lem-k-hom-fin},
    the operation $(-) \cap {\wedge}_{-z} (E)$
    is well-defined on $K_\circ (X; \mathbb{Q})$.
    We have the relation
    \begin{equation}
        {\wedge}_{-z} (E + F) = {\wedge}_{-z} (E) \cdot {\wedge}_{-z} (F)
    \end{equation}
    whenever the right-hand side is defined.
\end{para}

\begin{theorem}
    \label{thm-k-va}
    There is a functor of paracategories
    \begin{align*}
        V \colon \mathcal{T}'
        & \longrightarrow \mathsf{VS}_\mathbb{Q}^\mathrm{mul} \ ,
        \\
        X
        & \longmapsto
        K_\circ (X; \mathbb{Q}) \ ,
        \\
        \smash{(X \overset{f}{\to} Y)}
        & \longmapsto
        f_* \circ
        D_X (z) \circ
        \bigl( (-) \cap \wedge_{-z} (-\bbnu_f)^\vee \bigr) \ ,
    \end{align*}
    where~$K_\circ (X; \mathbb{Q})$ is equipped with the translation operator~$D_X$
    defined in \cref{para-translation-op-mul}.

    In particular, we have the following:
    \begin{enumerate}
        \item
            \label{item-k-va}
            Let~$X$ be as in \cref{para-va-setting}.
            Then the assignment
            \begin{equation}
                \label{eq-def-k-va}
                a_1 (z_1) \cdots a_n (z_n) =
                (\oplus_{(n)})_* \circ D (z)
                \Bigl(
                    (a_1 \boxtimes \cdots \boxtimes a_n) \cap
                    {\wedge}_{-z} (-\bbnu_{\smash{(n)}})^\vee
                \Bigr)
            \end{equation}
            defines a weak multiplicative vertex algebra structure
            on $K_\circ (X; \mathbb{Q})$,
            with poles along the locus
            $z_1^{\smash{\lambda_1}} \cdots z_n^{\smash{\lambda_n}} = 1$
            for non-zero weights $\lambda \in \mathbb{Z}^n$
            appearing in $\oplus_{(n)}^* (\mathbb{T}_X)$.
            In particular, if the class
            $\oplus^* (\mathbb{T}_X) \in K^\circ (X^2)$
            only has $\mathrm{U} (1)^2$-weights
            $(-1, 1),$ $(0, 0),$ and~$(1, -1),$
            then $K_\circ (X; \mathbb{Q})$ is a multiplicative vertex algebra.

        \item
            \label{item-k-vm}
            Let~$X, Y$ be as in \cref{para-vm-setting}.
            Then the assignment
            \begin{equation}
                \label{eq-def-k-vm}
                a_1 (z_1) \cdots a_n (z_n) \, m =
                (\diamond_{(n)})_* \circ
                D (z)
                \Bigl(
                    (a_1 \boxtimes \cdots \boxtimes a_n \boxtimes m)
                    \cap {\wedge}_{-z} (-\bbnu_{\diamond, \smash{(n)}})^\vee
                \Bigr)
            \end{equation}
            defines a weak module structure on
            $K_\circ (Y; \mathbb{Q})$
            for the weak multiplicative vertex algebra
            $K_\circ (X; \mathbb{Q})$,
            with poles along the locus
            $z_1^{\smash{\lambda_1}} \cdots z_n^{\smash{\lambda_n}} = 1$
            for non-zero weights $\lambda \in \mathbb{Z}^n$
            appearing in $\diamond_{(n)}^* (\mathbb{T}_Y)$.

        \item
            \label{item-k-vi-stack}
            Let~$\mathcal{X}$ be as in \cref{thm-vi-stack}.
            Then there is a functor
            \begin{align}
                \mathsf{Face}^\mathrm{sp} (\mathcal{X})^\mathrm{op}
                & \longrightarrow \mathsf{VS}^\mathrm{mul}_\mathbb{Q} \ ,
                \notag
                \\
                \alpha
                & \longmapsto
                K_\circ (\mathcal{X}_\alpha; \mathbb{Q}) \ ,
                \notag
                \\
                \smash{(\alpha \overset{f}{\to} \beta)}
                & \longmapsto
                (f^\circ)_* \circ
                D_{\mathcal{X}_\beta} (z) \circ
                \bigl( (-) \cap {\wedge}_{-z} (-\bbnu_{f^\circ})^\vee \bigr) \ .
            \end{align}
    \end{enumerate}
\end{theorem}

\begin{proof}
    For the main statement,
    we need to verify that the construction preserves composition.
    Let $X \overset{\smash[t]{f}}{\to} Y \overset{\smash[t]{g}}{\to} Z$
    be composable morphisms in~$\mathcal{T}'$.
    Then for any $a \in K_\circ (X)$, we have
    \begin{align*}
        &
        g_* \circ D_Y (z) \Bigl[
            f_* \circ D_X (w)
            \Bigl(
                a \cap {\wedge}_{-w} (-\bbnu_f)^\vee
            \Bigr)
            \cap {\wedge}_{-z} (-\bbnu_g)^\vee
        \Bigr]
        \\
        = {}
        &
        g_* \circ f_* \circ D_X (f^\sharp (z))
        \Bigl[
            D_X (w)
            \Bigl(
                a \cap {\wedge}_{-w} (-\bbnu_f)^\vee
            \Bigr)
            \cap {\wedge}_{-z} (f^* (-\bbnu_g))^\vee
        \Bigr]
        \\
        = {}
        &
        g_* \circ f_* \circ D_X (f^\sharp (z)) \circ D_X (w)
        \Bigl(
            a \cap {\wedge}_{-w} (-\bbnu_f)^\vee
            \cap {\wedge}_{-f^{\smash{\sharp}} (z) - w} (f^* (-\bbnu_g))^\vee
        \Bigr)
        \\
        = {}
        &
        (g \circ f)_* \circ D_X (f^\sharp (z) + w)
        \Bigl(
            a \cap {\wedge}_{-f^{\smash{\sharp}} (z) - w} (-\bbnu_{g \circ f})^\vee
        \Bigr) \ ,
    \end{align*}
    where the third line uses
    \cref{lem-swap-k} below,
    and the last step uses the fact that
    $\bbnu_f$ has weight~$0$ with respect to~$T_Y$.

    It is straightforward to check that
    the other statements follow from the main statement.
\end{proof}

\begin{lemma}
    \label{lem-swap-k}
    Let $(T, X, \odot) \in \mathcal{T}$,
    and let~$E \in K^\circ (X)$
    such that $E_0$ is the class of a vector bundle.
    Then for any $a \in K_\circ (X)$, we have
    \begin{equation}
        D (w) (a)
        \cap {\wedge}_{-z} (E)
        = \iota_{z - 1, w - 1}
        \bigl(
            D (w) (a \cap {\wedge}_{-zw} (E))
        \bigr) \ .
    \end{equation}
\end{lemma}

This is a $K$-theory version of
\cref{lem-swap}.
See also \textcite[Lemma~3.3.12]{liu-k-wall-crossing}.

\begin{proof}
    As ${\wedge}_{-z} (E)$ is multiplicative in~$E$,
    it is enough to prove this
    when~$E$ is of pure weight~$\lambda$.
    If $\lambda = 0$, then both sides are equal to
    $a \cap {\wedge}_{-1} (E)$.
    If $\lambda \neq 0$,
    by \cref{eq-translation-op-mul},
    it is enough to show that
    \begin{equation}
        \label{eq-pf-swap-k}
        a \cap ( w^{\mathrm{wt}} \cdot {\wedge}_{-z} (E) )
        = \iota_{z - 1, w - 1} (a \cap {\wedge}_{-z w} (E))
    \end{equation}
    for any $a \in K_\circ (X; \mathbb{Q})$,
    where $w^{\mathrm{wt}}$ is the operator that acts as~$w^\mu$
    on the weight~$\mu$ component.
    To prove~\cref{eq-pf-swap-k},
    we may restrict both sides to a finite subcomplex $U \subset X$
    where~$a$ is supported,
    as in the proof of \cref{lem-k-hom-fin}.
    Both sides are, by definition,
    expansions of rational functions in~$z$ and~$w$ valued in
    $K_\circ (U; \mathbb{C})$,
    and these rational functions are equal to the series
    $\sum_{k \geq 0} {} (-z^\lambda w^\lambda)^k \cdot (a \cap {\wedge}^k (E))$
    in the region of $(z, w) \in \mathbb{C}^{2n}$
    where the series converges.
    But ${\wedge}^k (E) |_U$ has polynomial growth in~$k$,
    which can be seen by writing it in terms of the
    $\vee^j (E) |_U$ as in \cref{para-wedge-z},
    which vanish for $j \gg 0$.
    The series thus converges in the region
    $|z^\lambda w^\lambda| < 1$,
    which implies that the two rational functions are equal.
\end{proof}

\section{Examples}

\label{sec-ex}

\subsection{Classifying stacks}

\begin{para}
    In this section, we discuss the following examples
    of our construction:

    \begin{enumerate}
        \item
            \label{item-eg-perf}
            A vertex algebra structure on the homology of
            $|\calPerf| \simeq \mathrm{BU} \times \mathbb{Z}$,
            the classifying stack of perfect complexes over~$\mathbb{C}$.
            This was due to \textcite{joyce-hall}
            and discussed in \textcite[\S 2.6.29]{latyntsev-thesis},
            and is related to the Heisenberg vertex algebra.

        \item
            \label{item-eg-perf-osp}
            Twisted module structures on the homology of
            $\mathrm{BO} \times \mathbb{Z}$
            and $\mathrm{BSp} \times \mathbb{Z}$.

        \item
            \label{item-eg-bg}
            Vertex induction for the homology of
            $\mathrm{B} G$ for reductive groups~$G$ over~$\mathbb{C}$.
    \end{enumerate}
\end{para}

\begin{para}[For reductive groups]
    \label{para-eg-bg}
    We begin with the example~\cref{item-eg-bg} above.

    Let~$G$ be a linearly reductive algebraic group over~$\mathbb{C}$,
    with maximal torus $T \subset G$.
    Let $\mathrm{Z} (G)^\circ \subset G$
    be the neutral component of the centre of~$G$.
    Consider the object $(\mathrm{Z} (G)^\circ, \mathrm{B} G, \odot) \in \mathcal{T}$,
    where $\mathrm{B} \mathrm{Z} (G)^\circ$ acts on $\mathrm{B} G$ by translation.
    We identify
    \begin{equation*}
        \mathrm{H}^\bullet (\mathrm{B} G; \mathbb{Q}) \simeq
        \mathbb{Q} [x_1, \dotsc, x_n]^W \ ,
    \end{equation*}
    where $x_i \in \mathrm{H}^2 (\mathrm{B} T; \mathbb{Q})$
    is the $i$-th standard generator
    upon an identification $T \simeq \mathbb{G}_\mathrm{m}^n$,
    and $W = \mathrm{N}_G (T) / \mathrm{Z}_G (T)$
    is the Weyl group.
    Dually, we have
    \begin{equation*}
        \mathrm{H}_\bullet (\mathrm{B} G; \mathbb{Q}) \simeq
        \mathbb{Q} [X_1, \dotsc, X_n]_W \ ,
    \end{equation*}
    where $(-)_W$ denotes taking coinvariants,
    and regarding $x_i$, $X_i$
    as (co)homology classes of~$\mathrm{B} T$,
    $x_i$ acts as $\partial / \partial X_i$ via the cap product.

    We set $\mathbb{T}_{\mathrm{B} G} = -[\mathfrak{g}]$,
    with the adjoint action of~$G$ on its Lie algebra~$\mathfrak{g}$.

    In fact, we are now in the situation of \cref{eg-symp-stacks},
    where $* / G$ admits a $2$-shifted symplectic structure
    given by choosing a Weyl-invariant inner product on~$\mathfrak{g}$.
    This symplectic structure is orientable if and only if
    $\uppi_0 (G)$ acts trivially on $\det (\mathfrak{g})$
    via the adjoint action of~$G$,
    which is always true if~$G$
    has an odd number of connected components.

    Assuming such an orientation exists,
    for a Levi subgroup $L \subset G$, we have
    \begin{align}
        \label{eq-sqrt-ez-bl-bg}
        \sqrte_z^{\,-1} (\bbnu_{\mathrm{B} L \to \mathrm{B} G})
        & =
        \pm \prod_{\lambda \in \Phi'_G} {}
        \lambda (z - x) \ ,
        \\
        \label{eq-ez-bl-bg}
        e_z^{-1} (\bbnu_{\mathrm{B} L \to \mathrm{B} G})
        & =
        (-1)^{(\dim G - \dim L) / 2} \cdot
        \sqrte_z^{\,-1} (\bbnu_{\mathrm{B} L \to \mathrm{B} G})^2 \ ,
    \end{align}
    where $\Phi'_G$ is the set of roots~$\lambda$ of~$G$
    such that $\langle \lambda, \mu \rangle > 0$
    for a cocharacter $\mu \colon \mathbb{G}_\mathrm{m} \to L$
    with $L = \mathrm{Z}_G (\mu)$,
    and the `$\pm$' sign depends on the choice of orientations
    on~$\mathrm{B} L$ and~$\mathrm{B} G$ and the choice of~$\mu$.
    For example, reversing~$\mu$ will change the sign by
    $(-1)^{|\Phi'_G|} = (-1)^{(\dim G - \dim L) / 2}$.

    We can then write down the vertex induction map
    \begin{equation}
        \mathrm{H}_{\bullet - \dim L} (\mathrm{B} L; \mathbb{Q})
        \longrightarrow
        \mathrm{H}_{\bullet - \dim G} (\mathrm{B} G; \mathbb{Q})
        \llbr z_1, \dotsc, z_k \rrbr \,
    \end{equation}
    defined in \cref{thm-real-joyce-va} explicitly
    using \cref{eq-sqrt-ez-bl-bg},
    where $z_1, \dotsc, z_k$
    is a set of coordinates on
    $\Lambda_{\mathrm{Z} (L)^\circ}$,
    the coweight lattice of the neutral component
    of the centre of~$L$,
    and the map has no poles in this case.
    Namely, it is given by
    \begin{equation}
        f \longmapsto
        \pm \exp \biggl( {}
            \sum_{i=1}^n
            z'_i X_i
        \biggr) \cdot
        \biggl( {}
            \prod_{\lambda \in \Phi'_G} {}
            \lambda \Bigl(
                z' - \frac{\partial}{\partial X}
            \Bigr)
        \biggr) \ f \ ,
    \end{equation}
    where~$f$ is a polynomial in $X_1, \dotsc, X_n$,
    and~$z'_i$ is a linear combination of the~$z_i$
    induced by the inclusion $\mathrm{Z} (L)^\circ \hookrightarrow T$,
    and the `$\pm$' sign depends on the choice of orientations.
    Similarly, this can be done for the version in \cref{thm-vi-stack}
    using \cref{eq-ez-bl-bg}.
\end{para}

\begin{para}[Perfect complexes]
    \label{para-eg-perf}
    Let~$\calPerf$ be the classifying stack
    of perfect complexes over~$\mathbb{C}$,
    and consider the space
    \begin{equation*}
        X = | \calPerf | \simeq \mathrm{BU} \times \mathbb{Z} \ ,
    \end{equation*}
    as a special case of \cref{eg-dg-cat}.
    Set $\mathbb{T}_X = -U^\vee \cdot U$,
    where $U \in K (X)$ is the universal class.
    The class~$\mathbb{T}_X$ agrees with the class of the tangent complex
    $\mathbb{T}_{\calPerf} = \mathcal{U}^\vee \otimes \mathcal{U} [1]$,
    where~$\mathcal{U}$ is the universal perfect complex on~$\calPerf$.
    We thus have $\vdim X_r = -r^2$,
    where $X_r \subset X$ denotes the component
    $\mathrm{BU} \times \{ r \}$ of rank~$r$.

    The class $\mathbb{T}_X$ lifts to a class
    $\mathbb{T}_X \in \mathit{KO} (X)$
    using the $2$-shifted symplectic structure on~$\calPerf$
    given by \textcite[\S 2.3]{pantev-toen-vaquie-vezzosi-2013},
    as discussed in \cref{eg-real-va-symp}.
    It admits an orientation given by a choice of identification
    $\det (\mathbb{T}_{\calPerf}) \simeq
    \det (\mathcal{U}^\vee)^r \otimes \det (\mathcal{U})^r
    \simeq \mathcal{O}_{\calPerf}$
    on the component of rank~$r$.

    We discuss two versions of Joyce vertex algebras in this case:
    The usual version from \cref{subsec-joyce-va},
    and the real version from \cref{subsec-real-va}.
    We have
    \begin{equation}
        \mathrm{H}_\bullet (X; \mathbb{Q})
        \simeq
        \bigoplus_{r \in \mathbb{Z}}
        \mathbb{Q} [s_1, s_2, \dotsc] \ ,
    \end{equation}
    where $s_i \in \mathrm{H}_{2 i} (\mathrm{BU}; \mathbb{Q})$
    are variables dual to the universal Chern characters,
    in that using the cap product,
    the $i$-th Chern character $\mathrm{ch}_i$
    acts as $\partial / \partial s_i$.
    We have
    \begin{equation}
        \bbnu_{(n)} =
        -\sum_{i \neq j} U_i^\vee \cdot U_j \ ,
    \end{equation}
    where $U_i \in K (X^n)$ is the pullback of
    the universal class of the $i$-th factor,
    so that using the universal relation
    \begin{equation}
        \sum_{i \geq 0} z^i \cdot c_i =
        \exp \biggl( {}
            \sum_{i > 0} {}
            (-1)^{i-1} \, (i - 1)! \,
            z^i \cdot \mathrm{ch}_i
        \biggr) \ ,
    \end{equation}
    where~$c_i$ is the $i$-th Chern class,
    and using a suitable choice of orientation,
    we have
    \begin{align}
        \label{eq-sqrt-ez-perf}
        \sqrte_z^{\,-1} (\bbnu_{(n)})
        & =
        \prod_{i < j} {} \biggl[
            (z_i - z_j)^{r_i r_j}
            \exp \biggl( {}
                \sum_{\substack{k, \ell \geq 0 \mathrlap{:} \\ k + \ell > 0}} {}
                (-1)^{k-1} (k+\ell-1)! \,
                (z_i - z_j)^{-k - \ell} \cdot
                \mathrm{ch}_{\smash{k}}^{(i)} \,
                \mathrm{ch}_{\smash{\ell}}^{(j)}
            \biggr)
        \biggr]
        \ ,
        \raisetag{3ex}
        \\
        e_z^{-1} (\bbnu_{(n)})
        & =
        (-1)^{\sum_{i < j} r_i r_j} \cdot
        \sqrte_z^{\,-1} (\bbnu_{(n)})^2 \ ,
    \end{align}
    where $i, j$ run over the range $1 \leq i < j \leq n$,
    and $(r_1, \dotsc, r_n) \in \mathbb{Z}^n$
    corresponds to a connected component of~$X^n$,
    and $\mathrm{ch}_{\smash{k}}^{(i)}$ is the $k$-th Chern character
    pulled back along the $i$-th projection $X^n \to X$.

    From these, one can write down explicit formulae
    for the vertex algebra structures
    given by
    \cref{thm-joyce-va,thm-real-joyce-va}.
    For example, the latter version is defined on the space
    \begin{equation*}
        \bigoplus_{r \in \mathbb{Z}}
        \mathrm{H}_{\bullet - r^2} (X_r; \mathbb{Q}) \ ,
    \end{equation*}
    whose vertex algebra structure is given by
    \begin{multline}
        a_1 (z_1) \cdots a_n (z_n)
        =
        (-1)^{\sum_{i < j} {} |a_i| \, r_j} \cdot
        \prod_{i < j} {}
        (z_i - z_j)^{r_i r_j} \cdot
        \exp \biggl[
            \sum_{i = 1}^n {}
            \biggl(
                z_i
                \sum_{k \geq 0}
                s_{\smash{k + 1}}^{(i)}
                \partial_{\smash{k}}^{(i)}
            \biggr)
        \biggr] \circ {}
        \\
        \exp \biggl[ {}
            \sum_{\substack{i < j; \ k, \ell \geq 0 \mathrlap{:} \\ k + \ell > 0}} {}
            \frac{(-1)^{k-1} (k+\ell-1)!}
            {(z_i - z_j) \lowerSup{k+\ell}} \cdot
            \partial_{\smash{k}}^{(i)}
            \partial_{\smash{\ell}}^{(j)}
        \biggr] \
        \Bigl[ a_1 (s_{\smash{k}}^{(1)}) \cdots a_n (s_{\smash{k}}^{(n)}) \Bigr] \
        \bigg|_{s_{\smash{k}}^{(i)} \, \mapsto \, s_k}
        \rlap{ ,}
    \end{multline}
    where each $a_i \in \mathrm{H}_\bullet (X; \mathbb{Q})$
    is homogeneous and supported on a single component $X_{r_i} \subset X$,
    and regarded as a polynomial in the variables~$s_{\smash{k}}^{(i)}$, and
    \begin{equation*}
        \partial_{\smash{k}}^{(i)} =
        \begin{cases}
            r_i
            & \text{if $k = 0$,}
            \\
            \partial / \partial s_{\smash{k}}^{(i)}
            & \text{if $k > 0$.}
        \end{cases}
    \end{equation*}
    As in \textcite[\S\S 2.6.25--29]{latyntsev-thesis},
    this is a \emph{lattice vertex algebra}
    on the lattice $\uppi_0 (X) \simeq \mathbb{Z}$,
    and the subalgebra $\mathrm{H}_\bullet (X_0; \mathbb{Q})$
    at $0 \in \mathbb{Z}$
    is isomorphic to the Heisenberg vertex algebra.
\end{para}

\begin{para}[Orthosymplectic complexes]
    \label{para-eg-perf-osp}
    \allowdisplaybreaks
    Let $X = \mathrm{BU} \times \mathbb{Z}$ as in \cref{para-eg-perf},
    and set
    \begin{equation*}
        Y =
        \mathrm{BO} \times 2 \mathbb{Z}
        \quad \text{or} \quad
        \mathrm{BO} \times (2 \mathbb{Z} + 1)
        \quad \text{or} \quad
        \mathrm{BSp} \times 2 \mathbb{Z} \ .
    \end{equation*}
    We refer to these cases as types~D, B, and~C, respectively.

    The space~$Y$ can be seen as a homotopy fixed locus
    of a $\mathbb{Z}_2$-action on
    $\mathrm{BU} \times 2 \mathbb{Z}$
    or $\mathrm{BU} \times (2 \mathbb{Z} + 1)$
    given by complex conjugation.
    See \textcite[Corollary~7.6]{dugger-2005}
    for the type~B and~D cases,
    and taking the $4$-fold loop space gives the case of type~C
    by Bott periodicity.

    Define $\mathbb{T}_Y \in \mathit{KO} (Y)$ by
    \begin{equation*}
        \mathbb{T}_{\mathrm{BO} \times \mathbb{Z}} =
        -{\wedge}^2 (U_{\mathrm{O}}) \ ,
        \qquad
        \mathbb{T}_{\mathrm{BSp} \times 2 \mathbb{Z}} =
        -\mathrm{Sym}^2 (U_{\mathrm{Sp}}) \ ,
    \end{equation*}
    where $U_{\mathrm{O}} \in \mathit{KO} (\mathrm{BO} \times \mathbb{Z})$ and
    $U_{\mathrm{Sp}} \in \mathit{KSp} (\mathrm{BSp} \times 2 \mathbb{Z})$
    are the universal classes,
    and we use the operation
    $\mathrm{Sym}^2 \colon \mathit{KSp} (-) \to \mathit{KO} (-)$.

    The class~$\mathbb{T}_Y$
    is orientable in types~B and~C,
    which can be checked via the $w_1$ class,
    where we have
    $w_1 ({\wedge}^2 (E)) = (\operatorname{rank} E - 1) \cdot w_1 (E)$
    for any~$E$, and
    $\mathrm{H}^1 (\mathrm{BSp}; \mathbb{Z}_2) = 0$.
    However, the same reasoning shows that it is not orientable in type~D.
    This is also similar to the orientability criterion
    that we obtained in \cref{para-eg-bg}.

    Let~$X$ act on~$Y$ via the map
    $\diamond \colon X \times Y \to Y$
    given by $x \diamond y = x \oplus y \oplus x^\vee$,
    with $x \oplus x^\vee$ equipped with the obvious orthosymplectic structure.
    More precisely, we define it as
    the homotopy $\mathbb{Z}_2$-fixed loci of the map
    $\oplus_{(3)} \colon X^3 \to X$,
    where $\mathbb{Z}_2$ acts on~$X^3$
    by complex conjugation followed by swapping the first and third factors.

    We have
    \begin{equation}
        \mathrm{H}_\bullet (Y; \mathbb{Q}) \simeq
        \bigoplus_{r \in \uppi_0 (Y)}
        \mathbb{Q} [s_2, s_4, \dotsc] \ ,
    \end{equation}
    where $s_{2i} \in \mathrm{H}_{4 i} (\mathrm{BO}; \mathbb{Q})$
    or $\mathrm{H}_{4 i} (\mathrm{BSp}; \mathbb{Q})$
    is dual to the $2i$-th Chern character,
    in a sense similar to that of \cref{para-eg-perf}.

    In types~B and~C, we have
    \begin{equation}
        \bbnu_{\diamond, \smash{(n)}} =
        -\sum_{i < j} {} (U_i + U_i^\vee) \cdot (U_j + U_j^\vee)
        -\sum_{i} {} \Bigl(
            (U_i + U_i^\vee) \cdot U_0
            + {\wedge}^2 (U_i)
            + {\wedge}^2 (U_i^\vee)
        \Bigr) \ ,
    \end{equation}
    and the twisted module structure given by
    \cref{thm-real-joyce-va}
    can be explicitly written as
    \begin{align}
        \hspace{1em}
        & \hspace{-1em}
        a_1 (z_1) \cdots a_n (z_n) \, m =
        (-1)^{\sum_{i < j} {} |a_i| \, r_j + \sum_i {} |a_i| \, r_0 (r_0 \mp 1) / 2} \cdot
        \prod_{i < j} {} (z_i^2 - z_j^2)^{r_i r_j} \cdot
        \prod_{i = 1}^n {} \bigl(
            z_i^{\smash{r_i s}} \cdot
            (2z_i)^{r_i (r_i \pm 1) / 2}
        \bigr) \cdot {}
        \notag \\*[-.5ex]
        &
        \exp \biggl[
            \sum_{i = 1}^n {}
            \biggl(
                z_i
                \sum_{k \geq 0}
                s_{\smash{k + 1}}^{(i)}
                \partial_{\smash{k}}^{(i)}
            \biggr)
        \biggr] \circ
        \exp \Biggl\{ {}
            \sum_{\substack{i < j; \ k, \ell \geq 0 \mathrlap{:} \\ k + \ell > 0}} {}
            \frac{(-1)^{k-1} (k+\ell-1)!}
            {(z_i - z_j) \lowerSup{k+\ell}} \cdot
            \partial_{\smash{k}}^{(i)}
            \partial_{\smash{\ell}}^{(j)}
        \notag \\[-1ex]
        &
            {} + \sum_{\substack{i < j; \ k, \ell \geq 0 \mathrlap{:} \\ k + \ell > 0}} {}
            \frac{(-1)^{k+\ell-1} (k+\ell-1)!}
            {(z_i + z_j) \lowerSup{k+\ell}} \cdot
            \partial_{\smash{k}}^{(i)}
            \partial_{\smash{\ell}}^{(j)}
            +
            \sum_{\substack{i; \ k, \ell \geq 0 \mathrlap{:} \\ k + \ell > 0}} {}
            \frac{(-1)^{k-1} (k+2\ell-1)!}
            {z_{\smash{i}} \lowerSup{k+2\ell}} \cdot
            \partial_{\smash{k}}^{(i)}
            \partial_{\smash{2 \ell}}^{(0)}
        \notag \\[-1ex]
        &
            {} + \frac{1}{2}
            \sum_{\substack{i; \ k, \ell \geq 0 \mathrlap{:} \\ k + \ell > 0}} {}
            \frac{(-1)^{k+\ell-1} (k+\ell-1)!}
            {(2z_i) \lowerSup{k+\ell}} \cdot
            \partial_{\smash{k}}^{(i)}
            \partial_{\smash{\ell}}^{(i)}
            \mp \frac{1}{2}
            \sum_{i; \ k > 0} {}
            \frac{(-1)^{k-1} (k-1)!}
            {z_{\smash{i}} \lowerSup{k}} \cdot
            \partial_{\smash{k}}^{(i)}
        \Biggr\}
        \notag \\*[-1.5ex]
        & \hspace{15em}
        \Bigl[
            a_1 (s_{\smash{k}}^{(1)}) \cdots a_n (s_{\smash{k}}^{(n)}) \cdot
            m (s_{\smash{2k}}^{(0)})
        \Bigr] \
        \bigg|_{\hspace{1.5em} \substack{
            \mathllap{s_{\smash{2k}}^{(i)}} \, \mapsto \, \mathrlap{2 s_{2k}} \\
            \mathllap{s_{\smash{2k+1}}^{(i)}} \, \mapsto \, \mathrlap{0} \\
            \mathllap{s_{\smash{2k}}^{(0)}} \, \mapsto \, \mathrlap{s_{2k}}
        } \hspace{1.2em}}
        \ ,
        \raisetag{3ex}
        \label{eq-osp-perf-expl}
    \end{align}
    where $a_i \in \mathrm{H}_\bullet (X; \mathbb{Q})$
    and $m \in \mathrm{H}_\bullet (Y; \mathbb{Q})$
    are homogeneous and supported on components
    $X_{r_i} \subset X$ and $Y_{r_0} \subset Y$,
    and `$\pm$' means `$+$' in type~B and `$-$' in type~C,
    and `$\mp$' vice versa, and
    \begin{equation*}
        \partial_{\smash{k}}^{(i)} =
        \begin{cases}
            r_i
            & \text{if $k = 0$,}
            \\
            \partial / \partial s_{\smash{k}}^{(i)}
            & \text{if $k > 0$,}
        \end{cases}
        \qquad
        \partial_{\smash{2k}}^{(0)} =
        \begin{cases}
            r_0
            & \text{if $k = 0$,}
            \\
            \partial / \partial s_{\smash{2k}}^{(0)}
            & \text{if $k > 0$,}
        \end{cases}
    \end{equation*}
    and~$i, j \in \{ 1, \dotsc, n \}$ throughout.
    The change of variables at the end of
    \cref{eq-osp-perf-expl}
    is the effect of pushing forward along
    $\diamond_{(n)} \colon X^n \times Y \to Y$.

    There is also the usual version given by \cref{thm-joyce-vm-tw},
    which is defined for all the types B, C, and D,
    with a similar expression to \cref{eq-osp-perf-expl}.
\end{para}

\begin{para}[Remark]
    \label{para-osp-remark}
    In~\cref{para-eg-perf-osp},
    we can also take~$Y$ to be the topological realization
    of classifying stacks of orthosymplectic perfect complexes.
    More precisely, define $\mathbb{Z}_2$-actions on
    $\calPerf_\mathrm{odd}$ or~$\calPerf_\mathrm{even}$
    by taking the dual complex.
    The data of a $\mathbb{Z}_2$-action
    includes a choice of identification
    $E^{\vee\vee} \simeq E$,
    which can be taken to be~$\pm 1$,
    which gives the distinction between types~C and~D.
    Define classifying stacks
    \begin{equation*}
        \calPerf_{\mathrm{B}} \ , \qquad
        \calPerf_{\mathrm{C}} \ , \qquad
        \calPerf_{\mathrm{D}}
    \end{equation*}
    of orthogonal (resp.~symplectic, orthogonal) perfect complexes
    of odd (resp.~even, even) rank,
    as derived fixed loci of these $\mathbb{Z}_2$-actions.
    We have an equivalence
    $|\calPerf_{\mathrm{B}}| \simto |\calPerf_{\mathrm{D}}|$
    given by $(-) \oplus \mathbb{C}$,
    with homotopy inverse
    $(-) \oplus \mathbb{C} [1] \oplus \mathbb{C} \oplus \mathbb{C} [-1]$.

    We expect that the topological realizations
    $| \calPerf_\mathrm{B} |$,
    $| \calPerf_\mathrm{C} |$, and~$| \calPerf_\mathrm{D} |$
    should agree with the spaces in~\cref{para-eg-perf-osp},
    although we cannot prove this yet.
    In any case,
    the constructions in~\cref{para-eg-perf-osp}
    still work when we replace~$Y$
    by these topological realizations.
\end{para}

\subsection{Principal bundles}

\begin{para}
    We describe the vertex induction
    for moduli stacks of principal $G$-bundles
    or perfect complexes on a variety.
    We consider the following versions:

    \begin{enumerate}
        \item
            An algebraic version,
            involving the moduli stack of
            principal $G$-bundles on a $\mathbb{C}$-variety.
        \item
            A topological version,
            using the topological mapping space
            from a $\mathbb{C}$-variety to~$\mathrm{B} G$.
        \item
            A version for moduli stacks of
            perfect complexes on a $\mathbb{C}$-variety,
            where we obtain the Joyce vertex algebra.
        \item
            A version for moduli stacks of
            orthosymplectic perfect complexes on a $\mathbb{C}$-variety,
            where we obtain twisted modules for the Joyce vertex algebra.
    \end{enumerate}
\end{para}

\begin{para}[The algebraic version]
    Let~$Z$ be a connected, smooth, projective variety over~$\mathbb{C}$,
    and let~$G$ be a linearly reductive algebraic group over~$\mathbb{C}$,
    with Lie algebra~$\mathfrak{g}$.
    Let
    \begin{equation}
        \mathcal{X} =
        \calBun_G (Z) =
        \calMap (Z, * / G)
    \end{equation}
    be the derived mapping stack,
    or the moduli stack of $G$-bundles on~$Z$.
    Its tangent complex~$\mathbb{T}_\mathcal{X}$
    satisfies that at a $\mathbb{C}$-point
    $[E] \in \mathcal{X} (\mathbb{C})$
    corresponding to a $G$-bundle~$E \to Z$, we have
    $\mathbb{T}_\mathcal{X} |_{[E]} \simeq
    \mathrm{H}^{\bullet + 1} (Z; \mathrm{Ad} (E))$,
    where~$\mathrm{Ad} (E) \to Z$ is the adjoint vector bundle of~$E$,
    with fibres isomorphic to~$\mathfrak{g}$.

    By \cref{thm-vi-stack}, there is a functor
    \begin{align*}
        \mathsf{Face}^\mathrm{sp} (\calBun_G (Z))^\mathrm{op}
        & \longrightarrow
        \mathsf{VS}_\mathbb{Q} \ ,
        \\
        \alpha
        & \longmapsto
        V_\alpha =
        \mathrm{H}_{\bullet + 2 \operatorname{vdim}} (\calBun_G (Z)_\alpha; \mathbb{Q}) \ .
    \end{align*}
    In particular, for each Levi subgroup
    $L \subset G$,
    meaning the centralizer of a cocharacter,
    there is a vertex induction map
    \begin{equation}
        \label{eq-vi-bun-g}
        \mathrm{H}_{\bullet + 2 \operatorname{vdim}} (\calBun_L (Z); \mathbb{Q})
        \longrightarrow
        \mathrm{H}_{\bullet + 2 \operatorname{vdim}} (\calBun_G (Z); \mathbb{Q})
        \llbr z_1, \dotsc, z_k \rrbr \,
        [ \lambda (z)^{-1} ] \ ,
    \end{equation}
    where $z_1, \dotsc, z_k$
    is a set of coordinates on
    $\Lambda_{\mathrm{Z} (L)^\circ}$,
    and we invert~$\lambda (z)$
    for non-zero elements
    $\lambda \in \Lambda^{\mathrm{Z} (L)^\circ}$
    that are images of roots of~$G$.
    The vertex induction map respects composition,
    in that for Levi subgroups $M \subset L \subset G$,
    the induction for $M \subset G$
    is the composition of the other two inductions in~$\mathsf{VS}_\mathbb{Q}$.

    If, moreover, $Z$ is Calabi--Yau of dimension $d \in 4 \mathbb{Z}$,
    then~$\mathcal{X}$ admits a $(2 - d)$-shifted symplectic structure,
    and we are in the situation of \cref{eg-real-va-symp}.
    For a Levi subgroup $L \subset G$,
    if we are given orientations of
    $\calBun_L (Z)$ and~$\calBun_G (Z)$,
    then there is a vertex induction map
    \begin{equation}
        \label{eq-vi-bun-g-real}
        \mathrm{H}_{\bullet + \operatorname{vdim}} (\calBun_L (Z); \mathbb{Q})
        \longrightarrow
        \mathrm{H}_{\bullet + \operatorname{vdim}} (\calBun_G (Z); \mathbb{Q})
        \llbr z_1, \dotsc, z_k \rrbr \,
        [ \lambda (z)^{-1} ] \ ,
    \end{equation}
    given by \cref{thm-real-joyce-va},
    with same notations as in~\cref{eq-vi-bun-g}.
\end{para}

\begin{para}[The topological version]
    Now let~$Z$ be a compact complex manifold,
    or more generally,
    a compact even-dimensional spin$^\mathrm{c}$ manifold,
    and let~$G$ be a linearly reductive algebraic group over~$\mathbb{C}$.
    Consider the topological mapping space
    \begin{equation}
        X =
        \mathrm{Bun}_G^{\smash{\mathrm{top}}} (Z) =
        \mathrm{Map} (Z, \mathrm{B} G) \ ,
    \end{equation}
    which is the space of topological $G$-bundles on~$Z$.
    Define the class
    \begin{equation}
        \mathbb{T}_X =
        -(\mathrm{pr}_X)_! \circ \mathrm{ev}^* ([\mathfrak{g}]) \ ,
    \end{equation}
    where $[\mathfrak{g}] \in K (\mathrm{B} G)$
    is the class of the adjoint representation,
    $\mathrm{ev} \colon Z \times X \to \mathrm{B} G$
    is the evaluation map,
    $\mathrm{pr}_X \colon Z \times X \to X$ is the projection, and
    $(\mathrm{pr}_X)_! \colon K (Z \times X) \to K (X)$
    is the Gysin map,
    defined in \textcite[\S IV.5.27]{karoubi-1978-k-theory}
    for spin$^\mathrm{c}$ manifolds,
    and defined here by taking the mapping space from~$X$
    to the Gysin map
    $\mathrm{Map} (Z, \mathrm{BU} \times \mathbb{Z}) \to \mathrm{BU} \times \mathbb{Z}$.

    Then for each Levi subgroup $L \subset G$,
    there is a vertex induction map
    \begin{equation}
        \label{eq-vi-bun-g-top}
        \mathrm{H}_{\bullet + 2 \operatorname{vdim}}
        (\mathrm{Bun}_L^{\smash{\mathrm{top}}} (Z); \mathbb{Q})
        \longrightarrow
        \mathrm{H}_{\bullet + 2 \operatorname{vdim}}
        (\mathrm{Bun}_G^{\smash{\mathrm{top}}} (Z); \mathbb{Q})
        \llbr z_1, \dotsc, z_k \rrbr \,
        [ \lambda (z)^{-1} ] \ ,
    \end{equation}
    similar to \cref{eq-vi-bun-g},
    given by \cref{thm-vi}.
    It respects composition in the same sense as in
    \cref{eq-vi-bun-g}.

    There is also a real version, similar to
    \cref{eq-vi-bun-g-real}.
    Let~$Z$ be a compact spin manifold
    of dimension~$8n$ for some $n \in \mathbb{N}$,
    such as a Calabi--Yau $4n$-fold,
    and let~$G$ be as above.
    Then using the class
    $[\mathfrak{g}] \in \mathit{KO} (\mathrm{B} G)$,
    and the Gysin map for such spin manifolds
    from \textcite[\S IV.5.27]{karoubi-1978-k-theory},
    we obtain a class
    $\mathbb{T}_X \in \mathit{KO} (X)$,
    where $X = \mathrm{Bun}_G^{\smash{\mathrm{top}}} (Z)$.
    If this class is orientable,
    then we obtain a real version of the vertex induction map
    \begin{equation}
        \label{eq-vi-bun-g-top-real}
        \mathrm{H}_{\bullet + \operatorname{vdim}}
        (\mathrm{Bun}_L^{\smash{\mathrm{top}}} (Z); \mathbb{Q})
        \longrightarrow
        \mathrm{H}_{\bullet + \operatorname{vdim}}
        (\mathrm{Bun}_G^{\smash{\mathrm{top}}} (Z); \mathbb{Q})
        \llbr z_1, \dotsc, z_n \rrbr \,
        [ \lambda (z)^{-1} ] \ ,
    \end{equation}
    similar to \cref{eq-vi-bun-g-real}.
\end{para}

\begin{para}[Perfect complexes]
    \label{para-perf-z}
    Let~$Z$ be a smooth proper $\mathbb{C}$-scheme,
    and let
    \begin{equation*}
        \mathcal{X} = \calPerf (Z)
        = \calMap (Z, \calPerf)
    \end{equation*}
    be the derived moduli stack of perfect complexes on~$Z$,
    as in \textcite[Definition~3.28]{toen-vaquie-2007-moduli},
    which admits a perfect tangent complex~$\mathbb{T}_\mathcal{X}$.
    Then

    \begin{itemize}
        \item
            \cref{thm-joyce-va} defines
            a vertex algebra structure on
            $\mathrm{H}_{\bullet + 2 \vdim} (\mathcal{X}; \mathbb{Q})$,
            originally due to \textcite{joyce-hall};
            see also \textcite[Theorem~4.4]{gross-homology}.

        \item
            \cref{thm-k-va} defines
            a multiplicative vertex algebra structure on
            $K_\circ (\mathcal{X}; \mathbb{Q})$,
            essentially due to
            \textcite{liu-k-wall-crossing}.

        \item
            If~$Z$ is Calabi--Yau of dimension $d \in 4 \mathbb{Z}$,
            then $\mathcal{X}$ is $(2 - d)$-shifted symplectic,
            and if we are given an orientation of~$\mathbb{T}_\mathcal{X}$,
            then \cref{thm-real-joyce-va} defines
            a vertex algebra structure on
            $\mathrm{H}_{\bullet + \vdim} (\mathcal{X}; \mathbb{Q})$,
            also originally due to \textcite{joyce-hall}.
    \end{itemize}

    Alternatively, we can consider a topological version of the above.
    Let~$Z$ now be a compact complex manifold, and let
    \begin{equation*}
        X = \mathrm{Map} (Z, \mathrm{BU} \times \mathbb{Z})
    \end{equation*}
    be the topological mapping space.
    Set $\mathbb{T}_X = - (\mathrm{pr}_X)_! (U^\vee \cdot U)$,
    where $U \in K (Z \times X)$ is classified by
    the evaluation map to $\mathrm{BU} \times \mathbb{Z}$,
    and $\mathrm{pr}_X \colon Z \times X \to X$ is the projection.
    Again,

    \begin{itemize}
        \item
            \cref{thm-joyce-va} defines
            a vertex algebra structure on
            $\mathrm{H}_{\bullet + 2 \vdim} (X; \mathbb{Q})$.

        \item
            \cref{thm-k-va} defines
            a multiplicative vertex algebra structure on
            $K_\circ (X; \mathbb{Q})$.

        \item
            If~$Z$ is Calabi--Yau of dimension $d \in 4 \mathbb{Z}$,
            then~$\mathbb{T}_X$ lifts to a class in $\mathit{KO} (X)$,
            and given an orientation,
            \cref{thm-real-joyce-va} defines
            a vertex algebra structure on
            $\mathrm{H}_{\bullet + \vdim} (X; \mathbb{Q})$.
    \end{itemize}

    By \textcite[Lemma~4.10 and Proposition~4.14]{gross-homology},
    if~$Z$ is of \emph{class D}
    in the sense of \cite[Definition~4.8]{gross-homology},
    which includes all curves and algebraic surfaces,
    then there is an equivalence $|\mathcal{X}| \simeq X$
    identifying $\mathbb{T}_\mathcal{X}$ and~$\mathbb{T}_X$,
    so the algebraic and topological versions agree.
\end{para}

\begin{para}[Orthosymplectic complexes]
    Let $\calPerf_{\mathrm{B/C/D}}$ be one of the three stacks
    in \cref{para-osp-remark}.

    Let~$Z$ be a smooth proper $\mathbb{C}$-scheme,
    and let~$\mathcal{X}$ be as in \cref{para-perf-z}.
    Let
    \begin{equation*}
        \mathcal{Y} = \calPerf_{\mathrm{B/C/D}} (Z)
        = \mathrm{Map} (Z, \calPerf_{\mathrm{B/C/D}})
    \end{equation*}
    be the derived mapping stack, using one of the three options,
    which admits a perfect tangent complex~$\mathbb{T}_\mathcal{Y}$.
    Let~$\mathcal{X}$ act on~$\mathcal{Y}$ by setting
    $x \diamond y = x \oplus y \oplus x^\vee$,
    similarly to \cref{para-eg-perf-osp}.
    Then
    \begin{itemize}
        \item
            \cref{thm-joyce-vm-tw}
            defines a twisted module structure on
            $\mathrm{H}_{\bullet + 2 \operatorname{vdim}} (\mathcal{Y}; \mathbb{Q})$
            for the vertex algebra
            $\mathrm{H}_{\bullet + 2 \operatorname{vdim}} (\mathcal{X}; \mathbb{Q})$.

        \item
            \cref{thm-k-va}
            defines a twisted module structure on
            $K_\circ (\mathcal{Y}; \mathbb{Q})$
            for the multiplicative vertex algebra
            $K_\circ (\mathcal{X}; \mathbb{Q})$.

        \item
            If~$Z$ is Calabi--Yau of dimension $d \in 4 \mathbb{Z}$,
            then $\mathcal{Y}$ is $(2 - d)$-shifted symplectic,
            and if we are given orientations of~$\mathcal{X}$ and~$\mathcal{Y}$,
            then \cref{thm-real-joyce-va}
            defines a twisted module structure on
            $\mathrm{H}_{\bullet + \vdim} (\mathcal{Y}; \mathbb{Q})$
            for the vertex algebra
            $\mathrm{H}_{\bullet + \vdim} (\mathcal{X}; \mathbb{Q})$.
    \end{itemize}

    Similarly, there are topological versions
    defined using topological mapping spaces to
    \begin{equation*}
        \mathrm{BO} \times (2 \mathbb{Z} + 1)
        \quad \text{or} \quad
        \mathrm{BSp} \times 2 \mathbb{Z}
        \quad \text{or} \quad
        \mathrm{BO} \times 2 \mathbb{Z} \ ,
    \end{equation*}
    and we avoid repeating the details here.
\end{para}

\subsection{Quivers}
\label{subsec-quiver}

\begin{para}[Quivers]
    \label{para-quiver}
    We very briefly describe
    the Joyce vertex algebra for quivers and its variants,
    which already exist in the literature.

    Let $Q = (Q_0, Q_1, s, t)$ be a quiver,
    where $Q_0$, $Q_1$ are finite sets
    of vertices and edges,
    and $s, t \colon Q_1 \to Q_0$
    are the source and target maps.

    Consider the moduli stack of representations of~$Q$,
    \begin{equation}
        \mathcal{X} =
        \coprod_{d \in \mathbb{N}^{Q_0}}
        V_d / G_d \ ,
    \end{equation}
    where $V_d = \prod_{a \in Q_1}
    \mathrm{Hom} (\mathbb{C}^{d (s (a))}, \mathbb{C}^{d (t (a))})$
    and $G_d = \prod_{i \in Q_0} \mathrm{GL} (\mathbb{C}^{d (i)})$.
    It is a smooth algebraic stack over~$\mathbb{C}$,
    and has a perfect tangent complex~$\mathbb{T}_\mathcal{X}$.
    Its topological realization is given by
    $|\mathcal{X}| \simeq \coprod_{d \in \mathbb{N}^{Q_0}} \mathrm{B} G_d$.

    We have a vertex algebra structure
    on $\mathrm{H}_{\bullet + 2 \dim} (\mathcal{X}; \mathbb{Q})$
    given by \cref{thm-joyce-va},
    due to \textcite{joyce-hall},
    and a multiplicative vertex algebra structure
    on $K_\circ (\mathcal{X}; \mathbb{Q})$
    given by \cref{thm-k-va},
    due to \textcite{liu-k-wall-crossing}.

    There is also a derived version,
    where we consider the derived moduli stack
    $\bar{\mathcal{X}} = \calPerf (Q)$
    of complexes of representations of~$Q$,
    as in \textcite[Definition~3.33]{toen-vaquie-2007-moduli}
    or \textcite[\S 4.3.2]{latyntsev-thesis}.
    Its topological realization is
    $|\bar{\mathcal{X}}| \simeq \prod_{i \in Q_0} {} (\mathrm{BU} \times \mathbb{Z})$.

    Again, we have Joyce's vertex algebra structure
    on $\mathrm{H}_{\bullet + 2 \vdim} (\bar{\mathcal{X}}; \mathbb{Q})$
    given by \cref{thm-joyce-va},
    and Liu's multiplicative vertex algebra structure
    on $K_\circ (\bar{\mathcal{X}}; \mathbb{Q})$
    given by \cref{thm-k-va}.

    By \textcite[\S 4.3.7 and Theorem~4.4.8]{latyntsev-thesis},
    the vertex algebra
    $\mathrm{H}_{\bullet + 2 \vdim} (\bar{\mathcal{X}}; \mathbb{Q})$
    is a lattice vertex algebra,
    and when~$Q$ is a Dynkin quiver of type A/D/E,
    it gives the simple quotient~$L_1 (\mathfrak{g})$
    of the Kac--Moody vertex algebra of the same type at level~$1$.
\end{para}

\begin{para}[Self-dual quivers]
    We now discuss orthosymplectic twisted modules
    for Joyce vertex algebras for quivers described above.

    Let~$Q$ be a quiver as above.
    A \emph{self-dual structure} on~$Q$ is the following data:
    \begin{enumerate}
        \item
            A \emph{contravariant involution}
            \begin{equation*}
                (-)^\vee \colon Q \longsimto Q^\mathrm{op} \ ,
            \end{equation*}
            where~$Q^\mathrm{op} = (Q_0, Q_1, t, s)$
            is the opposite quiver of~$Q$,
            such that $(-)^{\vee \vee} = \mathrm{id}$.

        \item
            Choices of signs
            \begin{equation*}
                u \colon Q_0 \longrightarrow \{ \pm 1 \} \ ,
                \qquad
                v \colon Q_1 \longrightarrow \{ \pm 1 \} \ ,
            \end{equation*}
            such that $u (i) = u (i^\vee)$ for all $i \in Q_0$,
            and $v (a) \, v (a^\vee) = u (s (a)) \, u (t (a))$
            for all $a \in Q_1$.
    \end{enumerate}
    For details, see \textcite[\S 6.1]{bu-osp-dt},
    \textcite{young-2015-self-dual,young-2016-hall-module,young-2020-quiver},
    and \textcite{derksen-weyman-2002}.

    As in \cite[\S 6.1.3]{bu-osp-dt},
    there is a moduli stack of \emph{self-dual representations} of~$Q$,
    \begin{equation}
        \mathcal{Y} = \coprod_{d \in (\mathbb{N}^{Q_0})^{\mathrm{sd}}}
        V_d^{\mathrm{sd}} / G_d^{\mathrm{sd}} \ ,
    \end{equation}
    where $(\mathbb{N}^{Q_0})^\mathrm{sd} \subset \mathbb{N}^{Q_0}$
    is the subset of dimension vectors~$d$ such that
    $d (i) = d (i^\vee)$ for all $i \in Q_0$,
    and~$d (i)$ is even if $i = i^\vee$ and $u (i) = -1$.
    The vector space~$V^\mathrm{sd}_d$
    and the group~$G^\mathrm{sd}_d$ are given by
    \begin{align}
        \label{eq-quiver-vsd}
        V^\mathrm{sd}_d & =
        \prod_{a \in Q_1^{\smash{\circ}} / \mathbb{Z}_2} {}
        \mathrm{Hom} (\mathbb{C}^{d (s (a))}, \mathbb{C}^{d (t (a))}) \times
        \prod_{a \in Q_1^{\smash{+}}} {}
        \mathrm{Sym}^2 (\mathbb{C}^{d (t (a))}) \times
        \prod_{a \in Q_1^{\smash{-}}} {}
        {\wedge}^2 (\mathbb{C}^{d (t (a))}) \ , 
        \\
        \label{eq-quiver-gsd}
        G^\mathrm{sd}_d & =
        \prod_{i \in Q_0^{\smash{\circ}} / \mathbb{Z}_2} {}
        \mathrm{GL} (\mathbb{C}^{d (i)}) \times
        \prod_{i \in Q_0^{\smash{+}}} {}
        \mathrm{O} (\mathbb{C}^{d (i)}) \times
        \prod_{i \in Q_0^{\smash{-}}} {}
        \mathrm{Sp} (\mathbb{C}^{d (i)}) \ ,
    \end{align}
    where~$Q_0^\circ$ is the set of vertices~$i$ with $i \neq i^\vee$,
    and~$Q_0^\pm$ the sets of vertices~$i$ with $i = i^\vee$ and $u (i) = \pm 1$.
    Similarly, $Q_1^\circ$ is the set of edges~$a$ with $a \neq a^\vee$,
    and $Q_1^\pm$ the sets of edges~$a$ with $a = a^\vee$ and
    $v (a) \, u (t (a)) = \pm 1$.
    There is a natural action
    ${\diamond} \colon \mathcal{X} \times \mathcal{Y} \to \mathcal{Y}$,
    given by
    $x \diamond y = x \oplus y \oplus x^\vee$,
    where $x^\vee$ is the dual representation of~$x$;
    see \cite[\S 6.1.2]{bu-osp-dt}.

    We have a twisted module
    $\mathrm{H}_{\bullet + 2 \dim} (\mathcal{Y}; \mathbb{Q})$
    for the Joyce vertex algebra
    $\mathrm{H}_{\bullet + 2 \dim} (\mathcal{X}; \mathbb{Q})$,
    given by \cref{thm-joyce-vm},
    and a twisted module
    $K_\circ (\mathcal{Y}; \mathbb{Q})$
    for the multiplicative vertex algebra
    $K_\circ (\mathcal{X}; \mathbb{Q})$,
    given by \cref{thm-k-va}.

    There is also a derived version,
    where we consider the derived stack
    $\bar{\mathcal{Y}} = \bar{\mathcal{X}}^{\mathbb{Z}_2}$,
    the fixed locus of the $\mathbb{Z}_2$-action on~$\mathcal{X}$
    given by the self-dual structure of~$Q$.
    Its topological realization is given by
    \begin{equation}
        \label{eq-quiver-bar-y}
        |\bar{\mathcal{Y}}| \simeq
        \prod_{i \in Q_0^\circ / \mathbb{Z}_2} {} (\mathrm{BU} \times \mathbb{Z})
        \times \prod_{i \in Q_0^+} {} |\calPerf_{\mathrm{O}}|
        \times \prod_{i \in Q_0^-} {} |\calPerf_{\mathrm{Sp}}|
        \ ,
    \end{equation}
    where $\calPerf_{\mathrm{O}} =
    \calPerf_\mathrm{B} \sqcup \calPerf_\mathrm{D}$
    and~$\calPerf_{\mathrm{Sp}} = \calPerf_\mathrm{C}$
    as in \cref{para-osp-remark}.

    We have a twisted module
    $\mathrm{H}_{\bullet + 2 \vdim} (\bar{\mathcal{Y}}; \mathbb{Q})$
    for the Joyce vertex algebra
    $\mathrm{H}_{\bullet + 2 \vdim} (\bar{\mathcal{X}}; \mathbb{Q})$,
    given by \cref{thm-joyce-vm},
    and a twisted module
    $K_\circ (\bar{\mathcal{Y}}; \mathbb{Q})$
    for the multiplicative vertex algebra
    $K_\circ (\bar{\mathcal{X}}; \mathbb{Q})$,
    given by \cref{thm-k-va}.

    As discussed in \cref{para-osp-remark},
    we expect to have
    $|\calPerf_\mathrm{O}| \simeq \mathrm{BO} \times \mathbb{Z}$
    and
    $|\calPerf_\mathrm{Sp}| \simeq \mathrm{BSp} \times 2 \mathbb{Z}$,
    but in any case,
    the above construction still works
    if we replace $|\bar{\mathcal{Y}}|$
    by the product~\cref{eq-quiver-bar-y}
    with $|\calPerf_\mathrm{O}|$,
    $|\calPerf_\mathrm{Sp}|$
    replaced by
    $\mathrm{BO} \times \mathbb{Z}$
    and $\mathrm{BSp} \times 2 \mathbb{Z}$.
\end{para}

\phantomsection
\addcontentsline{toc}{section}{References}
\sloppy
\setstretch{1.1}
\renewcommand*{\bibfont}{\normalfont\small}
\printbibliography

\par\noindent\rule{0.38\textwidth}{0.4pt}
{\par\noindent\small
\hspace*{2em}Chenjing Bu\qquad
\texttt{bucj@mailbox.org}
\\[-2pt]
\hspace*{2em}Mathematical Institute, University of Oxford, Oxford OX2 6GG, United Kingdom.}

\end{document}